\documentclass[aap,MSNbibl,dvips]{arximspdf}
\usepackage{graphicx}

\doi{10.1214/08-AAP591}
\volume{19}
\issue{5}
\pubyear{2009}
\firstpage{1719}
\lastpage{1780}

\makeatletter

\newproclaim{defn}{Definition}[section]
\newproclaim{rem}{Remark}[section]
\newtheorem{lemma}{Lemma}[section]
\newtheorem{corollary}{Corollary}[section]
\newtheorem{theorem}{Theorem}[section]
\newtheorem{proposition}{Proposition}[section]
\newtheorem{assumption}{Assumption}[section]
\newtheorem{conjecture}{Conjecture}[section]
\newproclaim{example}{Example}[section]

\makeatother

\begin{document}
\begin{frontmatter}

\title{State space collapse and diffusion approximation for a
network operating under a fair bandwidth sharing policy}
\runtitle{State space collapse and diffusion approximation}

\begin{aug}
\author[A]{\fnms{W. N.} \snm{Kang}\ead[label=e1]{weikang@andrew.cmu.edu}},
\author[B]{\fnms{F. P.} \snm{Kelly}\ead[label=e2]{F.P.Kelly@statslab.cam.ac.uk}\thanksref{T1}},
\author[C]{\fnms{N. H.} \snm{Lee}\ead[label=e3]{nhlee@jhu.edu}}
\and
\author[D]{\fnms{R. J.} \snm{Williams}\corref{}\ead[label=e4]{williams@math.ucsd.edu}\thanksref{T2}}
\thankstext{T1}{Supported in part by EPSRC Grant GR/586266/01.}
\thankstext{T2}{Supported in part by NSF Grant DMS-06-04537.}
\runauthor{Kang, Kelly, Lee and Williams}
\affiliation{Carnegie Mellon University, University
of Cambridge, Johns Hopkins University and University of California,
San Diego}
\address[A]{W. N. Kang\\ Department of Mathematical Sciences\\
Carnegie Mellon University \\
Pittsburgh, Pennsylvania 15213-3890\\
USA\\ \printead{e1}}

\address[B]{F. P. Kelly\\ Statistical Laboratory\\
Centre for Mathematical Sciences\\
University of Cambridge\\
Wilberforce Road\\
Cambridge CB3 0WB\\ United Kingdom\\ \printead{e2}}

\address[C]{N. H. Lee\\ Department of Applied Mathematics\\
\quad and
Statistics\\
Johns Hopkins
University\\
3400 North Charles Street\\
Baltimore, Maryland 21218-2682\\USA\\ \printead{e3}}

\address[D]{R. J. Williams\\ Department of Mathematics\\
University of California, San Diego,\\
9500 Gilman Drive\\
La Jolla, California 92093-0112\\ USA\\ \printead{e4}}
\end{aug}

\received{\smonth{11} \syear{2007}}
\revised{\smonth{8} \syear{2008}}

%
\begin{abstract}
We consider a connection-level model of Internet congestion control,
introduced by Massouli\'{e} and Roberts [\textit{Telecommunication
Systems} \textbf{15} (2000) 185--201], that represents
the randomly varying number of flows present in a network. Here,
bandwidth is shared fairly among elastic document transfers
according to a weighted $\alpha$-fair~band\-width sharing policy
introduced by Mo and Walrand [\textit{IEEE/ACM
Transactions on Networking} \textbf{8} (2000) 556--567] [$\alpha\in
(0,\infty)$].
Assuming Poisson arrivals and exponentially distributed document
sizes, we focus on the heavy traffic regime in which the average
load placed on each resource is approximately equal to its capacity.
A fluid model (or functional law of large numbers approximation) for
this stochastic model was derived and analyzed in a prior work
[\textit{Ann. Appl. Probab.} \textbf{14} (2004) 1055--1083] by two of
the authors. Here, we use the long-time behavior
of the solutions of the fluid model established in that paper to
derive a property called \emph{multiplicative state space collapse},
which, loosely speaking,
shows that in diffusion scale, the flow count process for the
stochastic model can be approximately recovered as a continuous
lifting of the workload process.

Under weighted proportional fair
sharing of bandwidth ($\alpha=1$) and a mild local traffic
condition, we show how multiplicative
state space collapse can be combined with
uniqueness in law
and an
invariance principle for the diffusion
[\textit{Theory Probab. Appl.} \textbf{40} (1995) 1--40], [\textit
{Ann. Appl. Probab.}
\textbf{17} (2007) 741--779] to establish a diffusion approximation for
the workload process and hence to yield an approximation for the
flow count process. In this case, the workload diffusion behaves
like Brownian motion in the interior of a polyhedral cone and is
confined to the cone by reflection at the boundary, where the
direction of reflection is constant on any given boundary face. When
all of the weights are equal (proportional fair sharing), this diffusion
has a product form invariant measure. If the latter is integrable,
it yields the unique
stationary distribution for the diffusion which has a strikingly simple
interpretation in terms of independent dual random variables, one for
each of the resources of the network.
We are able to extend this product form
result to the case where document sizes are distributed as
finite mixtures of exponentials and to models that include
multi-path routing.
We indicate some
difficulties related to extending the diffusion approximation result to
values of $\alpha\neq1$.

We illustrate our
approximation results for a few simple networks.
In particular, for a two-resource linear network, the
diffusion lives in a wedge that is a strict subset of the positive
quadrant. This geometrically illustrates the entrainment of
resources, whereby congestion at one resource may prevent another
resource from working at full capacity. For a four-resource network
with multi-path routing, the product form result under
proportional fair sharing is expressed
in terms of independent dual random variables, one for each of a set
of generalized cut constraints.
\end{abstract}

%
\begin{keyword}[class=AMS]
\kwd{60K30}
\kwd{90B15}.
\end{keyword}
\begin{keyword}
\kwd{$\alpha$-fair} \kwd{bandwidth sharing} \kwd{Brownian model}
\kwd{diffusion approximation} \kwd{flow-level Internet
congestion control} \kwd{fluid model} \kwd{invariant manifold}
\kwd{multi-path routing} \kwd{product form stationary distribution}
\kwd{proportional fair sharing} \kwd{reflected Brownian motion}
\kwd{simultaneous resource possession} \kwd{state
space collapse} \kwd{workload}.
\end{keyword}

\end{frontmatter}

\section{Introduction}

We consider a connection-level model of Internet congestion control
introduced and studied by Massouli\'{e} and Roberts \cite{RM}. This
stochastic model represents the randomly varying number of flows
present in a network where bandwidth is dynamically shared between
flows that correspond to continuous transfers of individual elastic
documents. This model, which we shall refer to as a \emph{flow-level
model}, assumes a ``separation of time scales'' such that the time
scale of the flow dynamics (i.e., of document arrivals and
departures) is much longer than the time scale of the packet level
dynamics on which rate control schemes such as the Transmission
Control Protocol (TCP) converge to
equilibrium. We consider the flow-level model operating under a
family of bandwidth sharing policies introduced by Mo and Walrand
\cite{MW}, called \textit{weighted $\alpha$-fair} policies.
Here,
$\alpha$ is a parameter lying in $ (0,\infty)$. The case $\alpha=1$
is called \textit{weighted proportional fair sharing} and the case
$\alpha
\to\infty$ corresponds to what is called \textit{weighted max--min fair}.

Assuming Poisson arrivals and exponentially distributed document
sizes, de Veciana, Lee and Konstantopoulos \cite{DLK} and Bonald and
Massouli\'{e} \cite{BM} studied stability of the flow-level model
operating under weighted $\alpha$-fair bandwidth sharing policies
(including limiting values of $\alpha$). Lyapunov functions
constructed in \cite{DLK} for weighted max--min fair and
proportionally fair policies, and in \cite{BM} for weighted
$\alpha$-fair policies [$\alpha\in(0,\infty)$],
imply positive recurrence of the Markov chain associated with the
model when the average load on each resource is less than its
capacity. Lin, Shroff and Srikant
\cite{LinSch2004,LinSchSri2006,Sri2004} have recently given sufficient
conditions for
stability of a Markov model under a back-pressure algorithm when the
assumption of time scale
separation is relaxed. For more general document size distributions,
there are a few results for specific values of $\alpha$ or for specific
distributions or topologies that provide sufficient conditions for
stability of the flow-level model operating under bandwidth sharing
policies \cite{Bra2005,CST2006,LaBeSr2004,Mas2005}.
A summary of
these results is provided in the introduction to \cite{growil06}. In
general, it remains an open question whether, with renewal arrivals
and arbitrarily (rather than exponentially) distributed document
sizes, the flow-level model is stable under an $\alpha$-fair
bandwidth sharing policy [$\alpha\in(0,\infty)$] when the nominal
load placed on each resource is less than its capacity
(see \cite{growil06,growil07} for some first steps in this direction).
Here, we
restrict our attention to the case of Poisson arrivals and exponential document
sizes, for which stability is well understood.

We are interested in using diffusion approximations to explore the
performance of the flow-level model operating under a weighted
$\alpha$-fair bandwidth sharing policy when the average load placed
on each resource is approximately equal to its capacity, that is, the
system is heavily loaded. We are particularly interested in
manifestations of the phenomenon of entrainment, whereby congestion
at some resources may prevent other resources from working at their
full capacity.

There are several motivations for our work. One
source of motivation lies in fixed-point approximations of network
performance for TCP networks (see \cite{BUT,GSV,REV}). These
approximations require, as input, information on the joint
distribution of the numbers of flows present on different routes,
where dependencies between these numbers may be induced by the
bandwidth sharing mechanism. Similarly, an understanding of such
joint distributions seems important if the performance models for a
single bottleneck described by Ben Fredj et al.~\cite{BF} are
to be generalized to a network.

Another motivation is that the
flow-level model typically involves the simultaneous use of several
resources. Due to the exponential document sizes, this model can be
equated (in distribution) with a stochastic processing network (SPN)
as introduced by Harrison \cite{HAR}. Open multiclass queueing
networks operating under head-of-the-line (HL) service disciplines
are a special case of SPNs without simultaneous resource possession.
For certain queueing networks of this type, it has been shown
\cite{BR,Wi} that suitable asymptotic behavior of critical fluid
models implies a property called \textit{multiplicative state space
collapse}, which, in turn, validates the use of Brownian model
approximations for these networks in heavy traffic. For more general
SPNs, investigation of the behavior of critical fluid models, of a
related notion of multiplicative state space collapse and of the
implications for
diffusion approximations are in the early stages of development.
The analysis in this paper can be viewed as a contribution to such
an investigation for models involving simultaneous resource
possession.
For another contribution, see the paper of Ye and Yao \cite{YY},
who consider a stochastic processing network with simultaneous resource
possession; in contrast to the fully heavily loaded,
multiple bottleneck situation
considered here, Ye and Yao consider the situation of a single heavily
loaded bottleneck.
A further recent contribution is the important paper of
Shah and Wischik \cite{SW}, who have proven multiplicative
state space collapse for a class of ``switched'' networks with multiple
bottlenecks operating under a family of scheduling policies
related to the maximum weight algorithm introduced by
Tassiulas and Ephremides~\cite{TE}.

Finally, although we restrict to exponential document
sizes in this paper, we would like to relax that assumption in
future work. Although this involves a significantly more elaborate
stochastic model (see \cite{growil06}) to keep track of residual
document sizes (because of the processor sharing nature of the
bandwidth sharing policy), knowing the results for exponential
document sizes is likely to be useful for such work.

\subsection{Overview}\label{sec1.1}

In this paper, we consider the flow-level model with Poisson
arrivals and exponentially distributed document sizes operating
under a weighted $\alpha$-fair bandwidth sharing policy for $\alpha
\in(0,\infty)$. We focus on the heavy traffic regime in which the
average load placed on each resource is approximately equal to its
capacity. We recall the definition of a critical fluid model from
the prior work \cite{KW} of two of the authors; this model is a
formal functional law of large numbers approximation to the
flow-level model. The asymptotic behavior of this critical fluid
model was studied in \cite{KW}. Here, we show how this behavior can
be used to prove a property called \textit{multiplicative state space
collapse}. Loosely speaking, this says that in diffusion scale, the
flow count process can be approximately recovered by a continuous
lifting of the lower-dimensional workload process.
Given the
asymptotic behavior of the critical fluid model,
our proof of multiplicative state
space collapse
follows a general line of argument
pioneered by Bramson
in \cite{BR}, where open multiclass queueing networks operating under
certain head-of-the-line (HL) service disciplines are treated.
There are some differences in setup and
proof details between our treatment and Bramson's \cite{BR}. These
are described in detail at the beginning of Section \ref{section:MSSC}.
However, we
wish to emphasize that
our main line of argument follows that of Bramson \cite{BR}.
It is interesting to note that, in
contrast to prior results on state space collapse for open
multiclass queueing networks, our
lifting map can be nonlinear (for $\alpha\not= 1$).

The multiplicative state space collapse result leads to a natural
conjecture for a diffusion approximation to the workload process.
For the case of weighted proportional fair sharing of
bandwidth ($\alpha=1$), we combine
multiplicative state space collapse
with uniqueness in law
for the diffusion \cite{DW} and an \textit{invariance principle} \cite
{KaWi} for
semimartingale reflecting Brownian motions living in domains with
piecewise smooth boundaries to obtain a diffusion approximation for
the flow count process under a mild local traffic condition. This
diffusion lives in a polyhedral cone. It behaves like Brownian
motion in the interior of the cone and is confined to the cone by
reflection (or regulation) at the boundary where the direction of
reflection is constant on any given boundary face.
We illustrate this diffusion approximation
result for a~simple two-resource linear network. Then, the diffusion
lives in a wedge that is a strict subset of the positive quadrant.
This geometrically illustrates the entrainment of resources, whereby
congestion at one resource may prevent another resource from working
at full capacity. We also observe how the wedge can
vary with the weights. Ongoing work is directed toward establishing
diffusion approximations for the workload process when
$\alpha\not=1$. We mention some of the difficulties associated with
this.
These center around the fact that when
$\alpha\not=1$ and the workload dimension is three or higher,
although the state space for the
proposed diffusion approximation for the workload process
is a cone, it is not
a polyhedral cone. Indeed, the cone has curved boundary faces that
intersect nonsmoothly and can even
meet in cusp-like singularities.
The current lack of a general existence and uniqueness theory (and
an associated invariance principle) for
reflecting Brownian motions in such domains
is a major obstacle to proving the conjecture for $\alpha\not=1$.

In the case of proportional fair sharing, that is,
when $\alpha=1$ and all of the
weights for the bandwidth sharing policy are equal,
we show that the approximating diffusion
has a product form
invariant measure.
When the latter is integrable over the state space,
our results suggest a strikingly simple approximation for
the joint stationary distribution of the number of flows present on
different routes under proportional fair sharing and the mild local
traffic condition.
In this, each of
the resources of the network has associated with it a
dual random variable. These dual variables are independent and
exponentially distributed, and the formal approximation to the number
of flows
on a route is proportional to the sum
of the dual variables along the route.

We also indicate an extension of the product form result
to the situation where document
sizes are finite mixtures of exponential distributions.
Under this extension, the formal approximation for the
joint stationary distribution for the number of flows
present on different routes is insensitive:
that is, the approximation does not depend on the distributions
of document sizes, other than through the means of these distributions,
provided that the distributions are finite mixtures of exponentials.
This result complements the known result~\cite{BM,BP,RM} that,
for proportional fair sharing
and a small class of topologies and parameters, the stationary
distribution for the number of flows present on different routes is
exactly insensitive: that is, the stationary distribution
does not depend on the distributions
of document sizes, other than through
the means of these distributions.

Finally, we indicate a relation to more general models with
routing.
There is considerable interest in multi-path routing in the Internet
and rate control schemes generalizing TCP have been proposed \cite{HSHST,KV}.
 It is known that the stability region for the flow-level model may be
strictly increased if multi-path routing is allowed \cite{HSHST,KM}.
We show that our results on diffusion approximations
under proportional fair sharing
extend to the multi-path case. The local traffic condition
becomes more difficult to verify in this setting,
but if it is satisfied, then our results suggest
a simple approximation for the stationary distribution of the numbers
of source--destination flows in terms of
independently distributed dual variables, one for each
generalized cut constraint.

A summary of some of the results of this paper (without proofs) was given
in the conference proceedings papers \cite{KKLW}
and \cite{KKLWM}.

\subsection{Notation and terminology}\label{sec1.2}

For each positive integer $d\geq1$, $\mathbb{R}^d$ will denote
$d$-dimensional Euclidean space and the positive orthant in this
space will be denoted by $\mathbb{R}^d_+ = \{x \in\mathbb{R}^d \dvtx
x_i\geq0
\mbox{ for } i = 1,\ldots, d\}.$ When $d=1$, we sometimes write
$\mathbb{R}$ instead of $\mathbb{R}^1$, and $\mathbb{R}_+$ instead
of $\mathbb{R}_+^1$. The
Euclidean norm of $x\in\mathbb{R}^d$ will be denoted by $|x|$.
Vectors in
$\mathbb{R}^d$ will be assumed to be column vectors unless specifically
indicated otherwise. The transpose of a vector or matrix will be
denoted by the use of a superscript ``$'$''.
Inequalities between vectors in
$\mathbb{R}^d$ will be interpreted componentwise. That is, for $x,
y\in
\mathbb{R}^d, x\leq y$ is equivalent to $x_i\leq y_i $ for $ i =
1,\ldots, d.$ For each $x\in\mathbb{R}^d$ and each set $S\subset
\mathbb{R}^d$,
the distance between $x$ and $S$ is denoted by
\[
d(x,S)=\inf\{|x-y|\dvtx y\in S\}.
\]
For $x, y\in\mathbb{R}$, $x\vee y=\max\{x,
y\}$. For each $x\in\mathbb{R}$, $\lfloor x \rfloor$ denotes the largest
integer less than or equal to $x$. Given a vector $x\in\mathbb{R}^d$, the
$d\times d$ diagonal matrix with the entries of $x$ on its diagonal
will be denoted by diag($x$). For positive integers $d_1$ and $d_2$,
the norm of a $d_1\times d_2$ matrix $A$ will be given by
\[
\Vert A\Vert=\Biggl(\sum_{i=1}^{d_1} \sum_{j=1}^{d_2}
A_{ij}^2\Biggr)^{1/2}.
\]
The set of nonnegative integers will be
denoted by $\mathbb{Z}_+$ and the set of points in $\mathbb{R}^d_+$
with all
integer coordinates will be denoted by $\mathbb{Z}^d_+$. A sum over an
empty set of indices will be taken to have a value of zero. The
cardinality of a finite set $S$ will be denoted by $|S|$.
For $0\leq s< t< \infty$, any integer $d\geq1$ and any bounded
function $x\dvtx[s,t]\to
\mathbb{R}^d$, let $\| x( \cdot )\|_{[s,t]}=\sup_{u\in[s,t]}
|x(u)|$ and when $s=0$, let $\|x(  \cdot )\| _t
=\|x( \cdot )\|_{[0,t]}$.

All stochastic processes in this paper will be assumed to have sample
paths that are right-continuous with finite left limits (r.c.l.l.).
We denote by $\mathbb{D}([0,\infty),\mathbb{R}^d)$ the space of
r.c.l.l.~functions
from $[0,\infty)$ into $\mathbb{R}^d$ and we endow this space
with the usual Skorokhod $J_1$-topology.
We denote by $\mathbb{C}([0,\infty),\mathbb{R}^d)$ the space
of continuous functions from $[0,\infty)$ into $\mathbb{R}^d$.
The Borel $\sigma$-algebra on either
$\mathbb{D}([0,\infty), \mathbb{R}^d)$ or $\mathbb{C}([0,\infty),
\mathbb{R}^d)$ will be
denoted by~$\mathcal B^d$.
Consider $ X, X^1,
X^2,\ldots,$ each of which is a $d$-dimensional process (possibly
defined on different probability spaces). The sequence
$\{X^n\}_{n=1}^{\infty}$ is said to be \textit{tight} if the
probability measures induced by the $X^n$ on the measurable space
$(\mathbb{D}([0,\infty), \mathbb{R}^d),{\mathcal{B}}^d)$ form a
tight sequence,
that is, they form a~weakly relatively compact sequence in the space of probability
measures on $(\mathbb{D}([0,\infty), \mathbb{R}^d),{\mathcal{B}}^d)$.
The notation ``$X^n\Rightarrow
X$'' will mean that as $n\rightarrow\infty$, the sequence of
probability measures induced on $(\mathbb{D}([0,\infty), \mathbb
{R}^d), \mathcal B^d)$
by $\{X^n\}$ converges weakly to
the probability measure induced on the same space
by~$X$. We shall describe this in words
by saying that $X^n$ converges weakly (or in distribution)
to~$X$ as $n\to\infty$. The sequence of processes
$\{X^n\}_{n=1}^{\infty}$ is called \textit{C-tight} if it is tight
and if each weak limit point, obtained as a weak limit along a
subsequence, almost surely has sample paths in $\mathbb{C}([0,\infty
),\mathbb{R}^d)$.

\section{Flow-level model}\label{secflow}

\subsection{Network structure}

We consider a network with finitely many \textit{resources}  labeled
by $j\in{\mathbb J}\not=\varnothing$. A \textit{route} $i$ is a nonempty
subset of ${\mathbb J}$ (interpreted as the set of resources used by route
$i$). We are given a finite, nonempty set ${\mathbb I}$ of allowed routes.
Let ${\mathbf J}=|{\mathbb J}|$, the total number of resources, and
${\mathbf I}=|{\mathbb I}|$,
the total number of routes. Let $A$ be the ${\mathbf J}\times{\mathbf
I}$ \textit{incidence matrix} which contains only zeros and ones, defined such
that $A_{ji}=1$ if resource $j$ is used by route~$i$
and $A_{ji}=0$ otherwise.
We assume that $A$ has rank ${\mathbf J}$ so that it has full row rank. We
further assume that resource (bandwidth) \textit{capacities} $(C_j:
j \in{\mathbb J})$ are given and that these are all strictly positive and
finite.

\subsection{Stochastic primitives}
An active flow on route $i$ corresponds to the
continuous transmission of a document through the resources used by
route~$i$. Transmission is assumed to occur simultaneously through
all resources on route $i$. It is assumed that a new document
arrives to route $i$ at each jump time of a Poisson process that
has rate parameter $\nu_i>0$ and that each such document has an
exponentially distributed size with mean $1/\mu_i$, where $\mu_i\in
(0,\infty)$. These document sizes are assumed to be independent of
one another and to be independent of all arrival times of documents.
The number of documents on route $i$ at time zero is assumed to be
independent of the remaining sizes of those documents and these
sizes are assumed to be independent and exponentially distributed
with mean $1/\mu_i$.
Initial numbers and sizes of documents, arrival times of new
documents and their sizes for different routes $i\in{\mathbb I}$ are
assumed to be mutually independent.

\subsection{Bandwidth sharing policy}
Bandwidth capacity is allocated dynamically to the routes according
to the following bandwidth sharing policy which was first introduced\vadjust{\goodbreak}
by Mo and Walrand \cite{MW}. The bandwidth for a route is shared
equally among all of the documents currently being transmitted
over that route. Given a fixed parameter $\alpha\in(0, \infty)$
and strictly positive weights $({\kappa}_i\dvtx i\in{\mathbb I})$, if $N_i(t)$
denotes the (random) number of flows on route $i$ at time~$t$ for
each $i\in{\mathbb I}$ and $N(t) =(N_i(t)\dvtx i\in{\mathbb I})$, then the
bandwidth
allocated to route $i$ at time
$t$ is
given by $\Lambda_i (N(t))$ and this bandwidth is shared equally
among all of the flows on route $i$, where the function
$\Lambda( \cdot )=(\Lambda_i( \cdot )\dvtx i\in{\mathbb I})$ is
defined as
follows (we define it on all of $\mathbb{R}_+^{{\mathbf I}}$ as we
shall later
apply it to rescaled versions of $N$).

Let $\Lambda\dvtx \mathbb{R}^{{\mathbf I}}_+\to\mathbb{R}^{{\mathbf
I}}_+$ be defined such that for
each $n\in\mathbb{R}^{{\mathbf I}}_+$, $\Lambda_i(n)
=0$ for $i\in{\mathbb I}_0(n) \equiv\{l \in{\mathbb I}\dvtx n_l =0\}$,
and when $ {\mathbb I}_+(n) \equiv\{l\in{\mathbb I}\dvtx n_l >0\}$ is nonempty,
$\Lambda^+(n) \equiv
(\Lambda_i(n) \dvtx i\in{\mathbb I}_+(n))$ is the unique value of $\Lambda^+
=(\Lambda_i\dvtx i\in{\mathbb I}_+(n))$ that solves the optimization problem
%
\begin{eqnarray}\label{eq:opt}
&&\mbox{maximize}\quad G_{n}(\Lambda^+)\nonumber
\\
&&\mbox{subject to}\quad  \sum_{i\in{\mathbb I}_+(n)} A_{ji}
\Lambda_i \leq C_j, \qquad j \in{\mathbb J},
\\
&&\mbox{over} \quad\hspace{24pt} \Lambda_i \geq0, \qquad i\in{\mathbb I}_+(n),\nonumber
\end{eqnarray}
where for $n \in\mathbb{R}_+^{{\mathbf I}}\setminus\{0\}$ and
$\Lambda^+
=(\Lambda_i\dvtx i\in{\mathbb I}_+(n)) \in\mathbb{R}_+^{|{\mathbb I}_+(n)|}$,
%
\begin{equation} \label{req:G}
G_{n} (\Lambda^+) =
\cases{\displaystyle
\sum\limits_{i\in{\mathbb I}_+(n)} {\kappa}_i n_i^\alpha
\frac{\Lambda_i^{1-\alpha}}{1-\alpha},
 &\quad \mbox{if} $\alpha\not= 1,$ \cr
\displaystyle \sum\limits_{i \in{\mathbb I}_+(n)} {\kappa}_i n_i\log\Lambda_i,
 &\quad \mbox{if } $\alpha=1$
 }
\end{equation}
and the value of the right member above is taken to be $-\infty$
if $\alpha\in[1,\infty)$ and $ \Lambda_i = 0$ for some $i\in
{\mathbb I}_+(n)$. The resulting bandwidth allocation is called a
\textit{weighted $\alpha$-fair allocation}.

The properties of the function $\Lambda$ are summarized in the
following proposition. This proposition is proved in the Appendix of
Kelly and Williams~\cite{KW}.

\begin{proposition} \label{r:Lambdapro}
For each $n\in\mathbb{R}_+^{{\mathbf I}}$,
\begin{longlist}
\item[(i)] $\Lambda_i(n)>0$ for each $i\in{\mathbb I}_+(n)$;

\item[(ii)] $\Lambda(r n) = \Lambda(n) $ for each $r>0$;
\item[(iii)] $\Lambda_i(\cdot)$ is
continuous at $n$ for those $i$ such that $n_i>0$;
\item[(iv)] there exists at least one $p\in\mathbb{R}_+^{{\mathbf
J}}$, depending on
$n$, such that
\end{longlist}
%
\begin{equation}
\label{eq:lm1}
\Lambda_i (n) = n_i \biggl( \frac{\kappa_i}{ \sum_{j\in
{\mathbb J}}
p_j A_{ji} } \biggr) ^{1/ \alpha} \qquad \mbox{for all } i\in
{\mathbb I}_+(n),
\end{equation}
where
%
\begin{equation}
\label{eq:lm2}
p_j\biggl(C_j - \sum_{i\in{\mathbb I}} A_{ji} \Lambda_i(n)\biggr) =0
\qquad\mbox{for all $j\in{\mathbb J}$}.
\end{equation}
\end{proposition}

The $(p_j \dvtx j\in{\mathbb J})$ are Lagrange multipliers (or {\it dual
variables})
for the optimization
problem, where there is one multiplier for each of the capacity
constraints.

\subsection{Stochastic process description}

$\!\!$The flow count process \mbox{$ N\,{=}\,(N_i\dvtx i\,{\in}\,{\mathbb{I}})$} is a Markov
process with state space $\mathbb{Z}_+^{{\mathbf I}}$. We use the following
(equivalent in distribution) representation for $N$ and the
cumulative unused capacity process $U=(U_j\dvtx j\in{\mathbb J})$:
%
\begin{eqnarray}\label{N}
N_i(t)& = & N_i(0) + E_i(t) -S_i(T_i(t)) , \qquad i\in{\mathbb I},
\\ \label{U}
U_j(t) & = & C_j t - \sum_{i\in{{\mathbb I}}} A_{ji} T_i(t), \qquad
j\in{\mathbb J},
\end{eqnarray}
where $E _i$ is a Poisson process with rate $\nu_i$, $S _i$ is a
Poisson process with ra\-te~$\mu_i$, $T _i(t)$ is the cumulative
amount of bandwidth allocated to route $i$ up to time $t$ and
%
\begin{equation} \label{T}
T_i(t) = \int_0^t \Lambda_i (N(s)) \,  ds.
\end{equation}
We assume that for each $i\in\mathbb{I}$, $E_i$ and $S_i$ are represented
by
%
\begin{equation}
E_i(t)=\operatorname{sup}\Biggl\{n\geq0\dvtx \sum_{l=1}^n\xi
_i(l)\leq
t\Biggr\}
\end{equation}
and
%
\begin{equation}
S_i(t)=\operatorname{sup}\Biggl\{n\geq0\dvtx \sum_{l=1}^n\zeta_i(l)\leq
t\Biggr\},
\end{equation}
respectively, where $\{\xi_i(l)\}_{l=1}^{\infty}$ is a sequence of
i.i.d.~exponential random variables with mean $1/\nu_i$ and
$\{\zeta_i(l)\}_{l=1}^{\infty}$ is a sequence of i.i.d.~exponential
random variables with mean $1/\mu_i$. It is assumed that
$\{\xi_i(l)\}_{l=1}^\infty,$ $\{\zeta_i(l)\}_{l=1}^\infty$,
$N_i(0),$ for $ i\in{\mathbb I}$, are mutually independent. We define an
(average) \textit{workload process}  by
%
\begin{equation}\label{req:W}
W(t ) = A M^{-1} N(t ) \qquad\mbox{for all } t\geq0,
\end{equation}
where $M=\operatorname{diag}(\mu)$ is the ${\mathbf I}\times
{\mathbf I}$ diagonal matrix with
the entries of $\mu$ on its diagonal.

\section{Sequence of systems and scaling}\label{r:sequence}

Consider an increasing sequence of positive scale parameters
$\{r_l\}_{l=1}^\infty$ which converges to infinity. To ease the
notation, we shall simply write $r$ in place of $r_l$, where it is
understood that~$r$ increases to infinity through a sequence. We
consider a sequence of flow-level models indexed by~$r$, where the
network structure with parameters~$A$ and $C$, and bandwidth sharing
policy with parameters~$\alpha$ and $\{\kappa_i, i\in\mathbb{I}\}$
do not
vary with $r$. Each member of the sequence is a stochastic system, as
described in the previous section. We append a superscript of $r$ to
any process, sequence of random variables or parameter associated
with the $r$th system that depends on $r$. Thus, we have
processes $N^r, W^r, U^r, T^r, E^r, S^r$, sequences of random
variables $\xi^r_i = \{\xi^r_i(l)\}_{l=1}^\infty$ and $\zeta^r_i =
\{\zeta^r_i(l)\}_{l=1}^\infty$ for $i\in{\mathbb I}$, parameters
$\nu^r$
and $ \mu^r$, and matrices $M^r$. Let $\rho_i^r =\nu_i^r/\mu_i^r$
for each $i\in{\mathbb I}$. We shall henceforth assume that the following
\textit{heavy traffic} condition holds.

\begin{assumption}[(Heavy traffic)] \label{HT}
There exist $\nu, \mu\in\mathbb{R}_+^{\bf I}$ and $\theta\in
\mathbb{R}^{\bf J}$
such that $\nu_i>0$ and $\mu_i>0$ for all $i\in{\mathbb I}$,
%
\begin{eqnarray}
\label{req:numu}
\nu^r &\rightarrow&\nu \quad \mbox{and}\quad  \mu^r\rightarrow\mu \qquad
\mbox{as }  r\rightarrow\infty,
\\ \label{req:Jplus}
r (A \rho^r - C)&\rightarrow&
\theta \qquad  \mbox{as }   r\rightarrow\infty.
\end{eqnarray}
\end{assumption}

Let $M =\operatorname{diag}(\mu)$ and $\rho_i =\frac{\nu_i}{\mu
_i}$ for all
$i\in{\mathbb I}$. We note that (\ref{req:numu})--(\ref{req:Jplus}) imply
that $\rho^r\to\rho$ as $r\to\infty$ and $A\rho=C$.

\begin{rem}
Assumption \ref{HT} thus implies that {\it all} resources are
heavily loaded. We do not consider the case where
some resources are underloaded; however, we conjecture that the diffusion
approximation in this case would be as if these underloaded resources
were deleted from the model.
\end{rem}

We define fluid scaled processes $\overline N^r, \overline W^r,
\overline U^r, \overline T^r, \overline E^r, \overline S^r$ as
follows. For each~$r$ and $t\geq0$, let
%
\begin{eqnarray}\label{rfluone}
\overline N^r(t) &= &N^r(r t) /r, \qquad
\overline W^r(t) = W^r(r t) /r,
\\\label{rflutwo}
\overline U^r(t) &=& U^r (rt)/r, \qquad\hspace{3pt}
\overline T^r(t) = T^r (rt)/r,
\\\label{rfluthree}
\overline E^r(t) &=& E^r(rt)/r, \qquad\hspace{5pt}\overline S^r(t) = S^r(rt)
/r.
\end{eqnarray}
We define diffusion scaled processes $\hat N^r , \hat W^r, \hat
U^r, \hat E^r, \hat S^r$ as follows. For each~$r$ and $t\geq0$, let
%
\begin{eqnarray}\label{eqn:N}
\hat N^r (t) &=& \frac{N^r(r^2t)}{r},
\\\label{eqn:WN}
\hat W^r(t) &=& \frac{W^r (r^2t)}{r} = A(M^r)^{-1} \hat N^r(t) ,
\\
\hat U^r(t)&=& \frac{U^r(r^2t)}{r},
\\
\hat E^r(t ) &=& \frac{E^r(r^2t) -\nu^r r^2 t}{r},
\\
\hat S^r(t) &=& \frac{S^r(r^2t) -\mu^r r^2 t }{r}.
\end{eqnarray}
As $ E^r_i, S_i^r$, $i\in{\mathbb I}$, are independent Poisson
processes with parameters satisfying the convergence conditions
$(\ref{req:numu})$, it follows that we have the following
well-known functional central limit result \cite{Bi99}:
%
\begin{eqnarray} \label{clt}
(\hat E^r, \hat S^r) \Rightarrow(\tilde E, \tilde S) \qquad\mbox{as
} r\to\infty,
\end{eqnarray}
where $\tilde E$ and $\tilde S$ are independent ${\mathbf I}$-dimensional
Brownian motions starting from the origin with zero drift and
covariance matrices $\operatorname{diag}(\nu)$ and $\operatorname
{diag}(\mu)$, respectively.

Finally, we assume that, independent of (\ref{clt}),
$\hat W^r (0)$ converges in distribution as $r\to\infty$ to a
${\mathbf J}$-dimensional random
variable.

\section{Fluid model}\label{sec:fluid}

In this section, we recall some definitions and results established
in the prior work \cite{KW}. These will be needed for our statement
and proof of
multiplicative state space collapse.

\subsection{Fluid model solution}

A fluid model solution can be thought of as a formal limit of the
sequence $\{\overline N^r\}$ as $r\to\infty$. In fact, if one
assumes that~$\overline N^r(0)$ converges in distribution as
$r\rightarrow\infty$ to a random variable taking values in $\mathbb{R}
_+^{{\mathbf I}}$,
then one can show (see the Appendix of \cite{KW})
that $\{(\overline T^r, \overline E^r, \overline S^r,\allowbreak \overline
N^r, \overline U^r)\}$ is $C$-tight and any weak limit point
$(\overline T, \overline E, \overline S, \overline N, \overline U)$
yields a fluid model solution $\overline N$ a.s.

The following
notions are used in the definition of a fluid model solution given
below. A function $f=(f_1,
\ldots, f_{{\mathbf I}})\dvtx [0,\infty)\to\mathbb{R}^{{\mathbf I}}_+$
is said to be absolutely\vadjust{\goodbreak} continuous
if each of its components $f_i\dvtx [0,\infty)\to\mathbb{R}_+$, $i=1,
\ldots,
{\mathbf I}$, is absolutely continuous. A \textit{regular point} for an
absolutely continuous function $f\dvtx [0,\infty)\to\mathbb{R}^{{\mathbf
I}}_+$ is a
value of $t\in(0,\infty)$ at which each component of $f$ is
differentiable. [Since $f$ is absolutely continuous, almost every
time $t\in(0,\infty)$ is a regular point for $f$. Furthermore, $f$
can be recovered by integration from its a.e.~defined derivative.]

\begin{defn} \label{defn:fluid}
A fluid model solution is an absolutely continuous function
$n\dvtx [0,\infty)\to\mathbb{R}_+^{{\mathbf I}}$ such that at each
regular point $t>0$
for $n( \cdot )$, we have, for each $i\in{\mathbb I}$,
%
\begin{equation}\label{eq:diff}
\frac{d}{dt}  n_i(t)=
\cases{
\nu_i -\mu_i \Lambda_i(n(t)), &\quad \mbox{if } $n_i(t)>0,$ \cr
0, &\quad \mbox{if } $n_i(t)=0$
}
\end{equation}
and for each $j\in{\mathbb J}$,
%
\begin{equation} \label{eq:constraint}
\sum_{i\in{\mathbb I}_+(n(t))} A_{ji} \Lambda_i (n(t)) +\sum_{i\in
{\mathbb I}_0(n(t))} A_{ji}
\rho_i
\leq C_j ,
\end{equation}
where ${\mathbb I}_+(n(t)) =\{ i\in{\mathbb I}\dvtx n_i(t)>0\}$ and
${\mathbb I}_0(n(t))=\{i\in{\mathbb I}\dvtx n_i(t) =0\}$.
\end{defn}

\begin{rem}
Note that we are not assuming uniqueness of fluid model
solutions given the initial state.
\end{rem}

\subsection{Invariant manifold}

\begin{defn} A state $n_0\in\mathbb{R}_+^{{\mathbf I}}$ is called
\textit{invariant} (for the
fluid model) if there is a fluid model solution $n( \cdot )$ such
that $n(t) =n_0$ for all $t\geq0$. Let ${\mathcal M}_\alpha$ denote
the set
of all invariant states. We call ${\mathcal M}_\alpha$ the \textit{invariant
manifold}.
\end{defn}

\begin{rem}
Although $\alpha\in(0,\infty)$ is fixed throughout, we indicate the
dependence of $\mathcal M_\alpha$ on
$\alpha$ explicitly here as it will be useful later on when we explain
how the state space
for the proposed workload diffusion approximation varies with~$\alpha$.
\end{rem}

Various characterizations of the invariant states were given in
\cite{KW}. We summarize these in Theorem \ref{invariant} below. For
this, we need the following definitions.

For each $n\in\mathbb{R}_+^{{\mathbf I}}$, define $w(n)=(w_j(n)\dvtx j\in
{\mathbb J})$, to be
given by
%
\begin{equation} \label{req:defw}
w_j(n) = \sum_{i\in{\mathbb I}} A_{ji} \frac{ n_i}{\mu_i}, \qquad j
\in{\mathbb J}.
\end{equation}
We call $w(n)$ the \textit{workload} associated with $n$.

For each $w\in\mathbb{R}_+^{{\mathbf J}}$, define $\Delta(w)$ to be
the unique
value of $ n\in\mathbb{R}_+^{{\mathbf I}}$ that solves the following
optimization
problem:
%
\begin{eqnarray}\label{req:minF}
&&\mbox{minimize}\quad\hspace{2pt} {F(n)}\nonumber
\\
&&\mbox{subject to}\quad \sum_{i\in{\mathbb I}} A_{ji} \frac{
n_i}{\mu_i} \geq w_j,
\qquad j \in{\mathbb J},
\\
&&\mbox{over}\quad\hspace{25pt} n_i \geq0, \qquad i\in{\mathbb I},\nonumber
\end{eqnarray}
where
%
\begin{equation} \label{eq:F}
F(n) = \frac{1}{\alpha+1} \sum_{i \in{\mathbb I}} \nu_i {\kappa}_i
\mu_i^{\alpha-1} \biggl(\frac{n_i}{\nu_i}\biggr)^{\alpha+1},
\qquad
n\in\mathbb{R}_+^{{\mathbf I}}.
\end{equation}
(This function $ F$ was introduced in \cite{BM} as a Lyapunov function
for the fluid model. In fact, it is a Lyapunov function for
the original flow count process~$ N$ and can be used to show
positive recurrence of~$N$ when the average load
on each resource is less than its capacity.)
The function $\Delta$ has the two properties stated in
the next proposition.

\begin{proposition} \label{scaling}
The function $\Delta\dvtx \mathbb{R}_+^{{\mathbf J}} \to\mathbb
{R}_+^{{\mathbf I}}$ is continuous.
Furthermore, for each $w\in\mathbb{R}_+^{{\mathbf J}}$ and $c >0$,
%
\begin{equation}
\Delta(c w) = c\Delta(w).
\end{equation}
\end{proposition}

\begin{pf} The first property is proved in Lemma 6.3 of \cite{KW}. For the second
property, note that for $w\in\mathbb{R}_+^{{\mathbf J}}$ and $c>0$,
$c\Delta(w)$\vadjust{\goodbreak}
satisfies the constraints in~(\ref{req:minF}) with $cw$ in place of $w$ and so
%
\begin{eqnarray}\label{scaleDelta}
c^{\alpha+1}F(\Delta(w))=F(c\Delta(w))\geq
F(\Delta(cw)).
\end{eqnarray}
On the other hand,
by writing $w/c$ in place of $w$ in (\ref{scaleDelta}), we find
that
\[
c^{\alpha+1}F(\Delta(w/c))\geq F(\Delta(w))
\]
and then, by
replacing $c$ by $1/c$ and rearranging, we obtain
%
\begin{equation}\label{req:calpha}
F(\Delta(cw))\geq c^{\alpha+1}F(\Delta(w)).
\end{equation}
On combining (\ref{scaleDelta}) and (\ref{req:calpha}), we conclude that
\[
F(\Delta(cw))=F(c\Delta(w))
\]
and by uniqueness of the solution to (\ref{req:minF}),
we obtain the second property.
\end{pf}

\begin{theorem}\label{invariant}
The following are equivalent for $n\in\mathbb{R}_+^{{\mathbf I}}$:
\begin{longlist}
\item[(i)] $n$ is an invariant state (for the fluid model), that is,
$n\in
{\mathcal M}_\alpha$;
\item[(ii)] $\Lambda_i (n) = \rho_i $ for all $i\in{\mathbb
I}_+(n)=\{l\in
{\mathbb I}\dvtx n_l >0\}$;
\item[(iii)] there exists some $q\in\mathbb{R}_+^{{\mathbf J}}$ such that
%
\begin{equation}
\label{eq:invnq}
n_i = \rho_i \biggl( \frac{ \sum_{j \in{\mathbb J}} q_j
A_{ji}}{{\kappa}_i}
\biggr) ^{1/\alpha}
\qquad\mbox{for all } i \in{\mathbb I};
\end{equation}
\item[(iv)]
$n = \Delta(w(n))$.
\end{longlist}
\end{theorem}

\begin{pf} This follows immediately from Lemma 5.1 and Theorems 5.1,
5.3 of \cite{KW}.
\end{pf}

\begin{rem}
Note that if the conditions of Theorem~\ref{invariant}
are satisfied, then $p=q$ satisfies conditions~(\ref{eq:lm1})
and~(\ref{eq:lm2}), and thus $(q_j \dvtx j\in{\mathbb J})$ are dual
variables for the optimization problem~(\ref{eq:opt});
note that we use the fact that $A\rho=C$ for this.
We have chosen to use $q$ to denote the dual variables associated
with the invariant states, to distinguish them from the
dual variables~$p$ associated with arbitrary states $n$.
This distinction will be useful in our proof of convergence
to a diffusion process (see
Lemma \ref{thm:ctrlWN}), where we need to
distinguish the dual variables associated with
actual system states from the dual variables associated with nearby points
on the invariant manifold.
It is important to make this distinction because
when a system state
is near an invariant state and some
component of the system state is near zero,
it need not follow that the
dual variables associated with the two states are close.
\end{rem}

\begin{proposition} \label{Deltaw}
For each $w\in\mathbb{R}_+^{\mathbf J}$, $\Delta(w)\in\mathcal
{M}_{\alpha}$, that is,
$\Delta(w)$ is an invariant state.\vadjust{\goodbreak}
\end{proposition}
\begin{pf} Let $w\in\mathbb{R}_+^{\mathbf J}$. Since $\Delta(w)$ is the
unique optimal
solution to (\ref{req:minF}), it follows from Lemma 6.4 of Kelly and
Williams \cite{KW} that there exists $q\in\mathbb{R}_+^{{\mathbf
J}}$ such that
$\Delta(w)_i = \rho_i ( { \sum_{j \in{\mathbb J}} q_j
A_{ji}/{\kappa}_i} ) ^{1/\alpha}$ for all $ i \in{\mathbb
I}$. Then, by
Theorem \ref{invariant}, we have that $\Delta(w)$ is an invariant
state.
\end{pf}

\subsection{Asymptotic properties of fluid model solutions}\label{r:asy}

$\!$The next three propositions note some properties of fluid model
solutions that follow
from the analysis in \cite{KW} and that are used in our
proof of multiplicative state space collapse (see Theorem \ref{req:SSC} below).

\begin{proposition} \label{compactfluid}
For each $R\in(0,\infty)$, there is a constant $D(R)\in[R,
\infty)$ such that for any fluid model solution $n( \cdot )$
satisfying $|n(0)| \leq R$, we have $|n(t) |\leq D(R) $ for all
$t\geq0$.
\end{proposition}
\begin{pf} The proof of this proposition is implicit in the proof of
Theorem~5.2 in \cite{KW}.
\end{pf}

The next proposition states that fluid model solutions converge
uniformly to the invariant manifold ${\mathcal M}_\alpha$, where the
uniformity applies across all fluid model solutions that start
inside a compact subset of $\mathbb{R}^{{\mathbf I}}_+$.

\begin{proposition} \label{thm:convinv} Fix $R\in(0,\infty)$ and
$\varepsilon>0$.
There is a constant $T_{R,\varepsilon} \in[1,\infty) $ such that
for each fluid model solution $n( \cdot )$ satisfying $|n(0)|\leq
R$, we have
%
\begin{equation}
d(n(t),{\mathcal M}_\alpha) < \varepsilon\qquad\mbox{for all } t\geq T_{R,
\varepsilon}.
\end{equation}
\end{proposition}

\begin{pf}
The content of this proposition is the same as that of Theorem~5.2 in
\cite{KW}.
\end{pf}

\begin{proposition} \label{closetoMalpha}
For each $R\in(0,\infty)$ and $\varepsilon>0$, there exists $\delta>0
$ such that for any fluid model solution $n( \cdot )$ satisfying
$|n(0)|\leq R$ and $d(n(0),\allowbreak \mathcal M_\alpha) <\delta$, we have
$d(n(t),\allowbreak\mathcal M_\alpha) <\varepsilon$ for all $t\geq0$.
\end{proposition}
\begin{pf} This proposition follows from the proof of
Theorem 5.2 in \cite{KW}. For completeness, we provide a few details.
Fix $R>0$ and $\varepsilon>0$. By Proposition \ref{compactfluid},
there exists a compact set $B(R)$ in $\mathbb{R}_+^{{\mathbf I}}$ such
that $n(t)
\in B(R)$ for all $t\geq0$
for any fluid model solution $n$ satisfying $|n(0)|\leq R$.
Let $D=\{ u\in B(R)\dvtx d(u, \mathcal M_\alpha) \geq\varepsilon\}$.
As shown in \cite{KW},
there is a continuous function $H:\mathbb{R}_+^{{\mathbf I}}\to
\mathbb{R}_+$ that is zero
on $\mathcal M_\alpha$ and
strictly positive off $\mathcal M_\alpha$ such that $H(n(\cdot))$ is
nonincreasing for each
fluid model solution $n$. Let $\delta_1=\inf\{H(u)\dvtx u\in D\}$.
Then, $\delta_1>0$ and, by the properties of $H$, there exists $\delta>0$
such that whenever $n$ is a fluid model solution satisfying $|n(0)|\leq
R$ and
$d(n(0), \mathcal M_\alpha)\leq\delta$, we have $H(n(0))<\delta_1$.
Since $H(n(\cdot))$ is a nonincreasing function,
it then follows that $H(n(t)) <\delta_1 $ for all $t\geq0$. The
latter implies that $n(t) \notin D$ for all $t\geq0$
and so $d(n(t) , \mathcal M_\alpha) < \varepsilon$ for all $t\geq0$.
\end{pf}

The following corollary shows that fluid model solutions
starting on the invariant manifold stay at their
starting points for all time. We note this for
the reader's interest. We
do not use this corollary in our proofs.

\begin{corollary}
Suppose that $n(\cdot)$ is a fluid model solution such that
$n(0)\in\mathcal M_\alpha$.
Then, $n(t) =n(0)\in\mathcal M_\alpha$ for all $t\geq0$.
\end{corollary}

\begin{pf}
By Proposition \ref{closetoMalpha}, since
$d(n(0),\mathcal M_\alpha)=0$, we have $d(n(t),\allowbreak \mathcal M_\alpha)=0$
for all $t\geq0$ and hence $n(t)\in\mathcal M_\alpha$ for all
$t\geq0$. It follows from Theorem~\ref{invariant}
that for each $t>0$,
$\Lambda_i (n(t) )= \rho_i =\nu_i/\mu_i$ for all $i $ such that $n_i(t)>0$.
Then, by the fluid model dynamics
(\ref{eq:diff}), for each regular point
$t>0$ of $n(\cdot)$, we have
%
\begin{equation}
\frac{d}{dt} n_i(t) =\nu_i -\mu_i \rho_i =0 \\
\end{equation}
if $n_i(t) >0$ and the last equality also holds if $n_i(t) =0$.
It follows, since the absolutely continuous function
$n(\cdot)$ can be recovered from its almost
everywhere defined derivative, that $n(t) =n(0) $ for all
$t\geq0$.
\end{pf}

\section{Main results}

In this section, we describe the main results of this paper.
We begin with our result on \textit{multiplicative state space collapse}.
This is established using the asymptotic behavior of fluid model
solutions described in Section \ref{r:asy}.
Loosely speaking, multiplicative state space collapse shows that an
approximation for
$\hat N^r$ can be derived from one for $\hat W^r$ via the continuous
lifting map
$\Delta$ [see (\ref{req:minF})--(\ref{eq:F}) for the definition of
this map].
This lifting map can be nonlinear (for $\alpha\not=1$).
The multiplicative state space collapse result leads to a natural conjecture
for a diffusion approximation to~$\hat W^r$.
In the case $\alpha=1$,
assuming a mild local traffic condition and suitable initial
conditions,
we prove that the conjectured diffusion approximation is valid.
(In Section \ref{notone}, we indicate some of the
challenges associated with
establishing this conjecture for $\alpha\not=1$.)
When $\alpha=1$ and all of the weights for the bandwidth sharing
policy are equal (proportional fair sharing),
we use results of Harrison and Williams
\cite{HaWi87} and Williams \cite{RSS} to show that the diffusion has
a product form
invariant measure.
When this measure has finite total mass, this result
suggests an approximation for the stationary distribution of the flow count
process which we are able to extend
to the case where the document size distributions are finite mixtures of
exponential distributions and to some models with multi-path routing.
So as not to disrupt the flow
of results and associated discussion,
we defer the rather lengthy proofs of multiplicative state
space collapse and of the diffusion approximation to
Sections \ref{section:MSSC} and \ref{alphaHTL}, respectively.

\subsection{Multiplicative state space collapse}

\begin{defn}[(Multiplicative state space collapse)]\label{msscdef}
Multiplicative\break state~space collapse holds
(for the sequence of flow-level models described in Section~\ref{r:sequence}),
if, for each $T>0$,
%
\begin{equation} \label{conveq}
\frac{\|\hat N^r( \cdot ) -\Delta(\hat
W^r( \cdot ))\|_{T }}{\|\hat N^r( \cdot ) \|_{T}
\vee1 } \to0
\end{equation}
in probability as $r\to\infty$.
\end{defn}

\begin{rem}
We note here that in our form of multiplicative state space collapse,
the normalization (in the denominator) is in terms of the
flow count process,
whereas in Bramson's version for multiclass queueing networks \cite{BR},
it is in terms of a workload process. Furthermore, the lifting maps in
\cite{BR}
are all linear, whereas here, $\Delta$ can be nonlinear (for $\alpha
\not=1$).
\end{rem}

\begin{rem}
If (\ref{conveq}) holds without the factor in the denominator,
then state space collapse is said to hold.
Multiplicative state space collapse is more convenient for the purpose
of verification and if $\{\hat N^r \}$ (or $\{\hat W^r\}$) satisfies a compact
containment condition, then state space collapse follows from mutiplicative
state space collapse.
As was the case for open multiclass HL queueing networks \cite{Wi},
in establishing our diffusion approximation result for $\alpha=1$
under a mild local traffic condition, we
will show for this case that multiplicative state space collapse implies
state space collapse.
\end{rem}

The following theorem is one of the main results of this paper. It is
proved in Section
\ref{section:MSSC}.

\begin{theorem} \label{req:SSC} Assume that
%
\begin{equation} \label{conveq:ini}
|\hat N^r(0) -\Delta(\hat W^r(0))| \to0
\end{equation}
in probability as $r\to\infty$. Multiplicative state space collapse
then holds.
\end{theorem}

\subsection{Conjectured diffusion approximation}

We are interested in obtaining a diffusion approximation for the scaled
workload process $ \hat
W^r$. The multiplicative state space collapse result
can then be used to obtain a diffusion approximation for the scaled
flow count process $\hat N^r$.

For each $r$, define the double fluid scaled bandwidth allocation process
%
\begin{eqnarray}\label{barbarT}
\bar{\hspace{-1pt}\bar T}{}^r(t) = \frac{T^r(r^2t)}{r^2}, \qquad t\geq0.
\end{eqnarray}
Using (\ref{N})--(\ref{T}) and (\ref{req:W}) for the $r$th system,
the definitions of rescaled
processes and (ii) of
Proposition \ref{r:Lambdapro}, after some simple manipulations,\vadjust{\goodbreak} we obtain,
for all $t\geq0$,
%
\begin{eqnarray}\label{W}
\hat W^r(t) =
\hat W^r(0)+\hat X^r(t)+ \hat U^r(t),
\end{eqnarray}
where
%
\begin{eqnarray}\label{rXr}
\hat X^r(t) &=& A(M^r)^{-1}\bigl(\hat{E}^r(t)-
\hat{S}^r(\bar{\hspace{-1pt}\bar{T}}{}^r(t))\bigr)+r (A\rho^r-C)t,
\\
(\hat S^r(\bar{\hspace{-1pt}\bar{T}}{}^r(t)))_i &=& \hat
S^r_i(\bar{\hspace{-1pt}\bar T}{}^r_{\!i}(t)), \qquad i\in{\mathbb I},
\\
\hat{U}^r(t)&=& r\bigl(Ct-A\bar{\hspace{-1pt}\bar{T}}{}^r(t)\bigr)\nonumber
\\\label{rU}
&=&r^{-1}\int_0^{r^2t}\bigl(C-A\Lambda
(N^r(s))\bigr)\, ds
\\
&=& r \int_0^{t}\bigl(C-A\Lambda
(\hat{N}^r(s))\bigr)\, ds. \nonumber
\end{eqnarray}

If we
postulate that multiplicative state space collapse
implies state space collapse,
then, by formally (nonrigorously) passing to the limit
in the expression (\ref{W}) for $\hat W^r$,
we can obtain a natural conjecture for a diffusion
approximation to $\hat W^r$.
Immediately below, we give an informal description of how one might
arrive
at this conjecture. Following that, we
give a precise mathematical description of the diffusion
process and of the conjecture.

For the following informal description, which is used
to motivate the form of the conjectured diffusion approximation, we
postulate that
the sequence of processes $\{(\hat W^r, \bar{\hspace{-1pt}\bar{T}}{}^r, \hat U^r, \hat
E^r, \hat S^r)\}$
converges in distribution to a 5-tuple of continuous processes
$(\tilde W,T^*, \tilde U, \tilde E, \tilde S)$.
We also postulate that state space collapse (SSC)
holds
(not just multiplicative state space collapse).
From the convergence of $\{\hat W^r\}$, SSC and continuity of the
lifting map $\Delta$,
it follows that $\hat N^r$ converges in distribution to a continuous
process $\tilde N= \Delta(\tilde W)$ that lives
on the invariant manifold $\mathcal M_\alpha$. The fact that
$\hat W^r=A(M^r)^{-1}\hat N^r$ will yield
in the limit that
$\tilde W=A{M}^{-1}\tilde N$. By the characterization
of invariant states given in Theorem \ref{invariant}, it will then follow
that
for each $t\geq0$ and realization $\omega$,
there exists
$ q(t, \omega)\in\mathbb{R}_+^{{\mathbf J}}$ such that
%
\begin{eqnarray}
\tilde N_i(t,\omega)= \rho_i \biggl(\frac{\sum_{j\in\mathbb{J}}
q_j(t,\omega) A_{ji}}{\kappa_i}\biggr)^{1/\alpha} \qquad\mbox{for
all } i\in\mathbb{I}.
\end{eqnarray}
Consequently, $\tilde W=AM^{-1}\tilde N$ will live in the space
%
\begin{equation}
{\mathcal W}_{\alpha} = AM^{-1}{\mathcal M}_\alpha,
\end{equation}
where
%
\begin{eqnarray}\label{req:Na}
 \quad\mathcal M_\alpha=\biggl\{ n\in\mathbb{R}_+^{{\mathbf I}}\dvtx
n_i = \rho_i
\biggl(\frac{(q' A)_{i}}{{\kappa}_i}\biggr)^{{1/\alpha}}
\mbox
{for all } i\in{\mathbb I},
\mbox{ some } q\in\mathbb{R}_+^{{\mathbf J}}\biggr\}.
\end{eqnarray}
We call $\mathcal W_\alpha$ the \textit{workload cone.}
This is the state space for the conjectured diffusion approximation
$\tilde W$.

By the assumption made at the end of Section
\ref{r:sequence}, we know that $\hat W^r(0)$ converges in distribution,
independently of the primitive arrival and service processes.
We denote the limit distribution of $\hat W^r(0)$ by $\eta$.
Under the postulated convergence of $\hat W^r$, $\eta$
will be concentrated
on $\mathcal W_\alpha$.

We now turn our attention to the term $\hat X^r$ in the expression
(\ref{W})
for~$\hat W^r$.
Given the functional central limit theorem result
(\ref{clt})
for the diffusion scaled arrival and service processes $(\hat E^r, \hat
S^r)$, and the heavy traffic Assumption~\ref{HT},
if we postulate that the double fluid scaled allocation processes
$\bar{\hspace{-1pt}\bar{T}}{}^r $ achieve the nominal levels given by~%
$T^*(t)\equiv\rho t$
in the heavy traffic limit,
then $\hat X^r $ given by (\ref{rXr})
will converge in distribution
to the Brownian motion\break
$ AM^{-1}( \tilde E(\cdot) - \tilde S(T^*(\cdot))) +\theta(\cdot),
$
where $(\tilde S(T^*(\cdot)))_i = \tilde S_i(T^*_i(\cdot))$ for $i\in
{\mathbb I}$,
$\theta(t)=\theta t $ for all $t\geq0$ and $\theta$ is defined
in the heavy traffic Assumption \ref{HT}.
This Brownian motion starts from the origin,
and has drift $\theta$ and covariance matrix $AM^{-1}\operatorname
{diag}(\nu+\nu
)M^{-1}A'$,
where we have used the facts that $M=
\operatorname{diag}(\mu)$ and $\mu_i\rho_i =\nu_i$ for all
$i\in\mathbb{I}$ to compute the second term in the diagonal part of
the covariance
matrix expression.

On examining the representation (\ref{W}) for $\hat W^r$,
we see that it remains to conjecture properties for the postulated limit
$\tilde U$ of the
scaled unused capacity process $\hat U^r$ as $r\to\infty$.
The limit $\tilde U$ will inherit the nondecreasing property from the
$\hat U^r$.
The main issue is to determine where each of the components
of $\tilde U$ can increase.
For the prelimit process,
$\hat U^r$, to determine where its components increase,
we see from (\ref{rU}) that
it suffices to identify where each of the
components of $C-A\Lambda(\hat N^r(\cdot))$ is strictly
positive.
From Proposition
\ref{r:Lambdapro}, if $j\in\mathbb{J}$ such that
$C_j -\sum_{i\in\mathbb{I}} A_{ji} \Lambda_i(\hat N^r(t,\omega))>0$,
then
there
is a Lagrange multiplier $p^r(t, \omega)\in\mathbb{R}_+^{{\mathbf J}}$
such that $p^r_j(t,\omega)=0$ and
%
\begin{eqnarray}\label{rnlambda}
\hat N^r_i(t,\omega)& =&\Lambda_i(\hat N^r(t,\omega))\biggl(
\frac{\sum_{k\in\mathbb{J}}p_k^r(t, \omega)A_{ki} }{\kappa
_i}
\biggr)^{1/\alpha}
\qquad\mbox{for all } i\in\mathbb{I}.
\end{eqnarray}
[For this, we note that both sides of (\ref{rnlambda})
are zero if $\hat N_i^r(t,\omega)\!=\!0$.]
If $\Lambda(\hat N^r(t, \omega))$ is close to the nominal allocation
$\rho$,
then, by (\ref{rnlambda}), $\hat N^r(t, \omega)$ is near
%
\begin{eqnarray}
\mathcal M^j_\alpha & =&  \biggl\{n\in\mathbb{R}_+^{{\mathbf I}}\dvtx n_i
=\rho_i
\biggl(\frac{(q' A)_{i}}{{\kappa}_i}\biggr)^{{1/\alpha}}
\mbox
{for all } i\in{\mathbb I},\nonumber
\\[-8pt]\\[-8pt]
&&\qquad\qquad\mbox{\ \ \  some } q\in\mathbb{R}_+^{{\mathbf J}} \mbox{ satisfying } q_j=0
\biggr\}, \nonumber
\end{eqnarray}
the subset of
$\mathcal M_\alpha$ obtained by setting $q_j$ equal to zero in
(\ref{req:Na}), and then $\hat W^r =A(M^r)^{-1} \hat N^r $ will be close
to
%
\begin{equation} \label{Wsupj}
\mathcal W^j_{\alpha} = \{ AM^{-1} n\dvtx n \in\mathcal M^j_\alpha\},
\end{equation}
which we refer to as the $j$th
face of the workload cone $\mathcal W_\alpha$.
Thus, one might conjecture that in the limit, $\tilde U_j$
can increase only when $\tilde W$ is on the face~$\mathcal W^j_\alpha$.

Combining all of the above considerations leads to the informal conjecture
that $\hat W^r$ converges in distribution to a ${\mathbf J}$-dimensional
diffusion process
$\tilde W $ of the form $\tilde W(0) + \tilde X +\tilde U$.
The state space for $\tilde W$
is the workload cone $\mathcal W_\alpha$ and the initial
distribution of $\tilde W$ is given by $\eta$.
In the interior of $\mathcal W_\alpha$,
the increments of $\tilde W$ are given by
the increments of a Brownian motion
$\tilde X$ with drift $\theta$ and
covariance matrix
%
\begin{equation} \label{r:gamma}
\Gamma= 2AM^{-1} \operatorname{diag}(\nu) M^{-1}A',
\end{equation}
that starts from the origin.
The process $\tilde W$
is confined to the cone $\mathcal W_\alpha$ by instantaneous
``pushing'' at
the boundary of $\mathcal W_\alpha$.
The direction of push allowed when $\tilde W$ is
on the boundary face $\mathcal W^j_\alpha$ is $\gamma^j$, the unit
vector parallel to the positive $j$th coordinate axis in $\mathbb
{R}_+^{{\mathbf J}}$,
and the cumulative amount of push in that direction is given by the
$j$th component of the continuous
nondecreasing process $\tilde U$. (At the intersections of boundary
faces, the combined effect of pushing on the
individual boundary faces is a push in
the direction of a~convex combination of the
pushing directions available from the intersecting boundary faces.)
The direction of push on a given
boundary face is usually called a \textit{direction of reflection},
whereas it might be better thought of as a~direction of regulation
for the process. (The term ``reflection'' comes from the fact that
in one dimension, when the drift is zero,
the construction of such a regulated process
from a Brownian motion can be achieved by a~mirror reflection; although
this type of
construction does not generally apply in higher dimensions, the
term ``reflection'' is still used.)

We now introduce a precise definition for the
conjectured diffusion approximation to the workload process $\hat W^r$.
(The delicate issue of existence and uniqueness for this process
is discussed below.)
This definition will be used
in giving a precise statement of our
conjecture.
Here, $\theta$ is a vector in~$\mathbb{R}^{{\mathbf J}}$, $\Gamma$
is given
by (\ref{r:gamma}), $\gamma^j$ is the unit vector parallel
to the positive $j$th coordinate axis and
$\eta$ is a Borel probability measure on $\mathcal W_\alpha$.

\begin{defn} \label{DefSRBM}
A \textit{semimartingale reflecting Brownian motion} that lives in the cone
$\mathcal W_{\alpha}$, has direction of reflection $\gamma^j$ on the
boundary face~$\mathcal W^j_{\alpha}$ for each $j\in{\mathbb J}$, has drift
$\theta$ and covariance matrix $\Gamma$, and has initial
distribution $\eta$ on $\mathcal W_{\alpha}$
is an adapted,
${\mathbf J}$-dimensional process $\tilde W$ defined on some filtered
probability space $(\Omega,{\mathcal{F}},\{{\mathcal{F}}_t\},P)$
such that
\begin{longlist}
\item[(i)] $P$-a.s., $\tilde W(t)=\tilde W(0)+\tilde X(t)+
\tilde U(t)$ for all $t\geq0$,
\item[(ii)] $P$-a.s., $\tilde W$ has continuous paths, $\tilde W(t)\in
\mathcal W_{\alpha}$ for all $t\geq0$ and $\tilde W(0)$ has
distribution $\eta$,
\item[(iii)] under $P$,
\begin{longlist}
\item[(a)] $\tilde X$ is a ${\mathbf J}$-dimensional Brownian motion
starting from
the origin with drift $\theta$ and covariance matrix $\Gamma$;
\item[(b)] $\{\tilde X(t) -\theta t, {\mathcal{F}}_t, t\geq
0\}$ is a martingale,
\end{longlist}
\item[(iv)] for each $j\in{\mathbb J}$, $\tilde U_j$ is an
adapted, one-dimensional process such that
$P$-a.s.,
\begin{longlist}
\item[(a)] $\tilde U_j(0)=0;$
\item[(b)] $\tilde U_j \mbox{ is continuous and
nondecreasing;}$
\item[(c)] $\tilde U_j(t)=\int_{0}^t1_{\{\tilde W(s)\in
\mathcal W^j_{\alpha} \}}\,d\tilde U_j(s)$ for all $t\geq0$.
\end{longlist}
\end{longlist}
\end{defn}

\begin{rem}
We call a process $\tilde W$ satisfying
the above properties an \textit{SRBM} associated with the
data
$(\mathcal W_\alpha,\theta,\Gamma, \{\gamma^j:j\in{\mathbb J}\},
\eta)$.
Here, ``adapted'' means adapted to the filtration
$\{\mathcal F_t\}$.
We note that this filtration
need not be the one generated by $\tilde X$---it can be larger.
However,
it can always be taken to be the filtration generated by $\tilde
W,\tilde X,\tilde U$.
The martingale condition on $\tilde X$ is included here as
we are using a ``weak'' definition of the process.
This martingale property is needed in establishing uniqueness in law
for an SRBM.
The term ``semimartingale'' refers to the fact that $\tilde W$ is the
sum of a continuous martingale and a continuous process that
is locally of bounded variation.
Condition (iv)(c) corresponds to the condition that $\tilde U_j$
can only increase when $\tilde W$ is on the boundary face $\mathcal
W^j_\alpha$.
\end{rem}

We now give a precise statement of our conjecture.

\begin{conjecture} \label{conjecture}
Suppose that the limit distribution of $\hat W^r(0)$
is $\eta$, a~probability distribution on $\mathcal W_\alpha$ endowed
with the
Borel $\sigma$-algebra, and
suppose that
\[
|\hat N^r(0) -\Delta(\hat W^r(0))|\to0 \qquad
\mbox{in probability as } r\to\infty.
\]
Then,
$\hat W^r$ converges in distribution as $r\rightarrow\infty$ to a
process $\tilde W$ that is an SRBM associated with
the data $(\mathcal W_\alpha,\theta, \Gamma, \{\gamma^j\dvtx j\in
{\mathbb J}\}
, \eta)$.
\end{conjecture}

We shall prove that Conjecture \ref{conjecture} holds when $\alpha
=1$, provided that a mild local
traffic condition holds. We state this result in the next subsection.
In the remainder of the current subsection, we indicate some of the
challenges associated with constructing a rigorous proof of the conjecture.

When $\alpha=1$ (corresponding to weighted proportional fair sharing),
we can express the cone $\mathcal W_\alpha$ in the
simple form
%
\begin{eqnarray}\label{simpleW} \mathcal{W}_1=\{ABA'q\dvtx q\in\mathbb
{R}^{\bf J}_+\},
\end{eqnarray}
where $B$ is an ${\bf I}\times{\bf I}$ diagonal matrix with the
$i$th diagonal entry being $\frac{\nu_i}{\mu_i^2\kappa_i}>0$.
Thus, $\mathcal W_1$ is a \textit{polyhedral} cone. Since
$A$ has full row rank and~$B$ is a diagonal matrix
with strictly positive diagonal entries,
$ABA'$ is a~linear bijection
between $\mathbb{R}_+^{{\mathbf J}}$ and $\mathcal W_1 $.
It follows from this that $\mathcal W_1$ is a~simple
polyhedral cone.
Necessary and
sufficient conditions for the existence and uniqueness in
law of SRBMs
living in simple polyhedral domains have been given by Dai and Williams
\cite{DW}. Under these conditions, the SRBM is a diffusion,
that is, a continuous strong Markov process. It will turn out that the
conditions of \cite{DW}
are satisfied by our data when $\alpha=1$. For $\alpha\not=1$
and ${\mathbf J}=2$, the workload cone is a wedge which is still a simple
polyhedral cone. However,
in general, for $\alpha\not=1$ and $\mathbf J >2$, the
workload cone $\mathcal W_\alpha$ is not a~polyhedral cone
(it has curved boundaries).
In this case, we have some partial (unpublished) results on existence
and uniqueness
for SRBMs.
The main impediment to obtaining a general result is that
boundary faces
can meet in cusp-like singularities, making it challenging to
even determine whether the process can escape from the cusp
and whether it can do so in a~unique manner
(see Section \ref{notone} for an example).

Even if one has existence and uniqueness of the SRBM,
for any proof of the conjecture,
there are a number of other challenges to overcome.
First, one needs to establish $C$-tightness of the sequence of triples
$\{(\hat W^r, \hat X^r, \hat U^r)\}$. This is largely an issue
of the $C$-tightness of $\{\hat U^r\}$.
One also needs to show that multiplicative state space collapse
implies state space collapse.
One of the
most challenging\vspace{1pt} aspects is to show that for any possible limit
$(\tilde W, \tilde X, \tilde U)$ of the sequence $\{(\hat W^r, \hat
X^r, \hat U^r)\}$,
for each $j\in\mathbb{J}$,
the process $\tilde U_j$ can only increase when $\tilde W$ is
on the boundary face $\mathcal W^j_\alpha$.
Indeed, in our informal use of (\ref{rnlambda}) to arrive at
our conjecture, we neglected the fact that
$\Lambda_i(\hat N^r(t,\omega))$
need not be near $\rho_i$ when some component of $\hat N^r(t,\omega)$
is near zero [recall that $\Lambda_i (n) $ need not be continuous
when $n_i$ is zero].
Also, the notions of ``nearness'' and ``closeness'' used loosely
in our informal description are not necessarily uniform.
However,
if $\hat N^r_i(t, \omega)$ is at or near zero when $p^r_j(t,\omega
)=0$ for
an $i $ such that $A_{ji}>0$,
then we can show that $\hat W^r(t, \omega)$
is near the boundary face~$\mathcal W^j_\alpha$.
To take advantage of this observation, in the case when
we prove the conjecture ($\alpha=1$), we will assume that a mild local
traffic condition holds.

In summary,
the main reasons that we are able to treat the case $\alpha=1$
are that the existence and uniqueness theory for the limit
diffusion process is in place \cite{DW} and there is an associated
invariance principle \cite{KaWi}
which, loosely speaking, is a perturbation result
telling us that processes such as~$\hat W^r$,
that satisfy perturbed versions of the defining conditions
for an SRBM, are close in distribution to an SRBM.
In particular, for $\alpha=1$, the invariance principle of \cite{KaWi}
takes care of establishing the $C$-tightness of $\{(\hat W^r,
\hat X^r, \hat U^r)\}$ and, in the presence of the uniqueness in
law of the SRBM \cite{DW}, it implies
convergence in distribution of $\hat W^r$ to an SRBM.
In the case $\alpha\not=1$ and ${\mathbf J}=2$, $\mathcal W_\alpha$
is a wedge (a polyhedral cone) and our proof
for $\alpha=1$ can
be extended to this case.
In as yet unpublished work,
we have been able to establish uniqueness in law of the
SRBM and to establish an invariance principle
for some cases where $\alpha\not=1$ and ${\mathbf J}>2$.
However, some cases, especially when boundary faces meet
in cusp-like singularities, are as yet unresolved.
We summarize the situation for $\alpha\not=1$ in
Section \ref{notone}. However, because
of the partial nature of our results so far, we leave
the description of these further developments
to future work.

\subsection{\texorpdfstring{Diffusion approximation for weighted proportional fair
sharing (\mbox{$\alpha\,{=}\,1$})}
{Diffusion approximation for weighted proportional fair
sharing (alpha=1)}}

The following condition is used in the next theorem. This condition
can be interpreted as a \textit{local traffic} assumption, under
which each resource has at least one route that only uses that
resource.

\begin{assumption}[(Local traffic)] \label{localtraffic}
For each $j\in{\mathbb J}$, there exists at least one $i\in{\mathbb
I}$ such that
$A_{ji}>0$ and
$A_{ki} =0 $ for all $k\not=j$.
\end{assumption}

The following theorem is proved in Section \ref{alphaHTL}.

\begin{theorem} \label{thm:HTL} Assume that $\alpha=1$ and that
the local traffic Assumption~\ref{localtraffic} holds. Suppose that
the limit distribution of $\hat W^r(0)$ as $r\to\infty$
is $\eta$ (a~probability measure on $\mathcal W_1$) and that
$|\hat
N^r(0) -\Delta(\hat W^r(0))
| \to0$ in probability as $r\rightarrow\infty$. Then,
$(\hat W^r,\hat N^r)$ converges in distribution as $r\rightarrow
\infty$ to a conti\-nuous
process $(\tilde W, \tilde N)$, where $\tilde W$ is an SRBM
with data $(\mathcal W_1,\theta,
\Gamma, \{ \gamma^j\dvtx j\in{\mathbb J}\},\allowbreak \eta)$
and $\tilde N=\Delta(\tilde W)$.
\end{theorem}

In the case $\alpha=1$, the lifting map $\Delta$ is in fact a linear
map on $\mathcal W_1$, given by
%
\begin{equation} \label{delwalpha}
\Delta(w) = \operatorname{diag}(\rho)\operatorname{diag}(\kappa
)^{-1} A'(ABA')^{-1} w,  \qquad
w\in\mathcal W_1.
\end{equation}
Indeed, for $\alpha=1$ and $w\in\mathcal W_1$,
if $q=(ABA')^{-1} w$ and $n$ is given by the right-hand side of (\ref
{delwalpha}),
then $w(n)=AM^{-1}n =w$, by the definition of $B=M^{-1}\operatorname
{diag}(\nu
)\operatorname{diag}(\kappa)^{-1}M^{-1}$,
and $n=\operatorname{diag}(\rho)\operatorname{diag}(\kappa
)^{-1}A'q$ so that (\ref{eq:invnq})
holds. Then, by Theorem~\ref{invariant}, we conclude that
\[
n=\Delta(w(n))=\Delta(w)
\]
and, hence, (\ref{delwalpha}) holds.
By the remark following Theorem \ref{invariant},
the
\mbox{$(q_j \dvtx j\,{\in}\,{\mathbb J})$} defined above are dual
variables for the optimization problem~(\ref{eq:opt}).
It follows that we can associate a \textit{process of dual random
variables} $\tilde Q$ with the SRBM~$\tilde W$ of Theorem~\ref
{thm:HTL}, as follows.
Given $\tilde W$ as in Theorem~\ref{thm:HTL},
define
%
\begin{equation}
\tilde Q= (ABA')^{-1} \tilde W.
\end{equation}
This process $\tilde Q$ inherits an SRBM structure from $\tilde W$.
In fact, $\tilde Q$ is a~semimartingale reflecting
Brownian motion living in $\mathbb{R}_+^{{\mathbf J}}$,
having the form
%
\begin{equation}\label{structQ}
\tilde Q (t) = \tilde Q(0) +(ABA')^{-1} \tilde X(t) +
(ABA')^{-1}\tilde U(t), \qquad t\geq0,
\end{equation}
where
the Brownian motion $(ABA')^{-1} \tilde X $ has
drift $(ABA')^{-1}\theta$ and covariance matrix $(ABA')^{-1} \Gamma
(ABA')^{-1}$,
and
$\tilde U_j$ can increase only when $\tilde Q_j$ is zero, $j\in
{\mathbb J}$.
The direction of reflection on the boundary face of $\mathbb
{R}_+^{{\mathbf J}}$ is
defined by the $j$th column of the matrix $(ABA')^{-1}$.
The initial distribution of $\tilde Q$ is
obtained by applying the linear transformation $(ABA')^{-1}$ to the
distribution~$\eta$.
(For the formal definition of such an SRBM, where
reflection directions are not in general parallel to
coordinate directions,
see
\cite{Wi98}.)

\begin{figure}[t]

\includegraphics{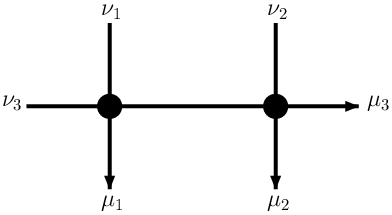}

\caption{A linear network with two resources and three routes.}
\label{fig:linear}
\end{figure}

As an illustration of Theorem \ref{thm:HTL}, consider a
two-resource linear network operating under a weighted proportional fair
sharing policy ($\alpha=1$). This network
is depicted in Figure \ref{fig:linear}. It has two resources
and three routes. Each resource has a route that passes only
through that resource and there is also a route that passes
through both resources.

\begin{figure}[t]

\includegraphics{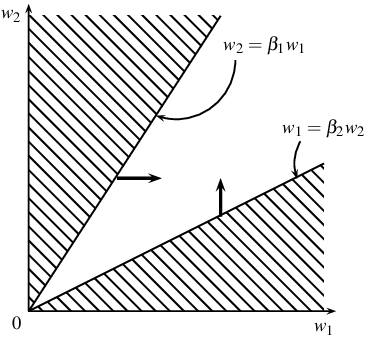}

\caption{The workload cone ${\mathcal W}_1$ for a linear network with two resources
and three routes is an infinite wedge in the positive quadrant of
$\mathbb{R}^2$.
A finite portion of the wedge is
shown above (between the two shaded regions).
Under the lifting map $\Delta$, points $(w_1, w_2)$ on the boundary
$w_2 = \beta_1 w_1$ of the wedge are mapped to points $(n_1, n_2, n_3)$,
where $n_1=0$ (and the corresponding $q\in\mathbb{R}_+^2$ has $q_1=0$);
similarly, points $(w_1, w_2)$ on the boundary $w_1 =\beta_2w_2$ are
mapped to points $(n_1, n_2, n_3)$, where $n_2 =0$ (and the corresponding
$q\in\mathbb{R}_+^2$ has $q_2=0$).}
\label{fig:wedge}
\end{figure}

The workload cone in this case is a wedge in
$\mathbb{R}^2_+$ that has the following representation:
%
\begin{equation}\mathcal{W}_1=
\cases{\left. \left[\matrix{w_1 \cr w_2}\right]
=\left[\matrix{
 \biggl(\dfrac{\nu_1}{\mu_1^2\kappa_1}+\dfrac{\nu
_3}{\mu_3^2 \kappa_3}\biggr) q_1+
\dfrac{\nu_3}{\mu_3^2\kappa_3}q_2
\cr \biggl(\dfrac{\nu_2}{\mu_2^2\kappa_2}+ \dfrac{\nu_3}{\mu_3^2
\kappa_3}\biggr)q_2+
\dfrac{\nu_3}{\mu_3^2\kappa_3}q_1
}\right]\dvtx q\in\mathbb{R}_+^2 \right\}}.
\end{equation}
Let
\[
\beta_1 =
\biggl(1 + \frac{\nu_2\mu_3^{2}\kappa_3 }{\nu_3\mu_2^{2}\kappa_2}
\biggr)
\quad\mbox{and }
\quad\beta_2 =
\biggl(1 + \frac{\nu_1\mu_3^{2} \kappa_3}{\nu_3\mu_1^{2} \kappa_1}
\biggr)
.
\]
It can be easily computed
that the two boundary faces of the wedge $\mathcal{W}_1$ have the
following expressions:
\[
\mathcal{W}_1^1 = \{w \in\mathbb{R}_+^2\dvtx
w_1\geq0, w_2 = \beta_1
w_1 \}
\]
and
\[
\mathcal{W}_1^2= \{ w\in\mathbb{R}_+^2\dvtx
w_2\geq0, w_1 = \beta_2
w_2 \}.\vadjust{\goodbreak}
\]
The wedge
is depicted in Figure \ref{fig:wedge}.
For a linear network, the local traffic condition holds
automatically. Thus, the conclusion of Theorem \ref{thm:HTL}
applies, provided
$\hat W^r(0)$ converges in distribution to
a random variable with distribution $\eta$ concentrated on $\mathcal
W_1$ and
$|\hat N^r(0) -\Delta(\hat W^r(0)
) | \to0$ in probability as $r\rightarrow\infty$.
For example, this holds if each of the systems indexed by $r$
starts empty and $\eta$ is the point mass at the origin
in $\mathbb{R}_+^{2}$.
The limiting SRBM $\tilde W$ lives in the wedge $\mathcal W_1$ and is
confined there by reflection (or pushing) at the boundary.
Reflection occurs in the horizontal direction
(corresponding to resource $1$ underutilizing capacity) on the bounding face
where $w_2=\beta_1w_1 $. The interpretation of this is that although
there is no work for
resource $1$ on route $1$, there is work for this resource on route
$3$. However,
congestion at resource $2$,
through the nature of the bandwidth sharing policy,
is preventing resource $1$ from working at its full capacity.
Similarly, vertical reflection (corresponding to resource $2$
underutilizing capacity)
on the bounding face $w_1=\beta_2w_2 $ is interpreted to mean that
congestion at resource $1$ is preventing resource $2$ from
working at its full capacity.
Thus, the shape of the workload space indicates the
entrainment of resources, whereby congestion at
some resources may prevent other resources from
working at their full capacity.
We note that when $\kappa_3\rightarrow\infty$, the
upper boundary of $\mathcal{W}_1$ tends to the vertical axis and the
lower boundary
tends to the horizontal axis, hence the wedge expands to the whole
quadrant, approaching the situation with full utilization of the resources.

\subsection{Product form stationary distribution for proportional fair sharing}

In this subsection, we
prove a result which shows
that when $\alpha=1$ and $\kappa_i=1$ for all $i\in\mathbb{I}$ (proportional
fair sharing), an SRBM $\tilde W$ with the
properties described in Theorem~\ref{thm:HTL} has a
product form invariant measure.
(This is a~result for SRBMs and so does not require the local traffic condition
a~priori.)

\begin{theorem} \label{productform}
Suppose that $\alpha=1$ and $\kappa_i=1$ for all $i\in\mathbb{I}$.
Let $\pi$ be the measure on $\mathcal W_1$
that is absolutely continuous with respect to Lebesgue
measure with
density given by
%
\begin{equation} \label{prodform}
p(w) = \exp(\upsilon\cdot w ), \qquad w\in\mathcal W_1,
\end{equation}
where
%
\begin{equation}
\upsilon= 2\Gamma^{-1} \theta.
\end{equation}
The product form measure $\pi$ is an invariant measure
for the SRBM having state space $\mathcal W_1$,
directions of reflection $ \{\gamma^j\dvtx j\in{\mathbb J}\}$,
drift
$\theta$ and covariance matrix $\Gamma$.
This measure is integrable over $\mathcal W_1$ if and
only if $\theta_j <0$ for all $j\in{\mathbb J}$
and then, after normalization, it defines the unique stationary
distribution for the SRBM.
\end{theorem}

\begin{pf}
Sufficient conditions for a reflecting Brownian motion in a
simple polyhedral domain
to have a product form invariant measure were determined by
Williams
in \cite{RSS}, building on solutions of a related analytic
problem obtained by Harrison and Williams in \cite{HaWi87}.
In these works, the covariance matrix for the
process is the identity matrix.
In order to apply these results, we need to
perform a linear transformation
to transform an SRBM with
covariance matrix
$\Gamma$
into one whose covariance
matrix is the identity matrix.
We perform this transformation below; the computations
are straightforward, though a little tedious, as we need to normalize the
resulting
directions of reflection to have inward normal components of unit length,
to facilitate use of the results in \cite{RSS}.
Similar manipulations were carried out in the proof of Theorem 23 in
\cite{HaWi87b}, where the reflection directions had a special form.

Before introducing the transformation,
we obtain an alternative representation for the simple convex polyhedron
$\mathcal W_1$.
By (\ref{simpleW}),
\[
\mathcal W_1 = \{ ABA' q\dvtx  q\in\mathbb{R}_+^{{\mathbf J}}\}=\{w\in
\mathbb{R}^{{\mathbf J}}\dvtx
(ABA')^{-1}w\in\mathbb{R}_+^{{\mathbf J}}\} .
\]
Here,
since $\kappa_i=1$ for all $i$, we have that
$B= M^{-1} \operatorname{diag}(\nu) M^{-1}$ and,
by (\ref{r:gamma}),
$ABA' = \frac{1}{2}\Gamma$.
Thus,
%
\begin{eqnarray} \label{altw}
\mathcal W_1 & = &\{ w\in\mathbb{R}^{{\mathbf J}}\dvtx \Gamma^{-1} w \in
\mathbb{R}_+^{{\mathbf J}}
\} ,
\end{eqnarray}
where we have used the fact that $q\in\mathbb{R}_+^{{\mathbf J}}$
if\vadjust{\goodbreak}
and only if $2
q\in\mathbb{R}_+^{{\mathbf J}}$.

We now define the linear transformation to be applied to the SRBM.
Let~$\Upsilon$ denote the diagonal matrix that has the same diagonal
entries as $\Gamma^{-1}$.
Let~$L$ be the rotation matrix whose rows are the orthonormal
eigenvectors of the covariance matrix $\Gamma$ and $D$
be the corresponding diagonal matrix of eigenvalues
such that $\Gamma= L'DL$, where $L'=L^{-1}$.
Let
$\tilde W$ be an SRBM with state space $\mathcal W_1$, drift
$\theta$, covariance matrix $\Gamma$ and directions of reflection
given by $\{\gamma^j\dvtx j\in{\mathbb J}\}$, with a decomposition
as in (i) of Definition \ref{DefSRBM}.
Let $V=D^{-{1/2}}L$ and define
$\tilde Z = V\tilde W$. Then, $\tilde Z$
is an SRBM in the simple convex polyhedron
\begin{eqnarray*}
\mathcal Z_1 &=&\{Vw\dvtx w\in\mathcal W_1 \}
= \{ z\in\mathbb{R}^{{\mathbf J}} \dvtx V^{-1} z \in\mathcal W_1\} \\
&=& \{z\in\mathbb{R}^{{\mathbf J}} \dvtx \Gamma^{-1} V^{-1} z \in\mathbb
{R}_+^{{\mathbf J}}\} \\
&=&\{ z\in\mathbb{R}^{{\mathbf J}}\dvtx \tilde n^j\cdot z \geq0 \mbox{
for all } j\in
{\mathbb J}\},
\end{eqnarray*}
where, for each $j\in{\mathbb J}$,
$\tilde n^j $ is given by the $j$th
row of the matrix
%
\begin{eqnarray} \label{req:Theta}
\Theta= \Upsilon^{-{1/2}} \Gamma^{-1} V^{-1} =\Upsilon
^{-{1/2}}\Gamma^{-1} L' D^{{1/2}}
\end{eqnarray}
and $\tilde n^j$ is
the unit inward normal to the $j$th face of
$\mathcal Z_1$. (The matrix $\Upsilon^{-{1/2}}$
is used to normalize the $\tilde n^j$ to be of unit length.)
The process $\tilde Z$ inherits the following decomposition
from $\tilde W$:
%
\begin{eqnarray}
\tilde Z(t) = \tilde Z(0) + V\tilde X(t) + R\tilde Y(t) , \qquad t\geq
0,   \mbox{ a.s.,}
\end{eqnarray}
where $R=VH^{-1} =D^{-{1/2}} L H^{-1} $, $H = \Upsilon^{
{1/2}}$ is the diagonal matrix with the same
diagonal entries as $\Theta V$
and $\tilde Y = H\tilde U$ is a continuous, nondecreasing process that
starts from zero and can increase only when $\tilde Z$ is on
the $j$th face of $\mathcal Z_1$.
The columns of the matrix $R$ are the directions of reflection for
$\tilde Z$
on the faces of $\mathcal Z_1$, normalized so that
the inward normal component of each direction of reflection
has unit length.
(The matrix $H^{-1}$ is used to achieve the normalization.)
The matrix $R$ has the form
%
\begin{eqnarray}
R= \Theta' + \Xi',
\end{eqnarray}
where
$\Theta$ is the matrix specified in (\ref{req:Theta}) (whose rows are
the inward unit
normals to the faces of $\mathcal Z_1$)
and $\Xi$ is the matrix whose rows consist of the components
of the (normalized) directions of reflection that are tangent to each
of the faces
of $\mathcal Z_1$. In particular, the diagonal entries of $\Theta\Xi'
$ are all zero.
In this context, the sufficient condition given in \cite{RSS} for
existence of a product form invariant measure for $\tilde Z$ is the
so-called ``skew symmetry
condition,''
%
\begin{eqnarray} \label{r:sss}
\Theta\Xi' + \Xi\Theta' = 0.
\end{eqnarray}
The vector form of (\ref{r:sss}) is
%
\begin{eqnarray}
\tilde n^i \cdot\tilde s^j + \tilde s^i \cdot\tilde n^j =0 \qquad
\mbox{for all } i, j \in{\mathbb J},
\end{eqnarray}
where, for each $i\in{\mathbb J}$,
$\tilde n^i $ is the inward unit normal to the $i$th
face of $\mathcal Z_1$ and~$\tilde s^i $ is the tangential
component of the (normalized) direction of reflection on that face.
A~geometric interpretation of the skew symmetry condition is provided in
\cite{HaWi87}. We refer the reader to that paper, especially Section~9, for
a full discussion. Briefly, for a simple convex polyhedron
as we have here, the skew symmetry condition
can be shown to be equivalent to a~local condition.
This local condition
can be stated in words as requiring that near the
intersection of any two distinct faces of $\mathcal Z_1$,
the component of the reflection vector on one face that is directed
toward the intersection of the faces is balanced
on the other face by the component of the direction
of reflection for that second face which
is of the same magnitude as~the component on the first face, but it is directed
away from the intersection of the faces.
Indeed, under the skew symmetry condition, with probability one,
the SRBM $\tilde Z$ will not hit the intersection
of any two faces when started away from that set, as is proved in \cite{RSS}.

The left-hand side of condition (\ref{r:sss})
is equal to
\begin{eqnarray*}
\Theta(R-\Theta') + (R'-\Theta)\Theta' & =&- 2 \Theta\Theta'
+\Theta R + R' \Theta'
\\
&=& - 2 \Upsilon^{-{1/2}}\Gamma^{-1} L' D^{{1/2}}
D^{{1/2}} L \Gamma^{-1} \Upsilon^{-{1/2}}
\\
&& {}+ \Upsilon^{-{1/2}}\Gamma^{-1} L' D^{{1/2}}
D^{-{1/2}} L H^{-1}
\\
&&{}+ H^{-1} L' D^{-{1/2}} D^{{1/2}} L \Gamma^{-1}
\Upsilon^{-{1/2}}
\\
&=& - 2\Upsilon^{-{1/2}} \Gamma^{-1} \Upsilon^{-
{1/2}}
\\
&&{}+ \Upsilon^{-{1/2}} \Gamma^{-1} H^{-1} + H^{-1}
\Gamma^{-1} \Upsilon^{-{1/2}}
\\
&=& 0,
\end{eqnarray*}
where we have used the facts that $\Gamma=L'DL$, $L'L=I$ and
$H=\Upsilon^{{1/2}}$.

Thus, the skew symmetry condition (\ref{r:sss})
holds and it follows from Theorem~1.2 of \cite{RSS} that
the SRBM $\tilde Z$ has a product form invariant measure
with a density relative to Lebesgue measure that is proportional to
$\exp( \beta\cdot z) $, $z\in\mathcal Z_1$, where
$\beta= 2(I-\Theta^{-1}\Xi)^{-1} V\theta=2V\theta$.
Note, for this, that
\begin{eqnarray*}
\Theta^{-1} \Xi & =& D^{-{1/2}}L\Gamma\Upsilon^{{1/2}}
(H^{-1} L' D^{-{1/2}} - \Upsilon^{-{1/2}} \Gamma^{-1}
L' D^{{1/2}}) \\
&=& D^{-{1/2}} L\Gamma\Upsilon^{{1/2}} (\Upsilon
^{-{1/2}}L'D^{-{1/2}} -\Upsilon^{-
{1/2}}L'D^{-1}LL' D^{{1/2}})
\\
&= & D^{-{1/2}} L \Gamma\Upsilon^{{1/2}} (\Upsilon
^{-{1/2}}L'D^{-{1/2}} -\Upsilon^{-{1/2}}L'D^{-{1/2}})
\\
&=& 0,
\end{eqnarray*}
where we have used the facts that $\Gamma^{-1}=L'D^{-1}L$ and $L'L=I$.
Furthermore, by Corollary 1.1 of \cite{RSS},
if the exponential density is integrable over $\mathcal Z_1$,
then, after normalization, it yields the unique stationary
distribution for $\tilde Z$.

By inverting the linear transformation $V$,
we can transform this result back to one for $\tilde W$.
Noting that for $z =Vw$,
%
\begin{equation}
\beta\cdot z =\beta\cdot Vw = 2 \theta' V' V w = 2\theta' \Gamma
^{-1} w ,
\end{equation}
we conclude that
(\ref{prodform}) is an invariant density for the original SRBM $\tilde W$
and if this is integrable over $\mathcal W_1$, then, after normalization,
it yields the unique stationary distribution for $\tilde W$.
Using the representation (\ref{altw})
for $\mathcal W_1$, we see that this density
will be integrable over $\mathcal W_1$
if and only if
$\exp( \theta\cdot q)$, $q\in\mathbb{R}_+^{{\mathbf J}}$ is
integrable over
$\mathbb{R}_+^{{\mathbf J}}$, which occurs if and only if $\theta
_j<0$ for each
$j\in{\mathbb J}$.
\end{pf}

\begin{rem}
The product form invariant measure of Theorem \ref{productform}
is remarkable,
yet the proof of the result gives little insight into why the reflection
directions and covariance matrix in the particular case $\alpha= 1,
\kappa_i = 1$ for $i=1,\ldots, \bf I$ should allow the result to hold.
We simply note here that the authors first suspected
that a product form result might be found after observing that
it is possible to describe a network of queues with a product form
stationary distribution whose conjectured Brownian model
approximation has the same directions of reflection and
covariance matrix as the Brownian model approximation
under study in this paper, provided that $\alpha= 1$ and $\kappa_i = 1$
for $i=1,\ldots, \bf I$, that is, the case of proportional fair sharing.
Earlier connections between product form queueing networks and
proportional fairness have been explored by
\cite{BP,MasRob2002}
and the relationship between these several connections seems
a rich area for further study.
\end{rem}

The product form of the density (\ref{prodform}) does not
imply that, when $\theta_j <0$ for all $j\in{\mathbb J}$,
the components of the SRBM $\tilde W$ are independent under the
stationary distribution for the SRBM since, in general, the
cone $\mathcal W_1$ is not an orthant. Independence can, however,
be deduced for the components of the SRBM $\tilde Q$ of dual
variables.

\begin{corollary} \label{productformq}
Suppose that the assumptions of Theorem \ref{thm:HTL} hold, that
$\kappa_i=1$ for all $i\in\mathbb{I}$
and that $\theta_j <0$ for all $j\in{\mathbb J}$.
Let $(\tilde W, \tilde N)$ be the process identified in
Theorem \ref{thm:HTL}. The SRBM $\tilde Q=2\Gamma^{-1} \tilde W$ of
dual variables
then has a unique stationary distribution and this distribution
has a density relative to Lebesgue measure that is proportional to
$\exp( \theta\cdot q) $, $q\in\mathbb{R}_+^{\bf J}$. Under this
stationary distribution, the components of $\tilde Q$
are independent and $\tilde Q_j$ is exponentially
distributed with parameter $- \theta_j$ for each $j\in{\mathbb J}$.
\end{corollary}

\begin{pf} This result is immediate from Theorem \ref{productform} upon
applying the linear transformation
$ (ABA')^{-1}= 2 \Gamma^{-1}$ to transform
$\tilde W$ into $\tilde Q$.
\end{pf}

Under the assumptions of Corollary \ref{productformq},
from Theorem \ref{thm:HTL}, (\ref{delwalpha}) and the definition of
$\tilde Q$, we have that
%
\begin{equation}
\tilde N = \operatorname{diag}(\rho) A' \tilde Q
\end{equation}
and it follows that the stationary distribution of $\tilde N$
can be expressed as a~linear combination of independent exponential random
variables.
Thus, resource $j$ has associated with it a
dual random variable $\tilde Q^s_j$, for $j\in{\mathbb J}$ (here, the
superscript
of $s$ signals that the random variable is associated with the stationary
distribution). These dual variables are independent and exponentially
distributed with parameters
$-\theta_j$ for $j\in{\mathbb J}$
and under its stationary distribution, the $i$th
component of $\tilde N$ is proportional to the sum of the
dual variables associated with the resources used by route $i$.
This suggests the following simple approximation
for the stationary distribution of the \textit{unscaled} network,
that is, the flow-level
model of Section \ref{secflow}. The~stationary approximation
is
%
\begin{equation} \label{approx}
N^s_i \approx\rho_i \sum_{j \in{\mathbb J}} Q^s_j A_{ji},
\end{equation}
where $Q^s_j, j \in{\mathbb J}$, are independent and $Q^s_j$ is
exponentially distributed with parameter
$C_j - \sum_{i\in{\mathbb I}} A_{ji} \rho_i$.
We want to emphasize that (\ref{approx}) above is merely a formal approximation
involving various implicit assumptions such as existence of a
stationary distribution for
$ N$ and formal unraveling of the heavy traffic scaling and limit
procedure. In particular, this involves an interchange of limits
and the approximation in (\ref{approx})
involves errors that may be substantial (e.g.,
of order $r$).
However, the following observation of Massouli\'{e} and Roberts
\cite{RM}, which yields an exact stationary
distribution for an unscaled linear network,
suggests that there is cause for optimism regarding the approximation
(\ref{approx}).

\begin{example}
Consider a linear network with ${\mathbf J}$ resources, where the set of
resources is
${\mathbb J}= \{ 1,2,\ldots, {\mathbf J}\}$. Let the
set of routes be labeled
${\mathbb I}= \{ 0, 1,2,\ldots, {\mathbf J}\}$, where we
use the symbol $i=0$ for the route
$ \{ 1,2,\ldots, {\mathbf J}\}$ through every resource
and, for $i=1,2,\ldots,{\mathbf J}$, we use the symbol $i$ for the
route $ \{ i \}$ through the single resource $i$.
The local traffic Assumption \ref{localtraffic} thus
holds and we assume $\alpha=1$ and $\kappa_i=1$ for all $i\in\mathbb{I}$.
Suppose that $C_j=1, j=1,2,\ldots,{\mathbf J}$, and that
$\rho_0+\rho_j < 1, j=1,2,\ldots,{\mathbf J}$. The stationary
distribution for $(N_0, N_1,\ldots,N_{\mathbf J})$ is then given by \cite{RM}
%
\begin{equation}
\pi(n_0, n_1,\ldots,n_{\mathbf J}) =
\frac{\prod_{j=1}^{{\mathbf J}} (1-\rho_0 - \rho_j)}{(1-\rho
_0)^{{\mathbf J}-1}}
\left(\matrix{\displaystyle\sum_{i=0}^{{\mathbf J}}n_i \vspace{6pt}\cr n_0}\right) \prod_{i=0}^{{\mathbf J}}
\rho_i^{n_i},
\end{equation}
where $n_0, \ldots, n_{\mathbf J}$ each run through the nonnegative integers.
Summing this formula over $n_0$ and using the negative binomial
expansion gives the marginal distribution for
$(N_1, N_2,\ldots,N_{\mathbf J})$ as
%
\begin{equation}
\pi(n_1, n_2,\ldots,n_{\mathbf J}) =
\prod_{i=1}^{{\mathbf J}} \frac{(1-\rho_0 - \rho_i)}{(1-\rho_0)}
 \biggl(\frac{\rho_i}{1-\rho_0} \biggr)^{n_i}.
\end{equation}
Thus, under the stationary distribution, $N_1, N_2,\ldots, N_{\mathbf J}$
are independent and~$N_i$ is geometrically
distributed
with mean $\rho_i / (1-\rho_0 - \rho_i)$. This accords
remarkably well with the approximation (\ref{approx}),
under which $N_i$ would be
approximated by an exponentially distributed random variable
with the same mean.

For this particular example, as observed by \cite{RM}, the Markov
chain $N$
is reversible and, consequently, the
stationary distribution $\pi$ for $N$ is insensitive to the document
size distributions, depending only on their means.
(For a description of general system structures to which such
insensitivity results
apply, see
\cite{SCH} and references therein.)
\end{example}

As a complement to the insensitivity result mentioned in the above
example, we
note that
the stationary distribution result of Corollary~\ref{productformq}
can be extended
to the situation where the document sizes are finite mixtures of exponentials.
Indeed,
such a flow-level model
may be realized by collapsing an extended (exponential) model. The
extended model
is obtained by
splitting each route in the original model into finitely many ``copies''
so that each copy uses the same
set of resources and has the same weights as the original route, but
where each copy
has its own Poisson arrival process
and distinct exponential document size distribution.
The relative arrival rates for the copies determine the proportions for
the mixture of
exponential distributions associated with the original route.
More precisely,
consider an extended (exponential) flow-level model in which,
for each $i=1, \ldots, \bf I$,
routes $i(1), i(2),\ldots,i(K_i)$ identify finitely many identical routes
(subsets of the
resources in ${\mathbb J}$)
whose associated weights are also identical.
In the following, parameters and processes associated with this extended
model will have a superscript of ${}^\dag$ appended.
In particular, $N^\dag$ will denote the flow count process.

Under the assumptions of Theorem \ref{thm:HTL} (heavy traffic, $\alpha=1$,
local traffic condition
and initial state space collapse) for the extended flow-level model, there
is an
SRBM
approximation $\tilde W^\dag$
for the workload and an attendant approximation $\tilde N^\dag=\Delta
(\tilde W^
\dag)$
for the flow count process.
Furthermore, there is an SRBM process $\tilde Q^\dag$
of dual random variables.
Under the additional assumptions of Corollary~\ref{productformq}
(all of the weights are one and $\theta^\dag_j<0$
for each $j\in{\mathbb J}$),
the stationary distribution of $\tilde Q^\dag$ is such that the components
are independent and $\tilde Q^\dag_j$ is exponentially distributed with
parameter $-\theta^\dag_j$ for each $j\in{\mathbb J}$.
Thus, under the assumptions
of Corollary~\ref{productformq}, the stationary distribution of
%
\begin{equation} \label{r:ndag}
\tilde N^\dag=\operatorname{diag}(\rho^\dag) (A^\dag)' \tilde
Q^\dag
\end{equation}
can be expressed as a linear combination of exponential random variables.

For the extended model,
the Poisson arrival rate for route $i(k)$ is $\nu^\dag_{i(k)}$ and the
exponential service time parameter for this route is
$\mu^\dag_{i(k)}$.\vadjust{\goodbreak}
Then, the $\bf I$-dimensional collapsed process $N$ defined by
%
\begin{equation} \label{r:collapse}
N_i =\sum_{k=1}^{K_i} N^\dag_{i(k)} ,\qquad i=1,\ldots, \bf I,
\end{equation}
has the same distribution as the flow count process in a flow-level model
with $\bf I$ routes, where the Poisson arrival rate for route
$i=1,\ldots, \bf I$
is
%
\begin{equation}\label{eq:nm1}
\nu_{i} = \sum_{k=1}^{K_i} \nu^\dag_{i(k)}
\end{equation}
and the document size distribution for route $i=1, \ldots, \bf I$ is a mixture
of $K_i$ exponentials with parameters $\mu^\dag_{i(k)}$ and proportions
$\nu^\dag_{i(k)}/\nu_i$ for $k=1, \ldots, K_i$. The mean of this
document size distribution is
$1/\mu_i$, where
$\mu_i$ is defined by
%
\begin{equation}\label{eq:nm2}
\frac{1}{\mu_{i}} = \sum_{k=1}^{K_i} \frac{\nu^\dag_{i(k)}}{\nu
_i}\frac{1}{\mu^\dag_{i(k)}}
\end{equation}
for $i=1,\ldots, \bf I$. Let $\rho_i = \nu_i / \mu_{i}$ for
$i=1,\ldots,
\bf I$.

Now, consider the exponential analog of the collapsed
flow-level model where the finite mixture of exponential distributions
for route $i$
is replaced by a single exponential distribution
with parameter $\mu_i$ for $i=1,\ldots, \bf I$.
The nominal load placed on resource $j$ in the extended model
is $(A^\dag\rho^\dag)_j$ and this is the same as the nominal load
$(A\rho)_j$ placed on $j$ in the exponential analog.
It follows that a sequence of extended flow-level models
satisfying the assumptions of Corollary \ref{productformq},
where the limiting value in (\ref{req:Jplus}) of
the heavy traffic Assumption \ref{HT}
is denoted by $\theta^\dag$,
will have a parallel sequence of exponential analogs
whose limiting value $\theta$ in (\ref{req:Jplus})
will be precisely the same as $\theta^\dag$.
Consequently,
the stationary distribution of $\tilde Q^\dag$ is the same as
that for the process $\tilde Q $ of dual
random variables associated with the sequence of exponential
analogs.

Then, combining (\ref{r:ndag}) and (\ref{r:collapse}),
we obtain the following as an approximation for the collapsed flow
count process:
\begin{eqnarray*}
\tilde N_i = \sum_{k=1}^{K_i}
\rho^\dag_{i(k)}\sum_{j\in{\mathbb J}} A^{\dag}_{ji(k)} \tilde
Q^\dag_j
= \rho_i\sum_{j\in{\mathbb J}} A_{ji} \tilde Q^\dag_j \qquad  \mbox{for }
i=1,\ldots, \bf I \rm,
\end{eqnarray*}
that is, $\tilde N= \operatorname{diag}(\rho) A' \tilde Q^\dag$.
Since the stationary distributions for $\tilde Q^\dag$ and $\tilde Q$
are the same,
it follows that the stationary distribution for $\tilde N$ is the same,
whether one considers the collapsed
flow-level model (where document sizes are distributed as finite
mixtures of exponentials)
or the associated exponential analog (where document sizes are
distributed as exponentials).
Similarly, the
formal approximation (\ref{approx}) for the unscaled network
is the same for both.
Note, however, that the stochastic processes $\tilde W^\dag$, $\tilde
N$ and
$\tilde Q^\dag$,
and, in particular, their covariance matrices are, in general, different
from those associated with the exponential analog.

It is natural to conjecture an extension of the above
insensitivity results to general document size distributions.
However, even an extension to phase-type distributions would
require generalization of the diffusion approximation results to
flow-level models with feedback and treatment of more
general distributions would appear to require
a significantly more elaborate
stochastic model (see \cite{growil06}) to keep track of residual
document sizes.

\subsection{Multi-path routing}\label{multipathroute}

In this subsection, we describe a generalization of the earlier model
that allows multi-path routing. In our initial description of
the model, we shall use a different notation for the
sets of routes and resources. This will allow
our eventual results to be expressed in a notation consistent with that
used elsewhere in the paper.

Suppose that we now interpret $i \in{\mathbb I}$ as a source--destination
pair and let~$\mathbb{K}$ be the set of routes.
Let ${\mathbf I}=|{\mathbb I}|$ and let ${\mathbf K} = |\mathbb{K}|$.
We suppose that~$\mathbb{K}$ is partitioned into ${\mathbf I}$ nonempty
subsets, each associated with a single source--destination pair. Let
$H_{ik} = 1$ if route $k \in{\mathbb{K}}$ is associated with
source--destination pair $i \in{\mathbb I}$ and let $H_{ik} = 0$ otherwise.
Thus, $H$ is an ${\mathbf I}\times\mathbf K$ matrix containing only
zeros and ones, and $\sum_{i \in{{\mathbb I}}} H_{ik} = 1$ for each
$k \in\mathbb{K}$.

There remain finitely many resources, but now labeled by
$l\in\mathbb{L}$, with capacities $(\bar{C}_l\dvtx l \in
\mathbb{L})$ that are all strictly positive and finite.
A route $k$ is a~nonempty subset of $\mathbb{L}$
(interpreted as the set of resources used by route~$k$). We assume
that $\mathbb{L}$ and $\mathbb{K}$ are both nonempty and finite.
Let $\mathbf L =|\mathbb{L}|$, the total number of resources.
Let $\bar{A}$ be the $\mathbf L\times{\mathbf
K}$ matrix containing only zeros and ones, defined such that
$\bar{A}_{lk}=1$ if resource $l$ is used by route $k$
and $\bar{A}_{lk}=0$
otherwise. Our assumption that each route $k$ identifies a
\textit{nonempty} subset of $\mathbb{L}$ implies that no column of
$\bar{A}$ is identically zero.

It is assumed that a new document
arrives to source--destination pair $i$ at each
jump time of a Poisson process that
has rate parameter $\nu_i>0$ and that each such document has an
exponentially distributed size with mean $1/\mu_i$, where $\mu_i\in
(0,\infty)$. These document sizes are assumed to be independent of
one another and to be independent of all arrival times of documents.
Let $\rho_i = \nu_i/\mu_i, i \in{\mathbb I}$.

Given a fixed parameter $\alpha\in(0, \infty)$ and strictly
positive weights \mbox{$({\kappa}_i\dvtx i\in{\mathbb I})$}, if $N_i(t)$ denotes the
(random) number of documents being transferred between
source--destination pair $i$ at
time $t$ for each $i\in{\mathbb I}$ and $N(t) =(N_i(t)\dvtx\allowbreak i\in{\mathbb
I})$, then the
bandwidth allocated to source--destination pair $i$ at time~$t$ is
given by $\Lambda_i (N(t))$ and this bandwidth is shared equally
among the~$N_i(t)$ transfers in progress between
source--destination pair $i$.
The function $\Lambda= (\Lambda_i\dvtx i\in{\mathbb I})\dvtx
\mathbb{R}^{{\mathbf I}}_+\to\mathbb{R}^{{\mathbf I}}_+$ is defined
such that for
each $n\in\mathbb{R}^{{\mathbf I}}_+$, $\Lambda_i(n) =0$ for $i\in
{\mathbb I}_0(n)
\equiv\{m \in{\mathbb I}\dvtx n_m =0\}$, and when
$ {\mathbb I}_+(n) \equiv\{m\in
{\mathbb I}\dvtx n_m >0\}$ is nonempty, $\Lambda^+(n) \equiv(\Lambda
_i(n) \dvtx i\in
{\mathbb I}_+(n))$
is the unique value\vadjust{\goodbreak} of $\Lambda^+
=(\Lambda_i\dvtx i\in{\mathbb I}_+(n))$ that solves the
optimization problem
\begin{eqnarray}\label{eq:opt1a}
&&\mbox{maximize}\quad G_{n}(\Lambda^+)\nonumber
\\
&&\mbox{subject to}  \quad  \sum_{k} \bar{A}_{lk}
y_k \leq\bar{C}_l, \qquad l \in\mathbb{L},\nonumber
\\[-8pt]\\[-8pt]
&&\hspace{33pt}\qquad \sum_{k} H_{ik} y_k = \Lambda_i ,
\qquad
i \in{\mathbb I}_+(n) ,\nonumber
\\
&&\mbox{over} \hspace{24pt}\quad y_k \geq0, \qquad k \in{\mathbb{K}}, \qquad
\Lambda_i \geq
0, \qquad i\in{\mathbb I}_+(n),\nonumber
\end{eqnarray}
where $G_{n} (\Lambda^+)$ is again given by the definition (\ref{req:G}).

Without loss of generality, we may assume that in the solution
of (\ref{eq:opt1a}), $y_k=0$ for those $k$ such that $H_{ik}=0$ for all
$i\in{\mathbb I}_+(n)$. With this convention, allocations
$y_k ,k\in\mathbb{K}$, associated with the unique optimal value
$\Lambda^+(n)$ can be interpreted as bandwidth allocations for the
routes which sum to give the bandwidth allocations to the
source--destination pairs.

The above optimization problem reduces to the earlier problem
(\ref{eq:opt})
in the case where ${\mathbf I}= \mathbf K$ and $H$ is the identity matrix,
that is, the case of a~single route for each source--destination pair.
More generally, the following proposition allows us to reduce
to the optimization problem (\ref{eq:opt}), even in the multi-path case.
This observation was previously made in
\cite{KEL}, Section 3.3.

Let
\[
{\mathcal Y} = \{y \in\mathbb{R}_+^{\mathbf K} \dvtx \bar{A}y \leq\bar
{C} \}.
\]

\begin{proposition} \label{r:represent}
There exists a representation
%
\begin{equation}\label{eq:rep}
H {\mathcal Y} = \biggl\{ (\Lambda_i, \enspace i\in{\mathbb I}) \dvtx
\Lambda_i \geq0, \enspace i\in{\mathbb I}, \quad
\sum_{i\in{\mathbb I}} A_{ji} \Lambda_i \leq
C_j, \enspace j \in{\mathbb J}\biggr\},
\end{equation}
where ${\mathbb J}$, $A$ and $C$ can be chosen so that $C$
has positive elements, $A$ has nonnegative elements and
no column of $A$ is identically zero.
\end{proposition}

\begin{pf}
The set
$\mathcal Y$ is the intersection of the half-spaces $\{y\in
\mathbb{R}^{\mathbf K} \dvtx
(\bar{A}y)_l \leq\bar{C}_l\}$, $ l \in\mathbb{L}$, and the
nonnegative orthant,
$\mathbb{R}_+^{\mathbf K}$. Since
no column of $\bar{A}$ is identically zero, the set $\mathcal Y$ is
bounded and is thus the convex hull of a finite number of extreme
points. Hence, $ H {\mathcal Y}$ is the convex hull of a finite
number of extreme points or, equivalently, the bounded
intersection of a finite set of half-spaces.

Next, we explore further the geometry of $ H {\mathcal Y}$.
First,
$ H {\mathcal Y} \subset\mathbb{R}_+^{\mathbf I}$
since the elements of $H$ are nonnegative. Also,
$ 0 \in H {\mathcal Y}$ since $ 0 \in\mathcal Y$.
Indeed, for $\delta$ positive and small enough,
$(\Lambda_i = \delta_i, i\in{\mathbb I}) \in H {\mathcal Y}$
for all $0 \leq\delta_i < \delta$, $i\in{\mathbb I}$, since
for $\varepsilon $ positive and small enough,
$(y_k = \varepsilon _k, k \in{\mathbb{K}}) \in\mathcal Y$
for all $0 \leq\varepsilon _k < \varepsilon $, $ k\in\mathbb{K}$,
and no row of $H$ is identically zero.
Thus, $ H {\mathcal Y}$
is bounded by the hyperplanes bounding the nonnegative orthant,
plus finitely many other hyperplanes, none of which\vadjust{\goodbreak}
contains the origin. Thus, there exists
a representation of the form~(\ref{eq:rep})
for some choice of ${\mathbb J}$, $A$ and $C$.
Since $ 0 \in H {\mathcal Y}$,
the elements of~$C$ are nonnegative.
Choose a minimal representation, where the set ${\mathbb J}$
is of minimal cardinality. Then, the intersection of the hyperplane
%
\begin{equation}\label{eq:rep1}
\biggl\{ (\Lambda_i, \enspace i\in{\mathbb I}) \dvtx
\sum_{i\in{\mathbb I}} A_{ji} \Lambda_i = C_j
\biggr\}
\end{equation}
with $ H {\mathcal Y}$ must have a nonempty
interior relative to the hyperplane, for each $j \in{\mathbb J}$.

Next, we show that if $\Lambda\in H {\mathcal Y}$ and
$0 \leq\Lambda' \leq\Lambda$,
then $\Lambda' \in H {\mathcal Y} $.
Let $y \in\mathcal Y$
be such that $\Lambda= H y$ and, for each $k \in\mathbb{K}$,
let $i(k)$ be the unique
index $i$ such that $H_{ik}=1$. Let $y'_k =
(\Lambda'_{i(k)}/\Lambda_{i(k)})y_k$ for $k \in\mathbb{K}$. Then,
$y' \in\mathcal Y$ and $H y' = \Lambda'$, so
$\Lambda' \in H {\mathcal Y} $.

We have seen that the intersection of the
hyperplane~(\ref{eq:rep1}) with $ H {\mathcal Y}$ has a nonempty
interior relative to the hyperplane. Choose a point $\Lambda$
in this relative interior. Then, $\Lambda_i >0$ for $i\in{\mathbb I}$.
Fix $i\in{\mathbb I}$ and choose $\Lambda'$ so that $0< \Lambda'_i <
\Lambda_i$,
with $\Lambda'_{i'}=\Lambda_{i'}$ for $i' \in{\mathbb I}\setminus\{
i\} $.
Then, $0 \leq\Lambda' \leq\Lambda$, so
$\Lambda' \in H {\mathcal Y}$. Hence, the
hyperplane~(\ref{eq:rep1}) must have $A_{ji} \geq0$.

Finally, note that $C_j>0$ for $j \in{\mathbb J}$ since
the hyperplane~(\ref{eq:rep1})
does not contain the origin, and no
column of $A$ is identically zero since $ H {\mathcal Y}$
is bounded.
\end{pf}

The represention~(\ref{eq:rep}) of $ H {\mathcal Y}$ allows us to
elide the variables $y$ from the optimization
definition of $\Lambda^+(n)$ and to deduce that the
unique solution $(\Lambda_i\dvtx i\in{\mathbb I}_+(n))$ to the
optimization problem (\ref{eq:opt1a}) is also the
unique solution to the optimization problem (\ref{eq:opt}).
With this generalization to multi-path routing, the matrix
$A$ may no longer
only contain zeros and ones.
However, the proofs of Theorem \ref{thm:HTL} and Corollary
\ref{productformq} still go through in this
enhanced generality, provided that the local traffic
condition is satisfied.

\begin{figure}

\includegraphics{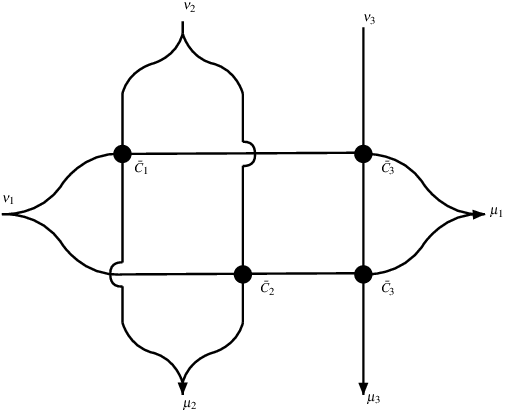}

\caption{A network with three source--destination pairs, four resources
and five routes. Each of the first two source--destination pairs
has two routes available to it.}
\label{fig:routing}
\vspace*{-3pt}
\end{figure}

As an illustration, consider the network depicted
in Figure \ref{fig:routing}, operating under the
proportional fair
sharing policy ($\alpha=1, \kappa_i =1$ for all $i$),
where the labeled resource capacities satisfy
$\bar{C}_3 < \bar{C}_1, \bar{C}_2$.
It has three source--destination pairs, four resources and five routes.
The matrices $A$ and $C$
can then be chosen as
\begin{eqnarray*}
A=\left(
\matrix{ 1 &1 & 0 \cr \tfrac{1}{2} & 0 & 1
}\right),
\qquad
C =\left(
\matrix{ \bar{C}_1 + \bar{C}_2 \cr \bar{C}_3
}\right).
\end{eqnarray*}
These matrices
may be viewed as expressing the \textit{generalized cut constraints}
%
\begin{eqnarray}\label{eq:gencut1}
\rho_1 + \rho_2 &\leq& \bar{C}_1 + \bar{C}_2,
\\\label{eq:gencut2}
\tfrac{1}{2} \rho_1 + \rho_3 &\leq& \bar{C}_3
\end{eqnarray}
apparent from Figure~\ref{fig:routing}.
In the first of these constraints, the capacities of
the first and second resources are pooled. Resource pooling
subject to generalized cut constraints is a common
phenomenon in stochastic processing networks with routing~\cite{LAWS}.
Theorem 3.1 of~\cite{LAWS} provides an alternative
form of the representation~(\ref{eq:rep}) and provides
references to algorithms to calculate matrices $A$ and $C$.

In the above example, we have been able to construct the matrix $A$
so that a subset of its columns forms the identity matrix and thus the
local traffic Assumption \ref{localtraffic} holds.
However,
in general, it is more difficult to verify the local traffic condition
for networks with multi-path routing
than for networks without routing choices.

Suppose that the local traffic Assumption \ref{localtraffic} holds
for the matrix $A$.
Define heavy traffic as in Assumption \ref{HT}
and assume that $\theta_{j} <0$ for $j \in{\mathbb J}$.
Corollary~\ref{productformq} then applies directly.
We can illustrate the
formal approximation~(\ref{approx})
for the example shown in Figure \ref{fig:routing}:
under the approximation,~$ Q^s_1$ and $ Q^s_2$,
the random dual variables associated with constraints
(\ref{eq:gencut1}) and~(\ref{eq:gencut2}),
respectively, are independent and exponentially
distributed, $ Q^s_1$ with parameter $\bar{C}_1 + \bar{C}_2 - \rho_1
- \rho_2$
and $ Q^s_2$ with parameter $\bar{C}_3 - \frac{1}{2}\rho_1 -
\rho_3$.

\subsection{\texorpdfstring{Discussion of the conjecture for $\alpha\not=1$}
{Discussion of the conjecture for alpha /=1}}\label{notone}
Consider the two-resource linear network depicted in Figure \ref{fig:linear},
operating under a weighted
$\alpha$-fair bandwidth sharing policy for $\alpha\neq1$. Then, the
workload cone $\mathcal W_\alpha$ is still a wedge in $\mathbb
{R}^2_+$ that
has the following
representation:
%
\begin{equation}\mathcal{W}_{\alpha}=
\cases{\left.\left[ \matrix{w_1 \cr w_2}\right]
=\left[
\matrix{ \dfrac{\nu_1}{\mu_1^2\kappa_1^{{1/\alpha
}}}q_1^{{1/\alpha}}+
\dfrac{\nu_3}{\mu_3^2\kappa_3^{{1/\alpha}}}(q_1+q_2)^{
{1/\alpha}}
\vspace{2pt}\cr \dfrac{\nu_2}{\mu_2^2\kappa_2^{{1/\alpha}}}q_2^{
{1/\alpha}}+
\dfrac{\nu_3}{\mu_3^2\kappa_3^{{1/\alpha}}}(q_1+q_2)^{
{1/\alpha}}
}
\right]\dvtx q\in\mathbb{R}_+^2\right\}}.\vadjust{\goodbreak}
\end{equation}
The two boundary faces of this wedge have the
following expressions:
\[
\mathcal W_\alpha^1 =
\cases{\left. \left[\matrix{ w_1 \cr w_2
}\right]
\dvtx w_1\geq0, w_2 = \biggl(1 + \dfrac{\mu_2^{-2} \nu_2}{\mu
_3^{-2} \nu_3}
\biggl(\dfrac{\kappa_3}{\kappa_2}\biggr)^{1/\alpha}\biggr)w_1
\right\}}
\]
and
\[
\mathcal W_\alpha^2 =
\matrix{\biggl\{ \left[\matrix{ w_1 \cr w_2}\right]
\dvtx w_2\geq0, w_1 = \biggl(1 + \dfrac{\mu_1^{-2} \nu_1}{\mu
_3^{-2} \nu_3}
\biggl(\dfrac{\kappa_3}{\kappa_1}\biggr)^{1/\alpha}\biggr)w_2
\biggr\}}.
\]
In this case, since the two-dimensional wedge $\mathcal W_\alpha$ is still
a polyhedral cone, the proof of Theorem \ref{thm:HTL} can be easily
modified to apply
for $\alpha\not=1$, with~$\mathcal W_\alpha$ used in place of
$\mathcal W_1$
there.
When the weights $\kappa_i, i=1,2,3,$ are all equal, the wedge
$\mathcal{W}_{\alpha}$ does not depend on $\alpha$.
However, when the $\kappa_i$'s are
not all equal, the wedge $\mathcal{W}_{\alpha}$ does depend on
$\alpha$.
However, as $\alpha$ goes to infinity, the quotients involving the
weights $\kappa_i$ in the expressions for $\mathcal{W}_\alpha^j,
j=1,2 $,
tend to one and there is a limiting wedge which
is the same as that obtained when the $\kappa_i$ are all equal. On
the other hand, as $\alpha$ tends to zero, the limit of the upper
boundary depends on whether $\kappa_3>\kappa_2$ (tends to the
vertical axis) or $\kappa_3<\kappa_2$ (tends to the 45 degree ray
from the origin).
 Similarly, for the lower boundary, if $\kappa_3>\kappa_1$, it
tends to the horizontal axis and if $\kappa_3<\kappa_1$, it tends to
the 45 degree ray from the origin. Hence, as $\alpha$ tends to zero,
if $\kappa_3> \max\{\kappa_1,
\kappa_2\}$, the wedge expands to the
whole quadrant and if $\kappa_3< \min\{\kappa_1, \kappa_2\}$, the wedge
contracts to the 45 degree ray from the origin.
Note that even when the wedge expands to the whole quadrant,
the components of the diffusion workload do not become independent
as the covariance matrix~$\Gamma$, which does not depend on $\kappa$
or $\alpha$, is not diagonal.

\begin{figure}

\includegraphics{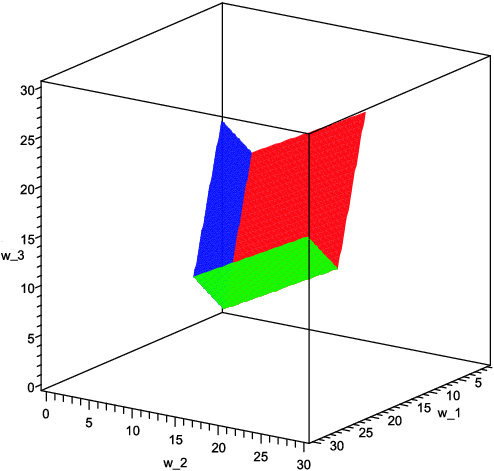}%
\vspace{-3pt}
\caption{A portion of the workload cone  $\mathcal W_{1}$ is shown for a linear network
with three resources and four routes with $\alpha=1$ and
$\nu_i=\mu_i=\kappa_i=1$ for all $i\in
\mathbb{I}$.}\label{fig:alpha1side}\vspace{-3pt}
\end{figure}

\begin{figure}[b]
\vspace{-3pt}
\includegraphics{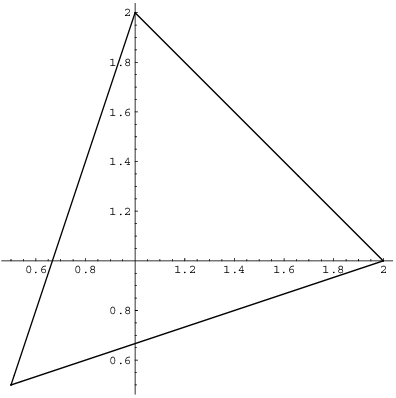}%
\vspace{-3pt}
\caption{A cross-section of the workload cone $\mathcal W_{1}$ depicted
in Figure \protect\ref{fig:alpha1side} taken at $w_3=1$.} \label{fig:alpha1cross}
\end{figure}

\begin{figure}[t]

\includegraphics{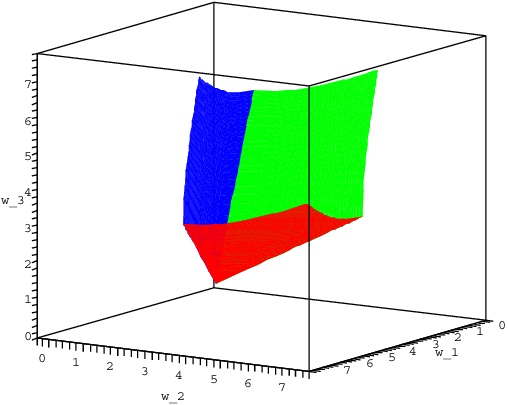}%
\vspace{-3pt}
\caption{A portion of the workload cone $\mathcal W_{2}$ is shown for a linear network
with three resources and four routes with $\alpha=2$ and $\nu_i=\mu_i=\kappa_i=1$ for
all $i\in \mathbb{I}$.}\label{fig:alpha2side}\vspace{-3pt}
\end{figure}

\begin{figure}[b]
\vspace{-3pt}
\includegraphics{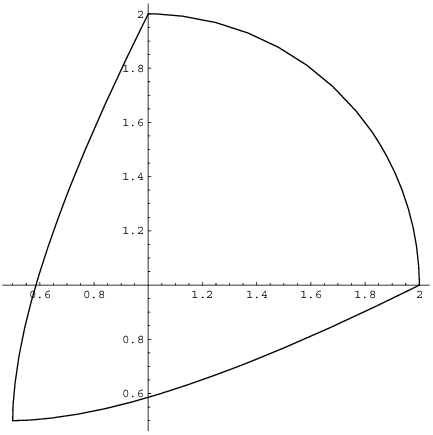}%
\vspace{-3pt}
\caption{A cross-section of the workload cone $\mathcal W_2$ depicted
in Figure \protect\ref{fig:alpha2side} taken at $w_3=1$.} \label{fig:alpha2cross}
\end{figure}

\begin{figure}[t]

\includegraphics{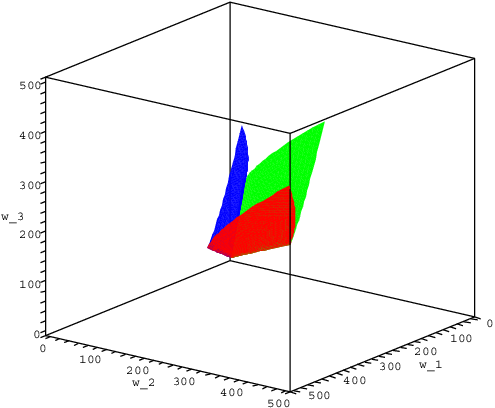}%
\vspace{-3pt}
\caption{A portion of the workload cone $\mathcal W_{0.5}$ is shown for a linear network
with three resources and four routes with $\alpha=0.5$ and $\nu_i=\mu_i=\kappa_i=1$ for
all $i\in \mathbb{I}$.} \label{fig:alpha5side}\vspace{-3pt}
\end{figure}

\begin{figure}[b]
\vspace{-3pt}
\includegraphics{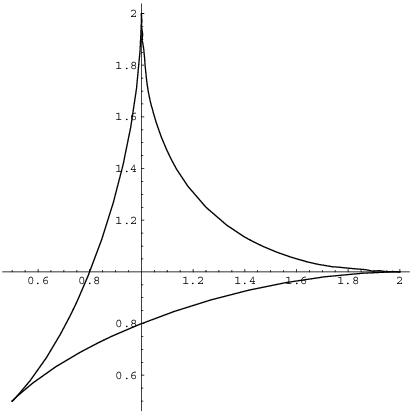}%
\vspace{-3pt}
  \caption{A cross-section of the workload cone
  $\mathcal W_{0.5}$ depicted in Figure \protect\ref{fig:alpha5side} taken at $w_3=1$.} \label{fig:alpha5cross}
\end{figure}

In general, for $\alpha\not=1$ with higher-dimensional workloads,
that is, ${\bf J}>2$, the
shape of the workload cone $\mathcal W_\alpha$ depends on $\alpha$, as
well as on $A$, $\nu$, $\mu$ and~$\kappa$.
This relationship appears to be quite complicated.
We are investigating Conjecture \ref{conjecture} in this case.
At this time, we only have some partial results.
The main difficulty in establishing the validity of the conjecture
concerns proving
$C$-tightness of $\{\hat W^r\}$ and
establishing uniqueness in law for any weak limit point.
To illustrate some of the difficulties involved,
we consider a
linear network with three resources and four routes.
Although this is a~particular example, it exemplifies
the challenges presented by the geometric data in
establishing existence, uniqueness and an invariance principle for the SRBM.
For this example, each of the three
resources has a route that only goes through that resource and there
is one additional route that goes through all three resources. Here,
the workload cone
is contained in the three-dimensional positive orthant.
When $\nu_i=\mu_i=\kappa_i=1, $ for all $i\in\mathbb{I}$, we
depict the associated workload cone and one of its cross-sections, for
$\alpha=2$, $\alpha=1$ and
$\alpha=0.5,$ in Figures \ref{fig:alpha1side}--\ref{fig:alpha5cross}.

%
%

For
the case depicted in Figures \ref{fig:alpha1side}--\ref{fig:alpha1cross},
$\alpha=1$ and
Theorem \ref{thm:HTL} applies~to justify the diffusion approximation.
The proof of this
theorem uses the uniqueness\vadjust{\goodbreak} of the diffusion
process \cite{DW} and the invariance principle developed in~\cite{KaWi}.
For the case depicted in Figures \ref{fig:alpha2side}--\ref{fig:alpha2cross},
$\alpha\!=\!2$ and
the workload~co\-ne has boundary faces that curve outward.
Away from the origin, the boundary and directions of reflection
for the proposed diffusion approximation locally satisfy conditions
required by the invariance principle given in~\cite{KaWi}. However,
the workload cone has a
``singular point''
at the origin where the conditions in \cite{KaWi} fail to be
satisfied. However, since this is an isolated point,
the invariance principle in \cite{KaWi} and the proof of
uniqueness starting in \cite{DW} can be adapted
to establish the SRBM diffusion approximation in this case.
For higher-dimensional analogs of this case,
we anticipate that a valid diffusion approximation can be established.
However,
as the dimension increases, it becomes more difficult
to compute the inward normals to all boundary faces and, as yet, we
do not have a systematic way to
verify geometric sufficient
conditions for uniqueness in law and a valid invariance principle.
For the case depicted\vadjust{\goodbreak} in~Figures \ref{fig:alpha5side}--\ref{fig:alpha5cross},
where $\alpha=0.5$, the situation
is more complex. Here,
the workload
cone has boundary faces that curve inward and any two faces
meet in a cusp, that is, the inward normals to any two boundary
faces point toward one another at the intersection of the faces.
At present, we do not have an existence and uniqueness result, nor
an invariance principle, to treat
this case because of the singular geometry at the intersections of
faces.
In particular, it has not been established that the diffusion
process can escape from the tip of a cusp or that there
is uniqueness in law for a~diffusion process starting there.

In a recent work,
Shah and Wischik \cite{SW} have proven multiplicative
state space collapse for a class of ``switched'' networks
operating under a family of scheduling policies
related to the maximum weight algorithm introduced
by Tassiulas and Ephremides \cite{TE}.
Based on this, Shah and Wischik have conjectured natural diffusion
approximations for the workload processes in these models. These
diffusions are SRBMs living in cones with piecewise smooth boundaries.
While these cones have some characteristics in common with those found
for the bandwidth sharing model considered here, there are
also some new features. Depending on a parameter
associated with the family of scheduling policies,
the cones of \cite{SW} include nonsimple convex
polyhedrons as well as cones with curved boundaries
where boundary faces can meet smoothly.
The validity of these conjectured diffusion approximations for
input-queued packet switches operating under a maximum
weight algorithm is being explored in \cite{KaWiS}.

\section{Proof of multiplicative state space
collapse}\label{section:MSSC}

In this section, we prove the multiplicative state space collapse
result, Theorem \ref{req:SSC}.
Our proof
follows a general line of argument
pioneered by Bramson
in \cite{BR}, where open multiclass queueing networks operating under
certain head-of-the-line (HL) service disciplines are treated.
However, there are some differences from \cite{BR}.
Here, we have the more general structure of
a stochastic processing network with simultaneous resource
possession
and our service discipline
is not work-conserving.
On the other hand, we have exponential interarrival and document
size distributions (rather than general distributions),
which lead to some simplifications in our proofs.
A particularly
interesting aspect here is that, in contrast
to prior results on state space collapse for open multiclass queueing networks,
our lifting map can be nonlinear (for $\alpha\not=1$).
In addition, unlike Assumption 3.1 in \cite{BR},
we do not require an
exponential rate of convergence of fluid model solutions
toward points on the invariant manifold;
we only use uniform convergence of fluid model
solutions toward the invariant manifold (starting in a compact set).
Despite these differences, our main line of argument follows
that of Bramson \cite{BR}.

We first
provide some preliminary results in Section \ref{section:preSSC}.
Our proof of multiplicative state space collapse is then given
in Section \ref{PSSC}.

We note (for the extension to multi-path routing described
in Section \ref{multipathroute})
that the proofs in this and the next section use the fact that
the entries in the matrix $A$ are nonnegative, rather
than the stronger condition that they are
zeros
or ones.

\subsection{Crucial estimates for fluid scaled processes}\label{section:preSSC}
Recall the definition of fluid scaled processes from Section
\ref{r:sequence}. For each $r>0$ and $m\in\mathbb{Z}_+$, let
%
\begin{eqnarray} \label{overlineNrm}
\overline n^r_m &=& |\overline
N^r(m)| \vee1
\end{eqnarray}
and define shifted and scaled processes $\overline
E^{r, m}, \overline S^{r, m}, \overline T^{r, m}, \overline U^{r,
m}, \overline N^{r, m}, \overline W^{r, m}$ as\vadjust{\goodbreak} follows. For
$i\in{\mathbb I}$, $ j\in{\mathbb J}$ and $t\geq0$,
%
\begin{eqnarray}
\overline E_i^{r, m} (t) &=& \frac{\overline E_i^{r}(m+ \overline n^r_m
t)-\overline E_i^r(m)}{\overline n^r_m},
\\
\overline S_i^{r, m} (t) &=&
\frac{\overline S_i^{r}(\overline T_i^r(m)+ \overline n^r_m
t)-\overline S_i^r(\overline T_i^r(m))}{\overline n^r_m},
\\\label{Tbar}
\overline T_i^{r, m} (t) &=& \frac{\overline T_i^{r}(m+ \overline n^r_m
t)-\overline T_i^r(m)}{\overline n^r_m},
\\
\overline U_j^{r, m} (t)
&=& \frac{\overline U_j^r(m+\overline n^r_mt)-\overline
U_j^r(m)}{\overline n^r_m},
\\\label{Nbar}
\overline N_i^{r, m} (t) &=& \frac{\overline N_i^{r}(m+ \overline n^r_m
t)}{\overline n^r_m},
\\
\overline W_j^{r, m} (t) &=&
\frac{\overline W_j^{r}(m+ \overline n^r_m t)}{\overline n^r_m} =
\sum_{i\in{\mathbb I}}A_{ji}(\mu_i^r)^{-1}\overline N_i^{r,m}(t).
\end{eqnarray}
The additional scaling by $\overline n^r_m $ used here is to ensure that
the starting values of $\overline N^{r,m}$ all lie in a
single compact set, namely the unit ball in $\mathbb{R}_+^{{\mathbf I}}$.
This facilitates use of the properties of fluid model solutions
described in Section~\ref{sec:fluid}.
It is easy to check using (\ref{N})--(\ref{U}),
(\ref{rfluone})--(\ref{rfluthree}) and the scaling property of
$\Lambda$
that for each $i\in\mathbb{I}$, $j\in\mathbb{J}$ and $t\geq0$,
%
\begin{eqnarray}\label{nam:N=E-ST}
\overline{N}^{r,m}_i(t) &=& \overline{N}^{r,m}_i(0) +
\overline{E}^{r,m}_i(t) - \overline{S}^{r,m}_i(\overline
{T}^{r,m}_i(t)),
\\\label{nam:U=C-AT}
\overline{U}_j^{r,m}(t)
&=& C_jt - \sum_{i\in\mathbb{I}} A_{ji} \overline{T}_i^{r,m}(t),
\\\label{nam:T=intLambda}
\overline{T}^{r,m}_i(t)
&=& \int_0^t \Lambda_i( \overline{N}^{r,m}(s) ) \,ds.
\end{eqnarray}

The following theorem summarizes essential properties of the above
processes needed for our proof of multiplicative state space collapse.
The proof of this uses arguments very similar to those in Sections 4, 5
and 6 of Bramson~\cite{BR}.

\begin{theorem} \label{fluidapprox}
Let
%
\begin{equation}\label{defK}
K= (1+|\mu|\vee|\nu|)\Bigl(1+\max_{j\in\mathbb{J}}C_j\Bigr).
\end{equation}
Fix $T>0$.
For each $L\geq1$, there exists a sequence of measurable sets
$\{\mathcal G^r_{L}\dvtx r >0\}$ and a family of positive constants
$\{r_{\varepsilon, L}\dvtx \varepsilon\in(0,1)\}$
such that for each $\varepsilon\in(0,1)$ and
$r\geq r_{\varepsilon, L}$,
\begin{longlist}
\item[(i)] $P(\mathcal G^r_{L}) \geq
1-\varepsilon$;

\item[(ii)] for each $\omega\in\mathcal G^r_{L}$ and
$m =0, 1, \ldots, \lfloor r T\rfloor$,
%
\begin{eqnarray} \label{nam:mainclaimone}
| \overline N^{r,m} (t, \omega) - \overline N^{r,m}
(s,\omega)| \leq K|t-s|+\varepsilon \qquad \mbox{for all } s, t
\in[0,L]
\end{eqnarray}
and there is a fluid model solution $\tilde
n( \cdot )$ satisfying
%
\begin{equation} \label{nam:mainclaimtwo}
\| \overline N^{r,m} ( \cdot  , \omega) - \tilde n( \cdot ) \|
_L < \varepsilon.
\end{equation}
\end{longlist}
\end{theorem}

Note that we have not indicated explicit dependence on $T$ in the above
as $T$ will always be fixed in the application of this result. Also, to simplify
notation, we have omitted explicit indication of the dependence of
$\tilde n( \cdot )$ on~$r$, $m$, $\omega$,
$L$ and $\varepsilon$.

Before proving this theorem, we establish the following preliminary lemma.
For this lemma and the proof of Theorem \ref{fluidapprox},
let $\tilde C=\max_{j\in\mathbb{J}}C_j$.

\begin{lemma}\label{nam:exists:delta}
Fix $T>0$. For each $L\geq1$, there exists a sequence of measurable
sets $\{ \mathcal G^r_L \dvtx r>0\}$ and a family of positive constants
$\{r_{\varepsilon,L}\dvtx \varepsilon\in(0,1)\}$ such that for each
$\varepsilon\in(0,1)$ and
$r\geq r_{\varepsilon,L}$, we have
%
\begin{equation}
P(\mathcal G^r_L ) \geq1-\varepsilon
\end{equation}
and on $\mathcal G^r_L$, for $m=0, 1, \ldots
, \lfloor rT \rfloor$,
%
\begin{eqnarray}
\label{nam:key epsilon approximation E}
\| \overline E^{r,m}(\cdot) -\nu^r(\cdot)\|_L &\leq& \varepsilon/4,
\\
\label{nam:key epsilon approximation S} \|\overline
S^{r,m}(\overline T^{r,m}(\cdot))-\mu^r(\overline
T^{r,m}(\cdot))\|_L &\leq& \varepsilon/4 ,
\end{eqnarray}
where $\nu^r(t) =\nu^r t$,
$(\overline S^{r,m}(\overline T^{r,m}(t)))_i =\overline
S^{r,m}_i(\overline T^{r,m}_i(t))$
and $(\mu^r(\overline T^{r,m}(t)))_i=\break \mu^r_i\overline T^{r,m}_i(t)$,
for all
$t\geq0$ and $i\in\mathbb{I}$.
\end{lemma}

\begin{pf} Fix $L\ge1$. Note that since the bandwidth allocations given
by~%
$\Lambda(\cdot)$
are bounded by~$\tilde C$, for each $r>0$ and $i\in\mathbb{I}$,
$T_i^r(\cdot)$ is Lipschitz continuous with Lipschitz constant
$\tilde C$. Since this property is unchanged by the (fluid) scaling in
(\ref{Tbar}), we have that $\overline T_i^{r,m}(\cdot)$ is also Lipschitz
continuous with the same Lipschitz constant.
It follows from this that
$\Vert \overline T^{r,m}(\cdot)\Vert_L\leq{\mathbf I}\tilde C L. $
On combining this with Assumption~\ref{HT}, the fact that the interarrival
times and document sizes are exponential and the memoryless property of
the associated
Poisson processes,
by an argument\vadjust{\goodbreak} similar to that used in verifying~(5.19) of
Proposition 5.1 in Bramson \cite{BR}, we have that
as $r\rightarrow
\infty$,
%
\begin{eqnarray} \label{claimS}
\max_{m=0}^{\lfloor rT \rfloor}\bigl( \|
\overline{S}^{r,m}(\overline{T}^{r,m}(\cdot))
-
\mu^r (\overline{T}^{r,m}(\cdot))
\|_L
\vee
\|
\overline{E}^{r,m}(\cdot)
-
\nu^r(\cdot)
\|_L\bigr)
\Rightarrow0.
\end{eqnarray}
It follows that there exists
a sequence $\{a_\ell\}_{\ell=0}^\infty$ satisfying $a_0 = 0$,
$a_\ell\in(\ell\vee a_{\ell-1}, \infty)$ for each integer $\ell
\geq1$, such that
for each $r\geq a_\ell$,
\begin{eqnarray*}\label{nam:sectiontwo:constructionone}
&&P\Biggl[ \bigcap_{m=0}^{\lfloor rT \rfloor} \biggl\{\|
\overline{S}^{r,m}(\overline{T}^{r,m}(\cdot)) - \mu^r (\overline
{T}^{r,m}(\cdot)) \|_L \vee \| \overline{E}^{r,m}(\cdot
) - \nu^r(\cdot) \|_L \le\frac{1}{2^{\ell+2}} \biggr\}
\Biggr]
\\
&&\quad \ge1 - \frac{1}{2^{\ell}}.
\end{eqnarray*}
For $r>0$, define
\[
\delta(r) = \sum_{\ell=0}^\infty\frac{1}{2^\ell} \texttt{\bf
1\rm}_{[a_{\ell}, a_{\ell+1})}(r)
\]
and
\[
\mathcal G_L^r=
\bigcap_{m=0}^{\lfloor rT \rfloor}
\biggl\{\|
\overline{S}^{r,m}(\overline{T}^{r,m}(\cdot))
-
\mu^r (\overline{T}^{r,m}(\cdot))
\|_L
\vee
\|
\overline{E}^{r,m}(\cdot)
-
\nu^r(\cdot)
\|_L
\le\frac{\delta(r)}{4}
\biggr\}.
\]
Then, since $\delta(r)\to0$ as $r\to\infty$,
for each $\varepsilon\in(0,1)$, there exists $r_{\varepsilon, L}>0$ such
that $\delta(r) \leq\varepsilon$ for all $r\geq r_{\varepsilon, L}$.
Combining the above, we conclude that for all $r\geq r_{\varepsilon, L}$,
%
\begin{equation}
\label{nam: buiding up a delta function}
P(\mathcal{G}_{L}^r) \geq1-\delta(r) \geq1-\varepsilon .
\end{equation}
\upqed\end{pf}

\begin{pf*}{Proof of Theorem \protect\ref{fluidapprox}} Fix $T>0$
and $L\geq1$. Also, fix a sequence of measurable sets $\{ \mathcal
G^r_L : r>0\}$ and a family of positive constants
$\{r_{\varepsilon,L}: \varepsilon\in(0,1)\}$, as in Lemma
\ref{nam:exists:delta}. Without loss of generality,
by Assumption \ref{HT}, we may further assume (by choosing
$r_{\varepsilon, L}$ slightly larger if necessary) that for all $r\geq
r_{\varepsilon , L}$,
\[
|\nu^r| \leq1+| \nu|\quad\mbox{and} \quad|\mu^r |\leq1+ |\mu|.
\]
By Lemma \ref{nam:exists:delta}, the uniform Lipschitz continuity of
$\overline T^{r,m}(\cdot)$
and
(\ref{nam:N=E-ST}), for each $\varepsilon\in(0,1)$ and $r\geq
r_{\varepsilon,L}$, we have that on $\mathcal{G}^r_L$, for each $s,t
\in[0,L]$ and $m =0,1, \ldots,\lfloor rT \rfloor$,
%
\begin{eqnarray}\label{nam:key epsilon approximation N}
|\overline{E}^{r,m}(t) - \overline{E}^{r,m}(s)| &\le&
\frac{ \varepsilon}{2} + |\nu^r| | t -s |,
\\\label{eq94}
|\overline{S}^{r,m}(\overline{T}^{r,m}(t))
-
\overline{S}^{r,m}(\overline{T}^{r,m}(s)) |
&\le&
\frac{\varepsilon}{2} +| \mu^r|   \tilde C   | t-s|,
\\\label{eq95}
|\overline{N}^{r,m}(t) - \overline{N}^{r,m}(s)|
&\le& \varepsilon+ (|\nu^r| \vee| \mu^r|)
(1 + \tilde C ) |t-s|.
\end{eqnarray}
Then, part (i) of Theorem \ref{fluidapprox} follows directly
from Lemma \ref{nam:exists:delta}.
Inequality~(\ref{nam:mainclaimone}) follows from
Lemma\vadjust{\goodbreak}
\ref{nam:exists:delta},
(\ref{eq95}) and the choice of
$r_{\varepsilon,L}$ and $K$.
As for the last part, (\ref{nam:mainclaimtwo}), of Theorem \ref{fluidapprox},
for a proof by
contradiction, suppose that there exists an $\varepsilon_0\in(0,1)$ such
that for each integer $\ell\geq r_{\varepsilon_0,L}$, there exists a~%
value $r_\ell$ of $r$
such that $r_{\ell} \ge\ell$, and such that there exists $m_\ell\in
\{0,1,\ldots,\lfloor r_\ell T\rfloor\}$ and $\omega_\ell\in
\mathcal{G}
^{r_\ell}_L$ so that for
any fluid model solution $\overline{n}$ we have
%
\begin{equation}\label{nam:contradiction}
\| \overline{N}^{r_\ell,m_\ell}(\cdot,\omega_\ell) - \overline
{n}(\cdot) \|_L
\ge\varepsilon_0.
\end{equation}
For the contradiction, we will show that there exists a subsequence
of\break
$\{\overline{N}^{r_\ell,m_\ell}(\cdot, \omega_\ell)\}_{\ell
=1}^\infty$ that
converges uniformly on $[0,L]$ to a fluid model solution. For this,
note that
$\{\overline{N}^{r_\ell,m_\ell}(0,\omega_\ell)\}_{\ell=1}^\infty
$ is bounded by one, and\break
$\{ \overline{T}^{r_\ell,m_\ell}(\cdot, \omega_\ell)\}_{\ell
=1}^\infty$ is
uniformly bounded and equicontinuous on $[0,L]$, by the uniform
Lipschitz continuity of
$\overline{T}^{r,m}$.
Thus, using the Ascoli--Arzel\`{a} theorem, we have that along a subsequence,
$\{\overline N^{r_\ell,m_\ell}(0,\omega_\ell)\}_{\ell=1}^\infty$
converges to a finite value and $\{\overline T^{r_\ell, m_\ell}
(\cdot, \omega_\ell)\}_{\ell=1}^\infty$ converges uniformly on
$[0, L]$.
Fix such a convergent subsequence and denote the respective limits by
$\overline n(0)$ and $ \tau(\cdot)$.
Then, using the convergence of
$\{(\overline E^{r_\ell,m_\ell}(\cdot, \omega_\ell), \overline
S^{r_\ell,m_\ell}(\overline T^{r_\ell,m_\ell}(\cdot, \omega_\ell
)\omega_\ell))\}_{\ell=1}^\infty$
implied by~(\ref{nam:key epsilon approximation E}) and (\ref{nam:key epsilon
approximation S}), we have that along the same subsequence,
$\{\overline N^{r_\ell,m_\ell}(\cdot, \omega_\ell)\}_{\ell
=1}^\infty$ converges
uniformly on $[0, L]$ to $\overline n(\cdot)$, defined by
%
\begin{eqnarray}
\overline{n}_i(t) = \overline{n}_i(0) + \nu_i t - \mu_i
{\tau}_i(t), \qquad t\in[0, L],  i\in\mathbb{I}.
\end{eqnarray}
Furthermore,
by following the pathwise argument in Appendix B of Kelly and Williams
\cite{KW}, one can show that at each $t \in(0, L]$ where the
derivative of
$\overline n(\cdot)$ exists (this means the left-hand derivative at $t=L$),
we have that~(\ref{eq:diff}) and (\ref{eq:constraint}) hold with
$\overline n$ in place of $n$ there.
It follows that $\overline n(\cdot)$ has the same properties as a fluid
model solution on $[0, L]$.
We can easily extend $\overline n$ to be a fluid model solution on $[0,
\infty)$ by
defining it on $(L, \infty)$ to equal
$
\overline{f}(\cdot-L),
$
where
$\overline{f}\dvtx [0,\infty) \rightarrow\mathbb{R}_+^{\mathbf{I}}$ is a
fluid model solution satisfying
$\overline{f}(0) = \overline{n}(L)$. (The existence of such
a fluid model solution follows from Theorem B.1 in \cite{KW}.)
This yields the desired contradiction of
(\ref{nam:contradiction}).
\end{pf*}

\subsection{Proof of multiplicative state space collapse}\label{PSSC}

We will use Theorem~\ref{fluidapprox} to show that
for $T>0$ fixed, for each $\overline L\geq
1$, there exists $L> \overline L$ such that for each $\varepsilon
\in(0, 1)$, for all $r$ sufficiently large, with probability at least
$1-\varepsilon$, the left member of (\ref{conveq}) is dominated by
%
\begin{eqnarray}\label{withL}
\|\overline
N^{r, 0}( \cdot  ) - \Delta(\overline W^{r, 0}( \cdot
))\|
_{[0, \overline L]} + \sup\limits^{\lfloor rT\rfloor-1}_{m=0}
\|\overline
N^{r, m}( \cdot ) - \Delta(\overline W^{r, m}( \cdot
))\|_{[\overline L,
L]},
\end{eqnarray}
where the
quantities
$\{ \overline N^{r,m}, m=0,1, \ldots, \lfloor rT\rfloor-1\}$ can be
approximated by fluid model solutions over $[0,L]$.
The
results of
Section \ref{sec:fluid} on the behavior of fluid model solutions, especially
Proposition \ref{thm:convinv}, can then be combined with
a suitable choice of $\overline L$ and the assumptions of Theorem \ref{req:SSC}
on the initial conditions to prove multiplicative state space collapse.
We now give the detailed proof.\vadjust{\goodbreak}

\begin{pf*}{Proof of Theorem \protect\ref{req:SSC}}
Fix $T>0$. For each $r >0$, let $ \mathcal H^r_T$ denote the left
member of (\ref{conveq}).
Since,
for each $r>0$ and
$t\geq0$, we have
%
\begin{eqnarray}\label{scaleNW}
\hat N^r(t) &=& \overline N^r(rt), \qquad\hat W^r(t) =
\overline
W^r(rt)= A(M^r)^{-1} \overline N^r(rt),
\end{eqnarray}
it follows that
%
\begin{eqnarray}\label{eq:1}
\mathcal H^r_T &=& \frac{ \|\overline N^r( \cdot
)-\Delta(\overline W^r ( \cdot  ))\|_{rT} }
{\|\overline N^r( \cdot )\|_{rT} \vee1}.
\end{eqnarray}
For $\overline L \ge1$ fixed, intervals of the form $[m,
m+\overline L  \overline n^r_m], m=0, 1, \dots, \lfloor rT\rfloor$,
cover $[0, rT]$ since $\overline n^r_m\geq1$
for each $m\in\mathbb{Z}_+$. Hence, for each $r>0$ and $t\in[0,r
T]$, the
following is well defined as a random variable taking values in $\{0,
1, \dots, \lfloor rT\rfloor\}$:
%
\begin{equation}\label{eq:2}
m^r(t) = \inf\{m\dvtx m\le t \le m+\overline L  \overline n^r_m\}.
\end{equation}
For any $r>0$ and $t\in[0, rT]$ such that $t > \overline L\overline
n^r_0$, we have $m^r(t) \ge1$ and then
%
\begin{equation}\label{eq:3}
t-\bigl(m^r(t)-1\bigr) > \overline L  \overline n^r_{m^r(t)-1}
\end{equation}
since, otherwise, $m^r(t)$ could be replaced by a value of $m$ strictly
less\break then~$m^r(t)$. Thus, for each $r>0$ and $t\in[0, rT]$ such that
$t>\overline L  \overline n^r_0$, we have
%
\begin{equation}\label{eq:4}
m^r(t)-1+ \overline L  \overline n^r_{m^r(t)-1} < t \le m^r(t) +
\overline L  \overline n^r_{m^r(t)},
\end{equation}
where $m^r(t) \in\{1, 2, \dots, \lfloor rT \rfloor\}$.
From (\ref{eq:1}) and (\ref{eq:4}), we have that $ \mathcal H^r_T$ is
dominated by
%
\begin{eqnarray} \label{rHHeqn}
\quad && \frac{\|\overline N^r( \cdot ) - \Delta(\overline
W^r(  \cdot  ))\|_{[0,\overline L  \overline
n^r_0]}} {\overline n^r_0} + \sup\limits^{\lfloor rT\rfloor}_{m=1}
\frac{\|\overline N^r( \cdot  ) - \Delta(\overline
W^r( \cdot  ))\| _{[m-1+ \overline L  \overline
n^r_{m-1}, m+\overline L  \overline n^r_m]}} {\overline
n^r_{m-1} }
\nonumber\\
&&\qquad = \|\overline N^{r, 0}( \cdot ) - \Delta(\overline
W^{r, 0}( \cdot ))\| _{[0, \overline L]}
\\
&&{}\qquad\quad +
\sup\limits^{\lfloor rT\rfloor}_{m=1}
\|\overline N^{r,
m-1}(\cdot) - \Delta(\overline W^{r, m-1}(\cdot
))\| _{[\overline L, {(1+\overline L
\overline n^r_m)/\overline n^r_{m-1}}]}, \nonumber
\end{eqnarray}
where we have used the scaling property of $\Delta(\cdot)$
in obtaining the equality (see Proposition \ref{scaling}).

Recall the definition of the constant $K\geq1$ from (\ref{defK}).
Given $\overline L\geq1$, let $L=4K\overline L$.
Below, we refer to the sequence of measurable sets $\{ \mathcal
G^r_L\dvtx
r>0\}$ introduced in
Theorem \ref{fluidapprox}.
Focusing on the right endpoints of the time intervals
appearing in the last term of (\ref{rHHeqn}),
we next show that, for all $r$ sufficiently large,
on $\mathcal G^r_L$,
for $m=1,\ldots, \lfloor rT\rfloor$, we have
%
\begin{eqnarray} \label{Lest}
\frac{1+\overline L \overline n^r_m}{\overline n^r_{m-1}}\leq L.
\end{eqnarray}
For this,
note that, by (\ref{Nbar}), we have, for $m=1, 2, \ldots,$
%
\begin{eqnarray} \label{rNbar}
\frac{\overline N^r(m)}{ \overline n^r_{m-1}}=\frac{\overline
N^r(m-1+1)}{\overline n^r_{m-1}}=\overline
N^{r,m-1}\biggl(\frac{1}{\overline n^r_{m-1}}\biggr).
\end{eqnarray}
On the other hand, by (\ref{overlineNrm}), for $m=1, 2, \ldots, $
%
\begin{equation} \label{rNbarzero}
|\overline
N^{r,m-1}(0)|=\biggl|\frac{\overline N^{r}(m-1)}{\overline
n^r_{m-1}}\biggr|\leq1.
\end{equation}
It follows from
(\ref{rNbar}), (\ref{rNbarzero}) and Theorem \ref{fluidapprox}
that for each $\varepsilon\in(0, 1)$, for each $r\geq r_{\varepsilon
, L}$
and $m =1, \ldots, \lfloor rT\rfloor$, on $ \mathcal G^r_{L}$, we have
%
\begin{eqnarray} \label{growthin}
\frac{\overline n^r_m}{\overline n^r_{m-1}} &=& \frac{|\overline
N^r(m)|}{\overline n^r_{m-1}}\vee\frac{1}{\overline n^r_{m-1}}
\nonumber
\\
&\leq& \biggl| \overline N^{r,m-1}\biggl(\frac{1}{\overline
n^r_{m-1}}\biggr)\biggr| \vee1
\\
&\leq& \biggl(| \overline N^{r, m-1}(0)| + K\biggl(\frac{1}{\overline
n^r_{m-1}}\biggr) +\varepsilon\biggr)\vee1\nonumber
\\
&\leq& (1 + K\cdot1 +1)\vee1\nonumber
\\
&\leq& 3 K \nonumber
\end{eqnarray}
and hence (\ref{Lest}) holds.
On combining (\ref{rHHeqn}) with (\ref{Lest}), we have that
for each $\varepsilon\in(0,1)$, for all $r\geq r_{\varepsilon, L}$, the
following inequality holds on $\mathcal G^r_{L}$:
\begin{eqnarray} \label{simpleH}
{\mathcal H}^r_T
&\leq& \|\overline N^{r, 0}( \cdot  ) - \Delta
(\overline W^{r, 0}( \cdot )) \| _{[0, \overline
L]}\nonumber
\\[-8pt]\\[-8pt]
&&{}+ \sup\limits^{\lfloor rT\rfloor-1}_{m=0}\
\|\overline N^{r, m}( \cdot ) - \Delta(\overline W^{r,
m}( \cdot  ))\| _{[\overline L,L]}. \nonumber
\end{eqnarray}

We will now estimate the two terms on the right-hand side of (\ref{simpleH}).
In brief, the idea is to use Theorem \ref{fluidapprox} to show that
for all $r $ sufficiently large, for each $\omega\in\mathcal G^r_L$,
$m=0,1, \ldots, \lfloor rT\rfloor$, there is
a fluid model solution $\tilde n$ that is uniformly close to $\overline
N^{r,m}(\cdot, \omega)$
over the time interval $[0, L]$.
For the first term in (\ref{simpleH}), we then use the assumptions on the
initial conditions and Proposition~\ref{closetoMalpha} to show that for all $r$ sufficiently large,
with high probability, $\tilde n$ is uniformly close to the invariant
manifold $\mathcal M_\alpha$
on $[0,\infty)$ and use this to control the first term in (\ref{simpleH}).
For the last term in (\ref{simpleH}), we use Proposition~\ref{thm:convinv}
to choose $\overline L$ (independent of $r, \omega$ or $m$)
so that $\tilde n$ is uniformly close to the
invariant manifold $\mathcal M_\alpha$ on $[\overline L,\infty)$
and so, for all $r$ sufficiently large, $\omega\in\mathcal G^r_L$ and
$m =0, 1, \ldots, \lfloor rT\rfloor-1$, for each $t\in[\overline L, L]$,
there exists a~point~$\tilde n^*_t$ (depending on $r, \omega, m, t$)
in $\mathcal M_\alpha$
that is uniformly close to~$\overline N^{r,m}(t, \omega)$.
Then,
\begin{eqnarray}
&& | \overline N^{r,m} (t,\omega) -\Delta(\overline W^{r,m} (t,\omega
)) | \nonumber
\\
&&\qquad \leq|\overline N^{r,m} (t,\omega) - \tilde n^*_t| + |\tilde
n^*_t- \Delta(AM^{-1}\tilde n^*_t)| \nonumber
\\
&&\qquad\quad{}+ |\Delta
(AM^{-1} \tilde n^*_t)
- \Delta(A(M^r)^{-1} \overline N^{r,m} (t,\omega)) |,
\nonumber
\end{eqnarray}
where each of the last three terms can be made small (uniformly for
$\omega\in\mathcal G^r_L$, $m=0, 1, \ldots, \lfloor rT\rfloor-1, t\in
[\overline L, L]$,
for all $r$\vadjust{\goodbreak} sufficiently large), by the choice of~$\tilde n_t^*$,
the fact that $n=\Delta(w(n))$ for a point $n$ on the invariant
manifold~$\mathcal M_\alpha$,
the continuity of $\Delta$ and the convergence of $(M^r)^{-1} $ to
$M^{-1}$ as $r\to\infty$.
We now give the full details of the argument. In the first paragraph below,
we develop estimates that will be used for both the first and the
second term on the
right-hand side of
(\ref{simpleH}). (Accordingly, we consider values of
$m$ that include $\lfloor rT\rfloor$, as this value may be zero.)

Fix $\eta\in(0,1)$. [In this proof only, we locally reuse the symbol
$\eta$
for a~positive constant in $(0,1)$. This is distinct
from the use of $\eta$ elsewhere as a probability measure.]
By Proposition~\ref{compactfluid},
there is a constant $D(9/8)\geq9/8$ such that any fluid model solution
$n( \cdot )$ satisfying $|n(0)|\leq9/8$ satisfies $|n(t) |
\leq D(9/8)$ for all $t\geq0$.
By the
uniform continuity of $\Delta(  \cdot )$ on compact sets in~$\mathbb{R}_+^{\bf J}$ (see Proposition \ref{scaling}) and the fact that
$(M^r)^{-1} \to M^{-1}$ as $r\to\infty$ (see Assumption \ref{HT}),
there are constants $r_\eta>0$ and $\gamma
\in(0, \eta/4)$ such that for all $r\geq r_\eta$,
%
\begin{equation} \label{deltacont}
|\Delta(AM^{-1} n ) -\Delta
(A(M^r)^{-1} n' )| <\frac{\eta}{8}
\end{equation}
whenever $n, n'\in\mathbb{R}_+^{\bf I}$, $|n|\vee|n'|\leq
1+D(9/8)$ and $|n-n'|< 2 \gamma$. By Proposition~\ref{closetoMalpha}, there exists $\varepsilon\in(0,\gamma/2) $
such that if
$n( \cdot )$ is a fluid model solution satisfying $|n(0)|\leq9/8$ and
$d(n(0),\mathcal{M}_{\alpha})<2\varepsilon$, then
%
\begin{equation}\label{zeroclose}
d(n(t), \mathcal M_\alpha) <\gamma\qquad\mbox{for all } t\geq0.
\end{equation}

Let $\overline L = T_{{9/8}, \varepsilon}$, where
$T_{R,\varepsilon}$
is defined in Proposition \ref{thm:convinv}.
Set
$L=4K\overline L$ as above. Then, by Theorem \ref{fluidapprox}, for
each $r\geq r_{\varepsilon, L}$, $m=0,1, \ldots, \lfloor rT\rfloor$ and
$\omega\in\mathcal G^r_{L}$, there exists a fluid model solution
$\tilde
n( \cdot )$ such that
%
\begin{equation} \label{closetom}
\| \overline N^{r,m} ( \cdot  , \omega) - \tilde n( \cdot ) \|
_L < \varepsilon.
\end{equation}
We fix such $r$, $m$ and $\omega$ for the remainder of this
paragraph.
Since $|\overline N^{r,m} (0,\allowbreak\omega)|\leq1$ and $\varepsilon<1/8$,
we have
$|\tilde n(0)|< 9/8 $. It then follows from Proposition~\ref{thm:convinv} and the choice of $\overline L$ that
%
\begin{equation}
d(\tilde n(t) , \mathcal M_\alpha) < \varepsilon \qquad  \mbox{for all }
t\geq\overline L.
\end{equation}
For each $t\in[\overline L, L]$, let $\tilde n^*_t\in\mathcal
M_{\alpha}$ such that
%
\begin{equation} \label{rnnstar}
|\tilde n(t) -\tilde n^*_t| <\varepsilon.
\end{equation}
On
combining this with (\ref{closetom}), we obtain, for each $t\in
[\overline L, L]$,
%
\begin{equation} \label{nclosetonstar}
|\overline N^{r,m} (t, \omega) - \tilde n^*_t| < 2\varepsilon< \gamma.
\end{equation}
By Proposition \ref{compactfluid} and the fact that $|\tilde
n(0)|\leq9/8$, we have $|\tilde n(t)|\leq D(9/8)$ for all $t\geq
0$. So, by (\ref{rnnstar}) and (\ref{closetom}), for each $t\in
[\overline L, L]$,
\[
|\tilde n_t^*|\leq\varepsilon+D(9/8)<1+D(9/8)
\]
and
\[
|\overline N^{r,m}
(t, \omega)|\leq\varepsilon+D(9/8)<1+D(9/8).\vadjust{\goodbreak}
\]

On combining the above and using
the fact that $\overline W^{r,m} = A(M^r)^{-1}\overline
N^{r,m}$, we have, for each $r\geq r_{\varepsilon, L}\vee r_\eta, $
$m =0, 1,
\ldots, \lfloor rT\rfloor-1,$ $\omega\in\mathcal G^r_L$ and $t\in
[\overline L, L]$,
\begin{eqnarray} \label{est:barNW}
&& | \overline N^{r,m} (t,\omega) -\Delta(\overline W^{r,m} (t,\omega
)) | \nonumber
\\
&&\quad \leq|\overline N^{r,m} (t,\omega) - \tilde n^*_t| + |\tilde
n^*_t- \Delta(AM^{-1}\tilde n^*_t)|\nonumber
\\[-8pt]\\[-8pt]
&&{}\qquad+ |\Delta
(AM^{-1} \tilde n^*_t)-
\Delta(A(M^r)^{-1} \overline N^{r,m} (t,\omega)) |
\nonumber
\\
&&\quad \leq 2\varepsilon+0+\frac{\eta}{8} <
\frac{\eta}{2}. \nonumber
\end{eqnarray}
For the second inequality, we have used
(\ref{nclosetonstar}), part (iv) of Theorem \ref{invariant} and~(\ref
{deltacont}).
This takes care of estimating the last term in (\ref{simpleH}).

To estimate
the first term on the right-hand side of
(\ref{simpleH}), we need to use the initial behavior of $\overline N^r$.
By the assumption in the
theorem,
%
\begin{equation} \label{convtozeroi}
|\overline N^{r,0}(0) -\Delta(\overline
W^{r,0}(0))|=|\hat N^r(0) -\Delta(\hat
W^r(0))| \to0
\end{equation}
in probability as $ r\to
\infty.$
By Proposition \ref{Deltaw}, we have $\Delta(\overline
W^{r,0}(0))\in\mathcal{M}_{\alpha}$ and so it follows that
%
\begin{equation}
d(\overline
N^{r,0}(0),\mathcal{M}_{\alpha})
\rightarrow0\qquad\mbox{in probability as }
r\rightarrow\infty.
\end{equation}
Let $r_\varepsilon>0$ such that
%
\begin{equation}\label{rdeltar}
P\bigl(d(\overline N^{r,0}(0), \mathcal M_\alpha)\geq\varepsilon
\bigr)<\varepsilon
\qquad\mbox{for all } r\geq r_\varepsilon.
\end{equation}
For each $r\geq r_\varepsilon\vee r_{\varepsilon,L}$ and $\omega\in
\mathcal G^r_L$ satisfying
$d(\overline N^{r,0}(0,\omega),\mathcal{M}_{\alpha})< \varepsilon$, by
(\ref{closetom}), we have
%
\begin{equation}
|\tilde n(0) |\leq1 +\varepsilon\leq\tfrac{9}{8} \quad\mbox{and} \quad
d(\tilde n(0),{\mathcal M}_{\alpha})< \varepsilon+\varepsilon=
2\varepsilon.
\end{equation}
For such $r, \omega$, it
follows from (\ref{zeroclose}) that for each $t\in[0, \infty)$,
there exists $\tilde n_t^*\in\mathcal M_\alpha$ such that $|\tilde n(t)
-\tilde n^*_t| <\gamma$. Hence, by (\ref{closetom}), for $t\in
[0,\overline L]\subset[0,L]$,
we have
%
\begin{eqnarray}
|\overline N^{r,0}(t,\omega) -n^*_t |<\varepsilon+\gamma< \frac
{3\gamma}{2} .
\end{eqnarray}
Then, in a manner similar to that used in showing (\ref{est:barNW}),
we have,
for all $t\in[0,\overline L]$, $r\geq r_\varepsilon\vee
r_{\varepsilon,L}\vee r_\eta$ and $\omega\in\mathcal G^r_L$ satisfying
$d(\overline N^{r,0}(0,\omega),\mathcal{M}_{\alpha})< \varepsilon$,
%
\begin{eqnarray}\label{est:N0}
|\overline N^{r, 0}(t ,\omega) - \Delta(\overline W^{r,
0}(t,\omega))| \leq(\varepsilon+ \gamma) + 0+
\frac{\eta}{8} < \frac{\eta}{2}.
\end{eqnarray}

By combining Theorem \ref{fluidapprox}, (\ref{simpleH}), (\ref{est:barNW}),
(\ref{rdeltar}) and (\ref{est:N0}), we
have that, for all $r\geq r_\varepsilon\vee r_{\varepsilon, L}\vee
r_\eta$,
%
\begin{eqnarray}\label{rfinalmss}
P( \mathcal H^r_T \geq\eta) &\leq&
P\bigl( \{\mathcal H^r_T \geq\eta\} \cap\mathcal G^r_L\cap\{
d(\overline
N^{r,0}(0),\mathcal{M}_{\alpha})<\varepsilon\}
\bigr)\nonumber
\\
&&{}+P((\mathcal G^r_L)^c)+P\bigl(d(\overline
N^{r,0}(0),\mathcal{M}_{\alpha})\geq\varepsilon\bigr)
\\
&\leq& 0+ \varepsilon+ \varepsilon=2\varepsilon<\frac{\eta}{4}.
\nonumber
\end{eqnarray}
Since $\eta\in(0,1)$ was arbitrary, it follows that
$\mathcal H^r_T\to0$ in probability as\break $r\to\infty$.\vspace*{-3pt}
\end{pf*}

\section{\texorpdfstring{Proof of diffusion approximation when $\alpha=1$}
{Proof of diffusion approximation when alpha=1}}\label{alphaHTL}

Throughout this section, we assume that $\alpha=1$ and that the local
traffic condition, Assumption~\ref{localtraffic}, holds.
To prove the diffusion approximation result,
Theorem \ref{thm:HTL}, we shall use
the invariance principle
in \cite{KaWi}.
A key
assumption for that theory to yield convergence (rather than just
$C$-tightness) is that there is
existence and uniqueness in law for the limit diffusion process $\tilde
W$, which follows in
the case $\alpha=1$ from work of Dai and Williams \cite{DW}.
We first verify that the basic assumptions of \cite{KaWi} and \cite
{DW} are satisfied by
the state space $\mathcal{W}_1$ [see~(\ref{simpleW})]
and directions of reflection $\{ \gamma^j\dvtx j\in{\mathbb J}\}$.
For this, we need the following definition. Also,
recall the definition of the matrix~$B$
following (\ref{simpleW}).

\begin{defn}
A $d\times d$ matrix $D$ is \textit{completely-}$S$ if and only if,
for each
principal submatrix $\widetilde D$ of $D$, there exists a vector
$\tilde
x\geq0$ such that $\widetilde D \tilde x>0$. (Here, a principal
submatrix of $D$ is a matrix obtained by deleting all rows and
columns of $D$ with indices in some strict subset of $\{1,2,\ldots,d\}$.)
\end{defn}

\begin{lemma} \label{lem:coms}
The ${\bf J}\times{\bf J}$ symmetric matrix $ABA'$ is positive
definite and invertible.
The inverse matrix $(ABA')^{-1}$ is positive definite and \mbox{completely-$S$}.
\end{lemma}

\begin{pf} Since $A$ has full row rank and $B$ is strictly
positive definite, the symmetric matrix $ABA'$ is also strictly
positive definite. Indeed, for a~nonzero vector $x\in\mathbb{R}^{\bf J}$,
$x'ABA'x=(A'x)'B(A'x)$. Since $A$ has full row rank,~$A'x$ is also
nonzero and then, since $B$ is strictly positive definite,
we have that $(A'x)'B(A'x)>0$.
Since $ABA'$ is strictly positive definite, it is
invertible and the inverse
matrix $(ABA')^{-1}$
is also strictly positive definite.
It then follows that $(ABA')^{-1}$ is
a $P$-matrix, that is, all principal minors are positive, and hence
$(ABA')^{-1}$ is completely-$S$ (see \cite{CPS}, especially Theorems~3.3.7, 3.9.11 and
Corollary 3.9.13 for the relationship between $P$-matrices and
completely-$S$ matrices).
\end{pf}

The results we need to apply from \cite{KaWi}
to prove Theorem \ref{thm:HTL} require that Assumptions
(A1)--(A5)  and 5.1 of \cite{KaWi} hold. We shall not
fully describe these general conditions.
However, we shall now indicate
the meaning of these assumptions for the context treated here.
Assumptions (A1)--(A3)  are restrictions on the state space
and Assumptions (A4)--(A5)  pertain to the directions of
reflection, all
for the limit diffusion.
As noted in \cite{KaWi}, Assumptions~(A1)--(A3)  are
satisfied by a convex polyhedron with nonempty interior that
is described as the intersection of a minimal
set of half-spaces. As we shall see in the next lemma, the state space
$\mathcal W_1$ has this convex polyhedral form.
Assumption (A4)  requires that the reflection vector field on
each
boundary face be uniformly Lipschitz continuous and of unit length.
Since\vadjust{\goodbreak} our reflection vectors are constant unit vectors on each boundary face,
this condition is trivially satisfied in our context.
Assumption (A5)  imposes geometric conditions on the directions
of reflection.
These conditions are generalizations for domains with piecewise smooth
boundaries
of conditions identified earlier
by Dai and Williams \cite{DW}
for existence and uniqueness of SRBMs living in convex
polyhedral domains with a constant direction of
reflection on each boundary face. The conditions of \cite{DW} are labeled
as Assumption 1.1 in \cite{DW} and as Assumption 5.1 in \cite{KaWi}.
(We shall
use the latter label here.)
For simple convex polyhedral cones (as we have here with $\mathcal W_1$),
this assumption can be expressed
as a completely-$S$ condition on a suitable matrix formed using the directions
of reflection and the normals to the boundary faces.
We will show that this condition is satisfied in our context,
and, as a consequence, Assumption (A5)  of \cite{KaWi} will follow
immediately. We now formally verify that the aforementioned assumptions
all hold
in our context.\looseness=-1\vspace*{-2pt}

\begin{lemma} \label{dwcond}
Assumptions \textup{(A1)--(A5)}  and Assumption 5.1 of \cite{KaWi}
are all satisfied
by $\mathcal{W}_1$ and $\{\gamma^j:j\in\mathbb{J}\}$.\vspace*{-2pt}
\end{lemma}

\begin{pf}
For $ j\in\mathbb{J}$, let $n^j$ be the vector given by the $j$th row of
$(ABA')^{-1}$. Since $(ABA')^{-1}$ is symmetric, $n^j$ is also the
$j$th column
of $(ABA')^{-1}$. By~(\ref{simpleW}) and (\ref{Wsupj}), we have
%
\begin{equation}
\mathcal{W}_1= \{ w\in\mathbb{R}^{{\mathbf J}}\dvtx n^j\cdot w \geq0
\mbox{ for all } j\in
\mathbb{J}\}
\end{equation}
and
%
\begin{equation}
\quad\ \mathcal{W}_1^j=\{ABA'q\dvtx q\in\mathbb{R}_+^{\bf J}, q_j=0 \} = \{
w\in\mathcal{W}_1\dvtx
n^j\cdot
w=0\}, \qquad j\in\mathbb{J}.
\end{equation}
Since the $n^j, j\in\mathbb{J}$, are linearly independent,
it follows that $\mathcal{W}_1$ is a simple convex polyhedral cone
with minimal
representation given
by the intersection of the half-spaces $\{ w\in\mathbb{R}^{{\mathbf
J}}\dvtx n^j\cdot w
\geq0\}$,
$j\in\mathbb{J}$. It is easy to see that~$\mathcal W_1$ has
nonempty interior. As noted in \cite{KaWi}, it follows that
Assumptions~(A1)--(A3)  of \cite{KaWi} are satisfied.
Since the $\{\gamma^j, j\in\mathbb{J}\}$ are constant vectors, they
are trivially uniformly Lipschitz continuous and so Assumption
(A4)  of \cite{KaWi} holds. As noted in
Section 5.2 of \cite{KaWi}, Assumption (A5) will hold
if Assumption 5.1
of \cite{KaWi} holds.
For this, since $\mathcal{W}_1$ is a simple convex polyhedron,
it is equivalent to verify that the matrix
$NR$ is completely-$S$, where~$N$ is the ${\mathbf J}\times{\mathbf
J}$ matrix whose rows
are given by the normals $n^j, j\in\mathbb{J}$, and $R$ is the matrix whose
columns are given by the vectors $\gamma^j, j\in\mathbb{J}$. (In
fact, $N$ is
normalized so that the rows have unit length in \cite{KaWi}, but one may
equivalently use the unnormalized matrix $N$.)
Now, $N=(ABA')^{-1}$ and $R=I$, the ${\mathbf J}\times{\mathbf J}$
identity matrix,
so $NR= (ABA')^{-1}$, which we know is completely-$S$ by Lemma \ref{lem:coms}.\vspace*{-2pt}
\end{pf}

Throughout this section, we will need various constants in inequalities.
Some of these constants are denoted by $\hat C_j$, for
$j=0,1,2,\ldots.$
The hat has been added to the notation here simply to distinguish these
constants from the
bandwidth capacities $C_j$.
The next two lemmas provide basic ingredients for the proof
of Theorem~\ref{thm:HTL}.\vadjust{\goodbreak}

\begin{lemma} \label{lem:WN} There exist constants $\hat C_1>0$, $\hat
C_2\geq1$ and $\bar r>0$ such
that for each $T>0$,
%
\begin{equation}\label{wnequal}
\hat C_1
\Vert \hat{W}^r(\cdot)\Vert_T\leq
\Vert\hat{N}^r(\cdot)\Vert_T\leq\hat C_2
\Vert\hat{W}^r(\cdot)\Vert_T  \qquad \mbox{for all
}r> \bar
r.
\end{equation}
\end{lemma}

\begin{pf} By Assumption \ref{HT}, there exists $\bar r>0$ such that
$\Vert A(M^r)^{-1}\Vert\leq
2\Vert AM^{-1}\Vert$ and
$|\mu^r|<2|\mu|$ for all $r>\bar r$.
For each $t\geq0$, by (\ref{eqn:WN}), we have that
%
\begin{equation}
|\hat{W}^r(t)|\leq
\Vert A(M^r)^{-1}\Vert |\hat{N}^r(t)|.
\end{equation}
Now,
$\Vert AM^{-1}\Vert >0$ and so by letting
$\hat C_1=(2\Vert AM^{-1}\Vert)^{-1}$,
we see that the left-hand inequality in (\ref{wnequal}) holds for any $T>0$,
for all $ r>\bar
r. $
On the other hand, for each $i\in
\mathbb{I}$, there is at least one $j_i\in\mathbb{J}$ such that
$A_{j_i i}>0$ and then
using (\ref{eqn:WN}) again, we obtain that, for each $t\geq0$,
\[
A_{j_i i}(\mu^r_i)^{-1}\hat{N}^r_i(t)\leq\hat{W}_{j_i}^r(t).
\]
By letting $\hat C_2=\max(1, 2|\mu| {\mathbf I}\max_{i\in\mathbb
{I}} (A_{j_i
i})^{-1})$, we have that
the right-hand inequality in (\ref{wnequal}) holds for any $T>0$ and
$r>\bar r.$
\end{pf}

\begin{lemma} \label{lem:hatXr}
The sequence of processes
$\{\hat{W}^r(0)+\hat{X}^r(\cdot),r>0\}$ is $C$-tight.
\end{lemma}

\begin{pf} Recall that for each $t\geq0$,
%
\begin{equation}\label{rXhat}
\hat{X}^r(t)=A(M^r)^{-1}\bigl(\hat{E}^r(t)-
\hat{S}^r(\bar{\hspace{-1pt}\bar{T}}{}^r(t))\bigr)+r (A\rho^r-C)t.
\end{equation}
Since each route must use at least one resource,
for each \mbox{$n\in\mathbb{R}_+^{{\mathbf I}}$} and \mbox{$i\in\mathbb{I}$}, the
bandwidth
allocation $\Lambda_i(n) $ must be bounded by $\tilde C=\max_{j\in
\mathbb{J}} C_j$.
It follows that for each $r>0$, $\bar{\hspace{-1pt}\bar{T}}{}^r(\cdot)$ is
uniformly
Lipschitz continuous with Lipschitz constant $\tilde C$ and
hence $\{\bar{\hspace{-1pt}\bar{T}}{}^r(\cdot)\}$ is $C$-tight. On combining this with
the functional central limit assumption (\ref{clt}), the
convergence of $\{r (A\rho^r-C)t,r>0\}$ and
$\{M^r,r>0\}$ and the convergence assumption
on $\hat W^r(0)$ made at the end of Section \ref{r:sequence},
the desired result follows. \end{pf}

With the above preliminaries in place, we are now
ready to address the main part of the proof of Theorem \ref{thm:HTL}.
In the following, we shall reuse some notation for
local proof purposes.
In particular, $\xi$ and $\zeta$ are used for different purposes here than
earlier in the paper.

Recall, from (\ref{W}), that for each $r>0$ and $t\geq0$,
%
\begin{equation} \label{rW}
\hat W^r(t) =
\hat W^r(0)+\hat X^r(t)+ \hat U^r(t)
\end{equation}
and that, by (\ref{eqn:WN}),
%
\begin{eqnarray} \label{rWp}
\hat W^r(t) =  A(M^r)^{-1} \hat N^r (t)
= \tilde W^r(t) +\hat\xi^r(t) ,\vadjust{\goodbreak}
\end{eqnarray}
where
%
\begin{eqnarray}\label{rWtilde}
\tilde W^r(t) &=& AM^{-1} \Delta(\hat W^r(t)),
\\\label{xidef}
\hat\xi^r(t) &=& AM^{-1} \bigl(\hat N^r(t) -\Delta(\hat W^r(t))\bigr) +
A\bigl((M^r)^{-1} - M^{-1}\bigr) \hat N^r(t).
\end{eqnarray}
By Proposition \ref{Deltaw} and the definition of $\mathcal W_1$ from
$\mathcal M_1$,
we have that
%
\begin{equation}
\tilde W^r (t) \in\mathcal W_1\qquad  \mbox{for all } t\geq0  \mbox{
and for all } r>0.
\end{equation}
Returning to (\ref{rW}), for $\delta>0$ fixed and each $r>0$, $t\geq
0, $ $j\in\mathbb{J}$,
%
\begin{equation}\label{zeta}
\hat U^r_j (t) = \int_0^t 1_{\{ d(\tilde W^r(s), \mathcal{W}_1^j)\leq
\delta
\}}  \, d\hat U^r_j(s)
+\hat\zeta^{r,\delta}_j(t),
\end{equation}
where
%
\begin{equation}\label{zetatwo}
\hat\zeta^{r,\delta}_j(t) =\int_0^t 1_{\{ d(\tilde W^r(s), \mathcal{W}
_1^j)>\delta\}} \, d\hat U^r_j(s).
\end{equation}

A main
step in establishing the diffusion approximation result, Theorem~\ref
{thm:HTL}, is to show
that for
each $j\in{\mathbb J}$, with probability tending to one as
$r\to\infty$, the process $\hat U^r$ increases only when $\hat W^r$
(or $\tilde W^r$)
is near the boundary portion~$\mathcal W^j_1$. The local traffic
condition (Assumption \ref{localtraffic}) is used in verifying this
in
the following lemma, which also
shows that
state space collapse holds and that for fixed $T>0$, with high probability,
we can obtain a uniform bound on $\|\hat W^r(\cdot)\|_T$ and on $\|
\hat U^r(\cdot)\|_T$.

\begin{lemma} \label{thm:ctrlWN} Suppose, in addition to the assumption
that $\alpha=1$ and Assumption~\ref{localtraffic}, we have that
%
\begin{equation} |\hat N^r(0)
-\Delta(\hat W^r(0))| \to0
\end{equation}
in probability as $r\to\infty$. Then, for each $T>0$, $\delta>0$,
there exist
constants $K(T, \delta)>0$ and $r(T,\delta)>0$ such that for each
$r\geq r(T,\delta)$,
\begin{eqnarray}\label{est:NW}
\quad && P\bigl(
\|\hat{N}^r(\cdot)-\Delta(\hat{W}^r(\cdot))\|_T
\leq\delta,
\|\hat\xi^r(\cdot)\|_T \leq
\delta,  \nonumber
\\[-8pt]\\[-8pt]\
&&\quad\, \|\hat\zeta^{r,\delta}(\cdot)\|_T =0, \|\hat W^r(\cdot)\|_T \leq K(T,\delta
),  \|\hat U^r(\cdot)\|_T\leq K(T,\delta)
\bigr)\geq1-\delta.\nonumber
\end{eqnarray}
\end{lemma}

A main aspect of the proof of
Lemma \ref{thm:ctrlWN} is to show that for fixed $T>0$, with high probability,
we can obtain a suitable uniform bound on
$\|\hat N^r (\cdot) \|_T$ [or equivalently, on $\|\hat W^r(\cdot)\|
_T$] for all $r$ sufficiently large.
We shall use an oscillation
inequality for the proof of this. A local version of this inequality
was established in Theorem 4.1 of Kang and Williams
\cite{KaWi}. For $\alpha\not=1,$ one would need to use
that local version. However, for the case $\alpha=1$ treated here,
one can choose $\rho$ in Theorem 4.1 of \cite{KaWi} to be
arbitrarily large and consequently obtain a global version
of the oscillation inequality.
Here, for the case of $\alpha=1$ only, rather than using Theorem 4.1
of \cite{KaWi},
we shall give an alternative proof of a global oscillation inequality
by first invoking a~linear transformation to transform $\mathcal{W}_1$
to the
orthant and then applying
an oscillation inequality developed earlier
by Williams \cite{Wi98} for that state space. This inequality has a
slightly simpler
form than that which would follow from \cite{KaWi}.

For the statement of the oscillation inequality, we need the following notation.
For any $0\leq s< t<\infty$, let $\mathbb{D}([s,t], \mathbb
{R}^{{\mathbf J}})$
denote the set of functions $x\dvtx [s, t]\to\mathbb{R}^{{\mathbf J}}$
that are
right-continuous on
$[s, t)$ and have finite left limits on~$(s, t]$, and for $x\in\mathbb{D}
([s,t], \mathbb{R}^{{\mathbf J}})$,
let
%
\begin{equation}
\operatorname{Osc}(x, [s,t]) = \sup\{|x(v)-x(u)|\dvtx  s\leq u<v \leq t\}.
\end{equation}

\begin{proposition}[(Oscillation inequality)] \label{thm:OI}
$\!\!\!$There exists a constant \mbox{$\hat C_0>0$} such that
for any $\delta\geq0$ and
any $0\leq s<t<\infty$, $w, x, y \in\mathbb{D}([s,t],\mathbb
{R}^{\bf J})$
satisfying
\begin{enumerate}[(iii)]
\item[(i)] $ w(u)=x(u)+\sum_{j\in
\mathbb{J}}y_j(u)\gamma^j \mbox{ for all }u\in
[s,t]$,
\item[(ii)]
$w(u)\in\mathcal{W}_1$ for all $u\in[s,t]$,
\item[(iii)] for each $j\in\mathbb{J}$,
\begin{longlist}
\item[(a)] $y_j(s)\geq0;$
\item[(b)] $y_j \mbox{ is
nondecreasing;}$
\item[(c)]
$y_j(u)=y_j(s)+\int_{(s,u]}1_{\{d(w(v),  \mathcal
W_{1}^j)\leq\delta\}}\,dy_j(v) \mbox{ for all }u\in[s,t],$
\end{longlist}
\end{enumerate}
the following
hold:
%
\begin{eqnarray}\label{ineq1}
\operatorname{Osc}(w,[s,t])&\leq& \hat C_0\bigl(\operatorname
{Osc}(x,[s,t])+\delta\bigr),
 \\\label{ineq2}
\operatorname{Osc}(y,[s,t])&\leq& \hat C_0\bigl(\operatorname
{Osc}(x,[s,t])+\delta\bigr).
\end{eqnarray}

\end{proposition}

\begin{pf}
For $w, x, y$ satisfying (i)--(iii) above, let
$\tilde w = (ABA')^{-1} w,$ $ \tilde x=(ABA')^{-1}x $, $\tilde y = y$
and
$ R= (ABA')^{-1}$.
Since $ABA'$ is a bijection from $\mathbb{R}_+^{{\mathbf J}}$ onto
$\mathcal{W}_1$,
the matrix with columns given by the vectors $\{\gamma^j, j\in\mathbb
{J}\}$
is the ${\mathbf J}\times{\mathbf J}$ identity matrix and for
$z\in\mathcal{W}_1$, the distance $d(z, \mathcal{W}_1^j) = n^j\cdot
z/|n^j|$ for
each $j\in\mathbb{J}$ [where $n^j$ is the vector given by the
$j${th} row
of $(ABA')^{-1}$],
it follows that
$\tilde w, \tilde x, \tilde y$ satisfy
\begin{enumerate}[(iii)]
\item[(i)] $ \tilde w(u)=\tilde x(u)+ R\tilde y(u)$
for all $u\in[s,t]$,\vspace{2pt}
\item[(ii)]
$\tilde w(u)\in\mathbb{R}_+^{{\mathbf J}}$ for all $u\in[s,t]$,
\item[(iii)] for each $j\in\mathbb{J}$,
\begin{longlist}
\item[(a)] $\tilde y_j(s)\geq0;$
\item[(b)] $\tilde y_j \mbox{ is
nondecreasing;}$
\item[(c)]
$\tilde y_j(u)=\tilde y_j(s)+\int_{(s,u]}1_{\{\tilde w_j(v)\leq\delta|n^j|
\}}\,d\tilde y_j(v) \mbox{ for all }u\in[s,t].$
\end{longlist}
\end{enumerate}
Since $R$ is completely-$S$, by Lemma \ref{lem:coms}, and $\tilde w,
\tilde x ,\tilde y$
satisfy (i)--(iii) above,
it is immediate from Theorem 5.1 of \cite{Wi98} that there is a constant
$ c_1>0$ depending only on $R$ such that
%
\begin{eqnarray}
\operatorname{Osc}(\tilde w,[s,t])&\leq& c_1\Bigl(\operatorname
{Osc}(\tilde x,[s,t])+\delta\max
_{j\in\mathbb{J}}|n^j|\Bigr),
\\
\operatorname{Osc}(\tilde y,[s,t])&\leq& c_1\Bigl(\operatorname
{Osc}(\tilde x,[s,t])+\delta\max
_{j\in\mathbb{J}}|n^j|\Bigr).
\end{eqnarray}
Applying the reverse linear transformation $ABA'$ and making $\hat C_0$
sufficiently
large to absorb the factor $\max_{j\in\mathbb{J}}|n^j|$, it follows that
(\ref{ineq1})--(\ref{ineq2}) hold where the constant $\hat C_0$ can
be chosen to depend
only on $ABA'$ (and its inverse).
\end{pf}

\begin{pf*}{Proof of Lemma \protect\ref{thm:ctrlWN}}
Fix $T>0$ and $\delta>0$.
By the convergence assumed for the initial random variables $\{\hat
W^r(0), r>0\}$,
the $C$-tightness of $\{\hat{W}^r(0)+\hat{X}^r(\cdot),r>0\}$ established
in Lemma \ref{lem:hatXr}, the multiplicative state space collapse established
in Theorem \ref{req:SSC} and the fact that $(M^r)^{-1}\to M^{-1}$ as
$r\to\infty$, we have that
for each $\varepsilon >0$, there are
constants $K_0\geq1$ (not depending on $\varepsilon $) and $r_0(
\varepsilon
)>0$ such
that for all $r\geq r_0( \varepsilon )$,
%
\begin{eqnarray}\label{est:BDX}
{P}\bigl( |\hat W^r(0)| \leq K_0,
\|\hat{X}^r(\cdot)\|_T\leq K_0\bigr)&\geq& 1-\frac
{\delta}{2},
\\\label{rmss}
P\bigl(\|\hat N^r(\cdot)- \Delta(\hat W^r(\cdot))\|_T\leq
\varepsilon
\bigl(\|\hat N^r(\cdot)\|_T\vee1\bigr)
\bigr)&\geq& 1-\frac{\delta}{2}
\end{eqnarray}
and
%
\begin{equation}\label{mtom}
\|(M^r)^{-1} - M^{-1} \| \leq\varepsilon .
\end{equation}
The constants $K_0$ and $r_0(\varepsilon )$ will depend on $T$ and
$\delta
$ as well,
but since these parameters are fixed throughout this proof,
we do not explicitly indicate that dependence here.

In the following, $\varepsilon >0$ will be fixed. A specific, suitably
small value of $\varepsilon $ will be chosen
later (as a function of $\delta, T$) to ensure that various
inequalities hold
[see (\ref{epsone})].
For $r>0$, let
\begin{eqnarray*}
O^{r,\varepsilon } & = & \bigl\{
|\hat W^r(0)|\leq K_0,
\|\hat{X}^r(\cdot)\|_T\leq K_0,
\\
&&\
\|\hat N^r(\cdot)- \Delta(\hat W^r(\cdot))\|_T\leq\varepsilon \bigl(\|
\hat
N^r(\cdot)\|_T\vee1\bigr) \bigr\}.
\end{eqnarray*}
From (\ref{est:BDX})--(\ref{rmss}), we have that for all $r\geq
r_0(\varepsilon )$,
%
\begin{equation}\label{oeps}
P( O^{r,\varepsilon } ) \geq1-\delta.
\end{equation}
Now, for $r\geq r_0(\varepsilon )$, on $ O^{r,\varepsilon }$, by
(\ref{rWp})--(\ref{xidef}) and
(\ref{mtom}), we have
%
\begin{equation}\label{estxi}
\| \hat\xi^r (\cdot)\|_T \leq\hat C_3\varepsilon \|\bigl(\hat N^r(\cdot
)\|
_T\vee1\bigr),
\end{equation}
where $\hat C_3= \| AM^{-1} \| +\|A\|\vee1$ and
%
\begin{equation}\label{rxiexp}
\hat\xi^r (t) =\hat W^r(t) -\tilde W^r(t) \qquad  \mbox{for all }
t\geq0.
\end{equation}
Since, by Proposition \ref{Deltaw},
$\Delta(\hat W^r(t, \omega))$ is an invariant state,
it follows from the characterization
of invariant states in Theorem \ref{invariant}
that
for $r\geq r_0(\varepsilon )$, $t\in[0,T]$ and $\omega\in
O^{r,\varepsilon
}$, there exists
$q^r(t,\omega)\in
\mathbb{R}^{\bf J}_+$ such that
%
\begin{equation}\label{est:2}
\Delta(\hat{W}^r(t,\omega))= \operatorname
{diag}(\rho)\operatorname{diag}(\kappa
)^{-1} A'
q^r(t,\omega)
\end{equation}
and so by the definitions of $\tilde W^r$ and the matrix $B$, we have
%
\begin{equation}\label{wq}
\tilde W^r(t,\omega) = AB A' q^r(t, \omega).
\end{equation}

We now turn to the behavior of $\hat U^r$.
By (\ref{rU}), for each $r\geq r_0(\varepsilon )$, $t\in[ 0, T]$
and $j\in\mathbb{J}$,
%
\begin{eqnarray} \label{req:U}
\hat{U}_j^r(t)&=&
r \int_0^{t}\bigl(C_j-(A\Lambda
(\hat{N}^r(s)))_j\bigr)\,ds.
\end{eqnarray}
We are interested in
where $\hat U^r_j$ increases.
Fix $r\geq r_0(\varepsilon )$, $j\in\mathbb{J}$ and $\omega\in
O^{r,\varepsilon }$.
Since the integrand in (\ref{req:U}) is
nonnegative, we concentrate on determining where the integrand
is strictly positive.
At an
instant $s\in[0,T]$ such that $C_j-(A\Lambda(\hat
N^r(s,\omega)))_j>0$, by Proposition \ref{r:Lambdapro}(iv),
there is $p^r(s,\omega)\in\mathbb{R}_+^{\bf J}$ such that $p^r_j
(s,\omega)=0$ and
%
\begin{equation}\label{req:hatN}
\hat N^r_i(s,\omega)=\Lambda_i (\hat N^r(s,\omega))
\biggl(\frac{ \sum_{k \in\mathbb{J}} p^r_k(s,\omega)
A_{ki}}{\kappa_i}
\biggr)
 \qquad \mbox{for all } i \in\mathbb{I}.
\end{equation}
By the local traffic
Assumption \ref{localtraffic}, there exists an index $i_j\in\mathbb{I}$
such that
$A_{ji_j}>0$ and $A_{ki_j}=0 $ for all $k\not=j$.
Using the fact that $p_j^r(s,\omega)=0$ in (\ref{req:hatN}), it
follows that
$\hat{N}^r_{i_j}(s,\omega)=0$.
Then, since $\omega\in O^{r,\varepsilon }$, we have
%
\begin{equation}
(\Delta(\hat W^r(s, \omega)))_{i_j} \leq\varepsilon \bigl(\|\hat
N^r(\cdot,
\omega)\|_T\vee1\bigr)
\end{equation}
and so, by (\ref{est:2}) and the local traffic condition,
%
\begin{equation}
q^r_j (s, \omega) \leq\varepsilon \bigl(\|\hat N^r(\cdot, \omega)\|
_T\vee1\bigr)\rho_{i_j}^{-1} \kappa_{i_j} .
\end{equation}
Thus, letting $n^j$ denote the vector given by the $j$th row of
$(ABA')^{-1}$ (see the proof of
Lemma \ref{dwcond}), using (\ref{wq}), we have that
\begin{eqnarray*}
d(\tilde W^r (s, \omega), \mathcal{W}_1^j) & = & n^j \cdot\tilde
W^r(s,\omega) /|n^j|
\\
&= & n^j\cdot ABA' q^r(s, \omega) /|n^j|
\\
&=& q^r_j (s,\omega) /|n^j|
\\
&\leq& \varepsilon \bigl(\|\hat N^r (\cdot, \omega)\|_T \vee1\bigr) \rho
^{-1}_{i_j}\kappa_{i_j}/|n^j|.
\end{eqnarray*}

Let $\hat C_4=\max_{j\in\mathbb{J}} \rho^{-1}_{i_j}\kappa_{i_j}/|n^j|$.
It follows from the reasoning above that for each $r\geq
r_0(\varepsilon
)$, $j\in\mathbb{J}$, on $O^{r,\varepsilon }$, we have, for all
$t\in[0, T]$,
that
%
\begin{equation}\label{usimple}
\hat U^r_j (t) = \int_0^t 1_{\{ d(\tilde W^r(s), \mathcal{W}_1^j)\leq
\hat
C_4 \varepsilon
(\| \hat N^r(\cdot)\|_T\vee1) \}} \,  d\hat U^r _j (s) ,
\end{equation}
that is,
$\hat U^r_j (\cdot)$ can increase only when $d(\tilde W^r(\cdot
),\mathcal{W}_1^j)
\leq\hat C_4 \varepsilon ( \|\hat N^r(\cdot)\|_T \vee1) $.
Furthermore,
\begin{eqnarray*}
\tilde W^r (t)& = & \hat W^r(t) -\hat\xi^r(t)
\\[2pt]
&=& \hat W^r(0) +\hat X^r(t) -\hat\xi^r(t) + \hat U^r(t) ,
\end{eqnarray*}
where
$\tilde W^r(t)\in\mathcal{W}_1$.
It then follows from the oscillation inequality in Proposition~\ref{thm:OI} that for $r\geq r_0(\varepsilon )$, on
$O^{r,\varepsilon }$,
\begin{eqnarray}\label{oscW}
\operatorname{Osc}(\tilde W^r , [0, T])& \leq& \hat C_0\bigl(
\operatorname{Osc}\bigl(\hat W^r(0)
+\hat X^r(\cdot) -\hat\xi^r(\cdot),[0, T]\bigr)  \nonumber
\\[2pt]
&& {}\qquad\quad\hspace{40pt}\, + \hat C_4\varepsilon \bigl(\|\hat N^r(\cdot)\|
_T\vee
1\bigr) \bigr)\nonumber
\\[-7pt]\\[-7pt]
&\leq& \hat C_0 \bigl( 2K_0 + (2\hat C_3+\hat C_4)\varepsilon \bigl(\|
\hat
N^r(\cdot)\|_T\vee1\bigr) \bigr)\nonumber
\\[2pt]
&\leq& \hat C_0 \bigl( 2K_0 + (2\hat C_3+\hat C_4)\varepsilon \hat
C_2\bigl(\|
\hat W^r(\cdot)\|_T\vee1\bigr) \bigr), \nonumber
\end{eqnarray}
where we have used the definition of $O^{r,\varepsilon }$ and (\ref
{estxi}) for the second
inequality
and we have used Lemma \ref{lem:WN} plus the fact that $\hat C_2\geq
1$ for the third inequality.
Similarly, we obtain an oscillation bound for $\hat U^r$ for $r\geq
r_0(\varepsilon )$ on~$O^{r,\varepsilon }$:
%
\begin{eqnarray} \label{oscU}
\operatorname{Osc}(\hat U^r , [0, T]) \leq\hat C_0 \bigl( 2K_0 +
(2\hat C_3+\hat
C_4)\varepsilon \hat C_2\bigl(\|\hat W^r(\cdot)\|_T\vee1\bigr) \bigr).
\end{eqnarray}
On combining (\ref{oscW})
with (\ref{rxiexp}), (\ref{estxi}) and Lemma \ref{lem:WN}, we have
that for $r\geq r_0(\varepsilon )$, on $O^{r,\varepsilon }$,
\begin{eqnarray*}
\|\hat W^r(\cdot)\|_T \vee1 &\leq& |\hat W^r(0)|\vee1 +
\operatorname{Osc}(\hat
W^r (\cdot), [0, T])
\\[2pt]
&\leq& |\hat W^r(0)|\vee1 + \operatorname{Osc}(\tilde W^r(\cdot),
[0, T]) +\operatorname{Osc}
(\hat\xi^r(\cdot), [0, T])
\\[2pt]
&\leq& K_0 + (\hat C_0 +1)\bigl( 2K_0 + (2\hat C_3+\hat
C_4)\varepsilon
\hat C_2\bigl(\|\hat W^r(\cdot)\|_T\vee1\bigr)\bigr).
\end{eqnarray*}

Now, choose
%
\begin{equation}\label{epsone}
\ \varepsilon = \min\biggl( \frac{1}{2(\hat C_0+1)(2\hat C_3 +\hat
C_4)\hat C_2}, \frac{\delta}{2K_0\hat C_2(3+2\hat C_0) (\hat C_3+\hat C_4)}
\biggr).
\end{equation}
Then,
from the above, we conclude that for $r\geq r_0(\varepsilon )$, on
$O^{r,\varepsilon }$,
%
\begin{equation}\label{restW}
\|\hat W^r(\cdot)\|_T\vee1 \leq2K_0( 3 +2\hat C_0)
\end{equation}
and, using the fact that $ \hat C_3 \geq1$,
we have
%
\begin{eqnarray}\label{restN}
\varepsilon \bigl(\|\hat N^r(\cdot) \|_T\vee1\bigr)&\leq&
\varepsilon \hat C_2 \bigl(\|\hat W^r(\cdot)\|_T\vee1\bigr)\nonumber
\\
&\leq& \frac{\delta}{\hat C_3 +\hat C_4 }
\\
&\leq& \delta. \nonumber
\end{eqnarray}
Further, by (\ref{oscU}), (\ref{restW}) and the definition of
$\varepsilon $, we have
%
\begin{equation}\label{restU}
\operatorname{Osc}(\hat U^r , [0, T]) \leq2 \hat C_0 ( K_0 + \delta) .
\end{equation}
Note that since $\hat U^r$ is nondecreasing and starts from zero,
the above also provides a bound for $\|\hat U^r(\cdot)\|_T$.
Now, by (\ref{restN}), for $r\geq r_0(\varepsilon )$, on $O^{r,
\varepsilon
}$, we have
\begin{eqnarray*}
\|\hat N^r (\cdot) -\Delta(\hat W^r(\cdot))\|_T &\leq& \varepsilon
\bigl(\|
\hat N^r(\cdot) \|_T\vee1\bigr)
\\
&\leq& \delta,
\end{eqnarray*}
 using (\ref{estxi}), we have
\begin{eqnarray*}
\|\hat\xi^r (\cdot)\|_T &\leq& \hat C_3 \varepsilon \bigl(\|\hat
N^r(\cdot
)\|_T\vee1\bigr)
\\
&\leq& \delta,
\end{eqnarray*}
using (\ref{usimple}), we have, for each
$j\in\mathbb{J}$ and $t\in[0, T]$,
%
\begin{eqnarray}
\hat U^r_j (t)
 =  \int_0^t 1_{\{ d(\tilde W^r(s,\omega), \mathcal{W}_1^j)\leq
\delta\}}  \, d\hat U^r _j (s),
\end{eqnarray}
and so by (\ref{zeta})--(\ref{zetatwo}),
%
\begin{equation}
\|\hat\zeta^{r,\delta}(\cdot)\|_T =0.
\end{equation}

Now, set $K(T,\delta) =\max(2K_0(3+2\hat C_0),2 \hat C_0 ( K_0 +
\delta) ) $ and $r(T,\delta) = r_0(\varepsilon )$, where
$\varepsilon >0$ is given by (\ref{epsone}).
Then, for all $r\geq r(T, \delta)$, by the above and~(\ref{oeps}), we have
%
\begin{eqnarray}
&&P\bigl(
\|\hat{N}^r(\cdot)-\Delta(\hat{W}^r(\cdot))\|_T
\leq\delta,
\|\hat\xi^r(\cdot)\|_T \leq
\delta,\nonumber
\\
&&\quad\ \|\hat\zeta^{r,\delta}(\cdot)\|_T =0,
  \|\hat W^r(\cdot)\|_T \leq K(T, \delta),   \|\hat U^r(\cdot)\|_T
\leq K(T, \delta)
\bigr)
\\
&&\quad\ \geq P(O^{r,\varepsilon }) \geq1-\delta.\nonumber
\end{eqnarray}
\upqed\end{pf*}

\begin{pf*}{Proof of Theorem \protect\ref{thm:HTL}}
Assume that the hypotheses of Theorem~\ref{thm:HTL} hold.
We will show that the conditions of Theorem 5.4 of \cite{KaWi}
hold, from which it will follow immediately that $\hat W^r $ converges
in distribution
as \mbox{$r\to\infty$} to an SRBM with data $(\mathcal W_1, \theta, \Gamma
, \{\gamma^j:j\in{\mathbb J}\}, \eta)$.
The joint convergence of~$\hat N^r$ with $\hat W^r$ will then follow by
the state space collapse
established in Lemma~\ref{thm:ctrlWN}.

The conditions of Theorem 5.4 fall into four groups. We treat each of these
groups separately below.

First, the diffusion state space $\mathcal W_1$ must be a convex
polyhedron having
nonempty interior described
as the intersection of a minimal set of half-spaces and the directions
of reflection $\{\gamma^j\dvtx j\in\mathbb{J}\}$
must satisfy Assumption~5.1 of \cite{KaWi}.
These properties were verified in the proof of Lemma \ref{dwcond}.

Second, one must verify that Assumption 4.1 of \cite{KaWi}
holds. This condition amounts to verifying that $\hat W^r, \hat X^r,
\hat U^r$
satisfy conditions similar to (i), (ii), (iv) on $\tilde W, \tilde X,
\tilde U$
in the definition of an SRBM, except that the prelimit processes are
allowed to have
r.c.l.l.~rather than continuous paths and small perturbations in the
conditions (i), (ii) and (iv) are allowed.
Furthermore, the sequence $\{\hat W^r(0) + \hat X^r(\cdot)\}$ is
required to be $C$-tight.
The latter follows from Lemma \ref{lem:hatXr} and the former follows
from the
properties in (\ref{rW})--(\ref{zetatwo}), once we show that
$\hat\xi^r \to0$ and $\hat\zeta^{r,\delta^r}\to0$ in probability
as $r\to\infty$
for a suitable sequence $\{\delta^r\} $ satisfying $\delta^r\to0$ as
$r\to\infty$.
We verify the latter properties in the next paragraph.

Recall the constants $K(T,\delta)$ and $r(T,\delta)$ from Lemma \ref{thm:ctrlWN}.
Choose strictly increasing sequences of positive
constants $\{K_k, k\geq1\}$
and
$\{r_k,k\geq1\}$
such that for each $k$, $K_k\geq K(k, \frac{1}{k})$ and $r_k \geq r(k,
\frac{1}{k}) $,
and
$r_k\rightarrow\infty$ as $k\rightarrow\infty$.
Define~$\delta^r$ such that $\delta^r=1$ when $r\leq r_1$ and
$\delta^r=\frac{1}{k}$ when $r\in(r_k,r_{k+1}]$ for $k\geq1$. Then,
$\delta^r\rightarrow0$ as $r\rightarrow\infty$ and, by Lemma \ref
{thm:ctrlWN},
for each $k\geq1$, for $r_k < r\leq r_{k+1} $,
\begin{eqnarray} \label{est:NWp}
 && P\biggl(
\|\hat{N}^r(\cdot)-\Delta(\hat{W}^r(\cdot))\|_k
\leq\frac{1}{k},
\|\hat\xi^r(\cdot)\|_k \leq
\frac{1}{k},
\nonumber
\\[-8pt]\\[-8pt]
&&\quad\ \ \|\hat\zeta^{r,\delta_r}(\cdot)\|_k =0, \|\hat W^r(\cdot)\|_k
\leq K_k,  \|\hat U^r(\cdot)\|_k\leq K_k
\biggr)\geq1-\frac{1}{k}.\nonumber
\end{eqnarray}
It follows from this that
$\hat\xi^r\to0 $ and $\hat\zeta^{r,\delta_r} \to0 $ in
probability as $r\to\infty$.
The conditions of Assumption 4.1 in
\cite{KaWi} are then satisfied with $\mathcal W_1$ in place of
$\overline G$,
$r$ in place of $n$, $j $ in place of $i$, $\hat W^r$ in
place of $W^n$, $\tilde W^r$ in place of $\tilde W^n$, $\hat W^r(0)
+\hat X^r$
in place of $X^n$, $\hat U^r$ in place of $Y^n$,
$\tilde Y_j^n(t)=\int_0^t
1_{\{d(\tilde W^r(s), \mathcal{W}^j_1)\leq\delta\}}\, d\hat U^r_j(s)$,
$\gamma^{j,n}=\gamma^j$,
$\alpha^n=\hat\xi^r$ and $\beta^n = \hat\zeta^{r,\delta_r}$.

Third,
one must show that $\{ \hat W^r(0) +\hat X^r(\cdot)\}$ converges in
distribution as $r\to\infty$
to a Brownian motion with drift $\theta$, covariance matrix $\Gamma$
given by~(\ref{r:gamma}) and initial distribution $\eta$.
[This amounts to verifying condition $\mbox{(vi)}'$ of Theorem 4.3 in
\cite{KaWi}
(with $\theta$ in place of $\mu$ and $\eta$ in place of $\nu$
there). This condition is needed
so that weak limit points of $\{(\hat W^r, \hat X^r, \hat U^r)\}$
will satisfy property (iii)(a) and the last part of (ii) in the Definition
\ref{DefSRBM} of an SRBM.]
To verify this condition, we first show that $\{\bar{\hspace{-1pt}\bar{T}}{}^r(\cdot
)\}$ converges in distribution
to the deterministic process $\{\rho( t), t\geq0\}$, where $\rho
(t)=\rho t$.
Indeed, by the fact that
$\bar{\hspace{-1pt}\bar{T}}{}^r(\cdot)$ is uniformly Lipschitz continuous with
a Lipschitz constant which does not depend on $r$, we have that $\{\bar
{\hspace{-1pt}\bar{T}}{}^r(\cdot)
,r>0\}$ is $C$-tight. Let $T^*(\cdot)$ be a weak limit point of this sequence,
obtained as a limit
in distribution along a suitable subsequence.
We will prove below that, almost surely,
%
\begin{equation}\label{Trho} T^*(t)=
\rho t \qquad  \mbox{for all } t\geq0,
\end{equation}
from which the desired convergence in distribution of $\{\bar{\hspace{-1pt}\bar
T}{}^r(\cdot)\}$ to
$T^*(\cdot)$ follows.

For the proof of (\ref{Trho}), by passing to a subsequence,
we may assume that~$\{\bar{\hspace{-1pt}\bar{T}}{}^r(\cdot)\}$
converges in distribution to $T^*(\cdot)$ as $r\to\infty$.
Now, by (\ref{eqn:N}) and~(\ref{N}), on dividing by $r$,
we have, for each $i\in\mathbb{I}$,
%
\begin{eqnarray} \label{hatNrr}
\frac{\hat N_i^r (t)}{r}
&=& \frac{\hat N_i^r(0)}{r} + \bar{\hspace{-2pt}\bar E}{}^r_{\!i}(t)
-\bar{\hspace{-1pt}\bar S}{}^{r}_{\!i}(\bar{\hspace{-1pt}\bar T}{}^r_{\!i}(t)),
\end{eqnarray}
where
%
\begin{equation}
\bar{\hspace{-2pt}\bar E}{}^r(t) = \frac{E^r(r^2t)}{r^2}, \qquad
\bar{\hspace{-1pt}\bar S}{}^r(t) = \frac{S^r(r^2 t)}{r^2}.
\end{equation}
From (\ref{est:NWp}), we can conclude, using the continuity of $\Delta
(\cdot)$, that
\[
\frac{\hat N^r(\cdot)}{r}\Rightarrow0 \qquad
\mbox{as }  r\to\infty
\]
and,
from (\ref{clt}), it follows that
\[
(\bar{\hspace{-2pt}\bar E}{}^r(\cdot), \bar{\hspace{-1pt}\bar S}{}^r(\cdot)) \Rightarrow(\nu
(\cdot), \mu(\cdot))
 \qquad \mbox{as } r\to\infty,
\]
where $\nu(t) =\nu t$ and $\mu(t) =\mu t$ for all $t\geq0$.
Thus, on letting $r\to\infty$ in (\ref{hatNrr}), we obtain that,
almost surely, for
each $i\in{\mathbb I}$,
\[
0 = \nu_it - \mu_i T_i^*(t) \qquad  \mbox{for all } t\geq0.
\]
The desired result (\ref{Trho}) follows immediately.

It now follows from (\ref{clt}) and the assumption on the convergence in
distribution of $\hat W^r(0)$ made at the end of Section \ref{r:sequence}
that
$(\hat W^r(0), \hat E^r(\cdot), \hat S^r(\cdot), \bar{\hspace{-1pt}\bar
T}{}^r(\cdot))$
converges in distribution to $(\tilde W(0), \tilde E(\cdot), \tilde
S(\cdot), \rho(\cdot))$
as $r\to\infty$, where $\tilde W(0)$ is independent of the Brownian
motion $(\tilde E(\cdot), \tilde S(\cdot))$
and $\tilde W(0)$ has distribution $\eta$ on $\mathcal W_1$.
On combining this with a random time change theorem and Assumption \ref{HT},
we conclude, using (\ref{rXhat}), that
$\{\hat W^r(0) +\hat X^r(\cdot), r>0\}$ converges in distribution
to a Brownian motion with drift $\theta$, covariance matrix~$\Gamma$
given by (\ref{r:gamma}) and initial distribution $\eta$, as desired.

Fourth, and finally,
we must verify condition (vii) of Theorem 4.3 of \cite{KaWi} (with
$\theta$ in place of $\mu$ there).
This condition requires that for any weak limit point $(W, X, U)$ of
$(\hat W^r, \hat X^r, \hat U^r)$,
$\{X(t)-\theta t , t\geq0\}$ is a martingale relative to the
filtration generated
by $(W, X, U)$.
[This condition is needed so that property (iii)(b) in Definition \ref
{DefSRBM} of an SRBM
will be satisfied by weak limit points of $\{(\hat W^r,\hat X^r, \hat
U^r)\}$.]
Let $\hat\theta^r =r (A\rho^r -C)$.
Then, by Proposition~4.1 of~\cite{KaWi},
to show that this final condition holds, it suffices to verify the
following properties
for the prelimit processes:
\begin{longlist}
\item[(a)] for each $r,$ $\{\hat X^r(t) -\hat\theta^r t, t\geq0\}$
is a martingale with respect to the
filtration generated by $(\hat W^r, \hat X^r, \hat U^r)$;
\item[(b)] $\hat\theta^r\to\theta$ as $r\to\infty$;
\item[(c)] $\{\hat X^r (t) , r>0\} $ is uniformly integrable for each
fixed $t\geq0$.
\end{longlist}
We now verify these three properties.

It is well known (see Theorem 4.1, Chapter 6 of \cite{EK}) that for
the continuous-time Markov chain
$N^r$, for each $r>0$ and $i\in\mathbb{I}$,
\[
\chi_i ^r (t) = N_i^r(t) -N_i^r(0) -\int_0^t \bigl(\nu_i^r -\mu_i^r
\Lambda_i(N^r(s))\bigr) \, ds, \qquad t\geq0,
\]
defines a martingale with respect to the filtration generated by $N^r$,
that is, with respect to $\{\mathcal F_t^r, t\geq0\}$,
where $\mathcal F_t^r=\sigma\{ N^r(s)\dvtx 0\leq s\leq t\}$.
Setting
\[
\hat\chi^r(t) =\frac{\chi^r(r^2t)}{r}, \qquad \hat{\mathcal
F}_t^r =\mathcal F^r_{r^2t},
\]
we have that for each $r>0$,
$\{\hat\chi^r(t), \hat{\mathcal F}^r_t,t\geq0 \}$
is a multidimensional martingale.
Now, using the expressions (\ref{eqn:WN}) and (\ref{W}) for $\hat
W^r$, and the
expression~(\ref{rU}) for $\hat U^r$, we see that
\begin{eqnarray*}
A(M^r)^{-1} \hat\chi^r(t) & = & \hat W^r(t) -\hat W^r(0) - A\rho^r rt
+Crt -\hat U^r(t)
\\
&=& \hat X^r(t) + r(C-A\rho^r ) t .
\end{eqnarray*}
Since $\hat W^r$ and $\hat U^r$ are adapted
to the filtration $\{\hat{\mathcal F}_t^r, t\geq0\}$, so is $\hat X^r$,
and it then follows that
\[
\hat X^r (t) -\hat\theta^r t = A(M^r)^{-1} \hat\chi^r (t) ,  \qquad
t\geq0,
\]
is a martingale with respect to the filtration generated by $(\hat W^r,
\hat X^r, \hat U^r)$.
By (\ref{req:Jplus}), $\hat\theta^r \to\theta$ as $r\to\infty$.
Finally, since the mean of a Poisson random variable
is the same as its variance and the rate of increase
of $\bar{\hspace{-1pt}\bar{T}}{}^r$ is uniformly bounded, it is
straightforward to verify, using Assumption \ref{HT},
that $\{\hat X^r (t) , r>0\}$ is uniformly integrable for each $t\geq0$.
Thus, we have verified conditions (a)--(c) above and
so condition (vii) of Theorem 4.3 in \cite{KaWi} holds.\looseness=1

We have now verified all of the hypotheses of
Theorem 5.4 of \cite{KaWi} and so it follows that
$\hat W^r$ converges in distribution as $r\to\infty$ to an SRBM
$\tilde W$ associated with the data
$(\mathcal W_1, \theta, \Gamma, \{\gamma^j\dvtx j\in{\mathbb J}\}, \eta)$.
From (\ref{est:NWp}), we have that $\hat N^r -\Delta(\hat W^r)$
converges in
distribution to the zero process as $r\to\infty$.
Since $\Delta(\cdot)$ is continuous, it follows that
we have the joint convergence of $(\hat W^r, \hat N^r)$ in distribution
as $r\to\infty$ to $(\tilde W, \Delta(\tilde W))$.\vadjust{\goodbreak}
\end{pf*}

\section{Conclusion}

A bandwidth sharing policy corresponds to a generalization of the
notion of a processor sharing discipline from a single resource to
a~network with several shared resources. In particular, weighted
$\alpha$-fair policies provide a tractable theoretical abstraction
of the bandwidth sharing effected by decentralized packet-based
end-to-end congestion control algorithms such as TCP. It is known
\cite{BM,DLK} that, for flow-level models with exponentially
distributed file sizes, weighted $\alpha$-fair policies are stable
when the average load on each resource is less than its capacity.

Weighted $\alpha$-fair policies can nevertheless suffer from
entrainment of resources, whereby congestion at some resources may
prevent others from working at full capacity: this is manifested
under diffusion scaling, where a~Brownian model for the workload
process lives in a cone which may be a~strict subset of the positive
orthant.

Under weighted proportional fair sharing ($\alpha=1$) and a mild
local traffic condition, this paper has shown how multiplicative
state space collapse can be combined with an invariance principle
to establish a diffusion approximation for the flow-level model. For
proportional fair sharing (equal weights), this diffusion has a
product form invariant measure which, when integrable, can be
normalized to yield
the unique stationary distribution for the diffusion. This result
extends to the case where file sizes are distributed
as finite mixtures of exponentials. Thus, in the diffusion limit,
more networks can have product form stationary distributions than are
known in the prelimit (see Bonald and Prouti\`{e}re \cite{BP}).
We also indicated some open problems for $\alpha\neq1$ and an
extension to models with routing.

Bandwidth sharing policies outside the class of weighted
$\alpha$-fair policies may avoid entrainment, although such policies
may not be easy to effect via decentralized packet-based end-to-end
congestion control algorithms.

%

\printaddresses


\begin{thebibliography}{44}

\bibitem{BF}
%
\begin{barticle}[auto]
\bauthor{\bsnm{Ben~Fredj},~\bfnm{S.}\binits{S.}},
\bauthor{\bsnm{Bonald},~\bfnm{T.}\binits{T.}},
\bauthor{\bsnm{Prouti\`{e}re},~\bfnm{A.}\binits{A.}},
\bauthor{\bsnm{Regnie},~\bfnm{G.}\binits{G.}} \AND
\bauthor{\bsnm{Roberts},~\bfnm{J.}\binits{J.}}
(\byear{2001}).
\btitle{Statistical bandwidth sharing: A study of congestion at flow level}.
\bjournal{Computer Communication Review}
\bvolume{31}
\bpages{111--121}.
\end{barticle}
%
\endbibitem

\bibitem{Bi99}
%
\begin{bbook}[msn]
\bauthor{\bsnm{Billingsley},~\bfnm{Patrick}\binits{P.}}
(\byear{1999}).
\btitle{Convergence of Probability Measures},
\bedition{2nd} ed.
\bpublisher{Wiley}, \baddress{New York}.
\bmrnumber{MR1700749}
\end{bbook}
%
\endbibitem

\bibitem{BM}
%
\begin{barticle}[unstr]
\bauthor{\bsnm{Bonald},~\bfnm{T.}\binits{T.}} \AND
\bauthor{\bsnm{Massouli\'{e}},~\bfnm{L.}\binits{L.}}
(\byear{2001}).
\btitle{Impact of fairness on Internet performance}.
\bjournal{ACM Sigmetrics Performance Evaluation Review}
\bvolume{29}
\bpages{82--91}.
\end{barticle}
%
\endbibitem

\bibitem{BP}
%
\begin{barticle}[msn]
\bauthor{\bsnm{Bonald},~\bfnm{T.}\binits{T.}} \AND
\bauthor{\bsnm{Prouti{\`e}re},~\bfnm{A.}\binits{A.}}
(\byear{2003}).
\btitle{Insensitive bandwidth sharing in data networks}.
\btitle{Queueing Syst.}
\bvolume{44}
\bpages{69--100}.
\bmrnumber{MR1989867}
\end{barticle}
%
\endbibitem

\bibitem{BR}
%
\begin{barticle}[msn]
\bauthor{\bsnm{Bramson},~\bfnm{Maury}\binits{M.}}
(\byear{1998}).
\btitle{State space collapse with application to heavy traffic limits for
multiclass queueing networks}.
\bjournal{Queueing Syst.}
\bvolume{30}
\bpages{89--148}.
\bmrnumber{MR1663763}
\end{barticle}
%
\endbibitem

\bibitem{Bra2005}
%
\begin{bmisc}[unstr]
\bauthor{\bsnm{Bramson},~\bfnm{M.}\binits{M.}}
(2009).
Network stability under max--min fair bandwidth sharing. Preprint.
\end{bmisc}
%
\endbibitem

\bibitem{BUT}
%
\begin{barticle}[unstr]
\bauthor{\bsnm{Bu},~\bfnm{T.}\binits{T.}} \AND
\bauthor{\bsnm{Towsley},~\bfnm{D.}\binits{D.}}
(\byear{2001}).
\btitle{Fixed point approximation for TCP behavior in
an AQM network}.
\bjournal{ACM Sigmetrics Performance Ecaluation Review}
\bvolume{29}
\bpages{216--225}.
\end{barticle}
%
\endbibitem

\bibitem{CST2006}
%
\begin{bmisc}[unstr]
\bauthor{\bsnm{Chiang},~\bfnm{M.}\binits{M.}},
\bauthor{\bsnm{Shah},~\bfnm{D.}\binits{D.}} \AND
\bauthor{\bsnm{Tang},~\bfnm{A.~K.}\binits{A.~K.}}
(2006). Stochastic stability under network utility maximization:
General file size distribution. In \textit{Proceedings of 44th Allerton
Conference on Communication, Control and Computing},
September 2006, Monticello, IL.
\end{bmisc}
%
\endbibitem

\bibitem{CPS}
%
\begin{bbook}[msn]
\bauthor{\bsnm{Cottle},~\bfnm{Richard~W.}\binits{R.~W.}},
\bauthor{\bsnm{Pang},~\bfnm{Jong-Shi}\binits{J.-S.}} \AND
\bauthor{\bsnm{Stone},~\bfnm{Richard~E.}\binits{R.~E.}}
(\byear{1992}).
\btitle{The Linear Complementarity Problem}.
\bpublisher{Academic Press}, \baddress{Boston, MA}.
\bmrnumber{MR1150683}
\end{bbook}
%
\endbibitem

\bibitem{DW}
%
\begin{barticle}[msn]
\bauthor{\bsnm{Dai},~\bfnm{J. G.}\binits{J. G.}} \AND
\bauthor{\bsnm{Williams},~\bfnm{R. J.}\binits{R. J.}}
(\byear{1995}).
\btitle{Existence and uniqueness of
semimartingale reflecting Brownian motions in convex polyhedrons}.
\bjournal{Theory Probab. Appl.}
\bvolume{40}
\bpages{1--40}.
\bmisc{Correction:}
\bvolume{50}
(\byear{2006}),
\bpages{346--347}.
\end{barticle}
%
\endbibitem

\bibitem{DLK}
%
\begin{barticle}[auto]
\bauthor{\bsnm{De~Veciana},~\bfnm{G.}\binits{G.}},
\bauthor{\bsnm{Lee},~\bfnm{T.~J.}\binits{T.~J.}} \AND
\bauthor{\bsnm{Konstantopoulos},~\bfnm{T.}\binits{T.}}
(\byear{2001}).
\btitle{Stability and performance analysis of networks supporting elastic
services}.
\bjournal{IEEE/ACM Transactions on Networking}
\bvolume{9}
\bpages{2--14}.
\end{barticle}
%
\endbibitem

\bibitem{EK}
%
\begin{bbook}[msn]
\bauthor{\bsnm{Ethier},~\bfnm{Stewart~N.}\binits{S.~N.}} \AND
\bauthor{\bsnm{Kurtz},~\bfnm{Thomas~G.}\binits{T.~G.}}
(\byear{1986}).
\btitle{Markov Processes:
Characterization and Convergence}.
\bpublisher{Wiley}, \baddress{New York}.
\bmrnumber{MR838085}
\end{bbook}
%
\endbibitem

\bibitem{GSV}
%
\begin{bmisc}[unstr]
\bauthor{\bsnm{Gibbens},~\bfnm{R.~J.}\binits{R.~J.}},
\bauthor{\bsnm{Sargood},~\bfnm{S.~K.}\binits{S.~K.}},
\bauthor{\bsnm{Van~Eijl},~\bfnm{C.}\binits{C.}},
\bauthor{\bsnm{Kelly},~\bfnm{F.~P.}\binits{F.~P.}},
\bauthor{\bsnm{Azmoodeh},~\bfnm{H.}\binits{H.}},
\bauthor{\bsnm{Macfadyen},~\bfnm{R.~N.}\binits{R.~N.}} \AND
\bauthor{\bsnm{Macfadyen},~\bfnm{N.~W.}\binits{N.~W.}}
\bhowpublished{(2000). Fixed-point models for the end-to-end performance
analysis of IP networks. In \textit{13th ITC Specialist Seminar: IP Traffic
Measurement, Modeling and Management}, Sept 2000, Monterey, California}.
\end{bmisc}
%
\endbibitem

\bibitem{growil06}
%
\begin{barticle}[msn]
\bauthor{\bsnm{Gromoll},~\bfnm{H.~Christian}\binits{H.~C.}} \AND
\bauthor{\bsnm{Williams},~\bfnm{Ruth~J.}\binits{R.~J.}}
(\byear{2009}).
\btitle{Fluid limits for networks with bandwidth sharing and general document
size distributions}.
\bjournal{Ann. Appl. Probab.}
\bvolume{19}
\bpages{243--280}.
\bmrnumber{MR2498678}
\end{barticle}
%
\endbibitem

\bibitem{growil07}
%
\begin{bmisc}[unstr]
\bauthor{\bsnm{Gromoll},~\bfnm{H.~C.}\binits{H.~C.}} \AND
\bauthor{\bsnm{Williams},~\bfnm{R.~J.}\binits{R.~J.}}
(2009). Fluid model for a data network with $\alpha$-fair
bandwidth sharing and general document size distributions: Two examples of
stability. In \textit{Proceedings of Markov Processes and Related
Topics---A Festschrift for Thomas G. Kurtz. IMS Collections}
\textbf{4} 253--265.
Institute of Mathematical Statistics, Beachwood, OH.
\end{bmisc}
%
\endbibitem

\bibitem{HAR}
%
\begin{barticle}[msn]
\bauthor{\bsnm{Harrison},~\bfnm{J.~Michael}\binits{J.~M.}}
(\byear{2000}).
\btitle{Brownian models of open processing networks: Canonical representation
of workload}.
\bjournal{Ann. Appl. Probab.}
\bvolume{10}
\bpages{75--103}.
\bmrnumber{MR1765204}.
\bnote{Correction: \textbf{13} (2003) 390--393}.
\bmrnumber{MR1952004}
\end{barticle}
%
\endbibitem

\bibitem{HaWi87}
%
\begin{barticle}[msn]
\bauthor{\bsnm{Harrison},~\bfnm{J.~M.}\binits{J.~M.}} \AND
\bauthor{\bsnm{Williams},~\bfnm{R.~J.}\binits{R.~J.}}
(\byear{1987}).
\btitle{Multidimensional reflected {B}rownian motions having exponential
stationary distributions}.
\bjournal{Ann. Probab.}
\bvolume{15}
\bpages{115--137}.
\bmrnumber{MR877593}
\end{barticle}
%
\endbibitem

\bibitem{HaWi87b}
%
\begin{barticle}[msn]
\bauthor{\bsnm{Harrison},~\bfnm{J.~M.}\binits{J.~M.}} \AND
\bauthor{\bsnm{Williams},~\bfnm{R.~J.}\binits{R.~J.}}
(\byear{1987}).
\btitle{Brownian models of open queueing networks with homogeneous customer
populations}.
\bjournal{Stochastics}
\bvolume{22}
\bpages{77--115}.
\bmrnumber{MR912049}
\end{barticle}
%
\endbibitem

\bibitem{HSHST}
%
\begin{barticle}[auto]
\bauthor{\bsnm{Han},~\bfnm{H.}\binits{H.}},
\bauthor{\bsnm{Shakkottai},~\bfnm{S.}\binits{S.}},
\bauthor{\bsnm{Hollot},~\bfnm{C.~V.}\binits{C.~V.}},
\bauthor{\bsnm{Srikant},~\bfnm{R.}\binits{R.}} \AND
\bauthor{\bsnm{Towsley},~\bfnm{D.}\binits{D.}}
(\byear{2006}).
\btitle{Multi-path TCP: A joint congestion control and routing scheme to
exploit path diversity on the Internet}.
\bjournal{IEEE/ACM Transactions on Networking}
\bvolume{14}
\bpages{1260--1271}.
\end{barticle}
%
\endbibitem

\bibitem{KaWi}
%
\begin{barticle}[msn]
\bauthor{\bsnm{Kang},~\bfnm{W.}\binits{W.}} \AND
\bauthor{\bsnm{Williams},~\bfnm{R.~J.}\binits{R.~J.}}
(\byear{2007}).
\btitle{An invariance principle for semimartingale reflecting {B}rownian
motions in domains with piecewise smooth boundaries}.
\bjournal{Ann. Appl. Probab.}
\bvolume{17}
\bpages{741--779}.
\bmrnumber{MR2308342}
\end{barticle}
%
\endbibitem

\bibitem{KaWiS}
%
\begin{bmisc}[unstr]
\bauthor{\bsnm{Kang},~\bfnm{W.~N.}\binits{W.~N.}} \AND
\bauthor{\bsnm{Williams},~\bfnm{R.~J.}\binits{R.~J.}}
(2009). Diffusion approximation for an input-queued packet
switch operating under a maximum weight algorithm. To appear.
\end{bmisc}
%
\endbibitem

\bibitem{KKLW}
%
\begin{bmisc}[unstr]
\bauthor{\bsnm{Kang},~\bfnm{W.}\binits{W.}},
\bauthor{\bsnm{Kelly},~\bfnm{F.~P.}\binits{F.~P.}},
\bauthor{\bsnm{Lee},~\bfnm{N.~H.}\binits{N.~H.}} \AND
\bauthor{\bsnm{Williams},~\bfnm{R.~J.}\binits{R.~J.}}
(2004). On fluid and Brownian approximations for an Internet
congestion control model. In \textit{Proceedings of the 43rd IEEE
Conference on
Decision and Control}. December 2004, 3938--3943.
\end{bmisc}
%
\endbibitem

\bibitem{KKLWM}
%
\begin{barticle}[auto]
\bauthor{\bsnm{Kang},~\bfnm{W.~N.}\binits{W.~N.}},
\bauthor{\bsnm{Kelly},~\bfnm{F.~P.}\binits{F.~P.}},
\bauthor{\bsnm{Lee},~\bfnm{N.~H.}\binits{N.~H.}} \AND
\bauthor{\bsnm{Williams},~\bfnm{R.~J.}\binits{R.~J.}}
(\byear{2007}).
\btitle{Product form stationary distributions for diffusion
approximations to
a flow-level model operating under a proportional fair sharing policy}.
\bjournal{ACM SIGMETRICS Performance Evaluation Review}
\bvolume{35}
\bpages{36--38}.
\end{barticle}
%
\endbibitem

\bibitem{KEL}
%
\begin{barticle}[msn]
\bauthor{\bsnm{Kelly},~\bfnm{F.~P.}\binits{F.~P.}}
(\byear{1991}).
\btitle{Loss networks}.
\bjournal{Ann. Appl. Probab.}
\bvolume{1}
\bpages{319--378}.
\bmrnumber{MR1111523}
\end{barticle}
%
\endbibitem

\bibitem{KV}
%
\begin{bmisc}[unstr]
\bauthor{\bsnm{Kelly},~\bfnm{F.~P.}\binits{F.~P.}} \AND
\bauthor{\bsnm{Voice},~\bfnm{T.}\binits{T.}}
(2005). Stability of end-to-end algorithms for joint
routing and
rate control. \textit{Computer Communication Review} \textbf{35} 5--12.
\end{bmisc}
%
\endbibitem

\bibitem{KW}
%
\begin{barticle}[msn]
\bauthor{\bsnm{Kelly},~\bfnm{F.~P.}\binits{F.~P.}} \AND
\bauthor{\bsnm{Williams},~\bfnm{R.~J.}\binits{R.~J.}}
(\byear{2004}).
\btitle{Fluid model for a network operating under a fair bandwidth-sharing
policy}.
\bjournal{Ann. Appl. Probab.}
\bvolume{14}
\bpages{1055--1083}.
\bmrnumber{MR2071416}
\end{barticle}
%
\endbibitem

\bibitem{KM}
%
\begin{barticle}[msn]
\bauthor{\bsnm{Key},~\bfnm{Peter}\binits{P.}} \AND
\bauthor{\bsnm{Massouli{\'e}},~\bfnm{Laurent}\binits{L.}}
(\byear{2006}).
\btitle{Fluid models of integrated traffic and multipath routing}.
\bjournal{Queueing Syst.}
\bvolume{53}
\bpages{85--98}.
\bmrnumber{MR2230015}
\end{barticle}
%
\endbibitem

\bibitem{LaBeSr2004}
%
\begin{bmisc}[unstr]
\bauthor{\bsnm{Lakshmikantha},~\bfnm{A.}\binits{A.}},
\bauthor{\bsnm{Beck},~\bfnm{C.~L.}\binits{C.~L.}} \AND
\bauthor{\bsnm{Srikant},~\bfnm{R.}\binits{R.}}
(2004). Connection level stability analysis of the Internet
using the sum of squares (SoS) techniques. In \textit{Conference on Information
Sciences and Systems, Princeton}.
\end{bmisc}
%
\endbibitem

\bibitem{LAWS}
%
\begin{barticle}[msn]
\bauthor{\bsnm{Laws},~\bfnm{C.~N.}\binits{C.~N.}}
(\byear{1992}).
\btitle{Resource pooling in queueing networks with dynamic routing}.
\bjournal{Adv. in Appl. Probab.}
\bvolume{24}
\bpages{699--726}.
\bmrnumber{MR1174386}
\end{barticle}
%
\endbibitem

\bibitem{LinSch2004}
%
\begin{bmisc}[unstr]
\bauthor{\bsnm{Lin},~\bfnm{X.}\binits{X.}} \AND
\bauthor{\bsnm{Shroff},~\bfnm{N.}\binits{N.}}
(2004). On the stability region of congestion control. In
\textit{Proceedings of the Allerton Conference on Communications,
Control and
Computing}, September 2004, Monticello, IL.
\end{bmisc}
%
\endbibitem

\bibitem{LinSchSri2006}
%
\begin{barticle}[msn]
\bauthor{\bsnm{Lin},~\bfnm{Xiaojun}\binits{X.}},
\bauthor{\bsnm{Shroff},~\bfnm{Ness~B.}\binits{N.~B.}} \AND
\bauthor{\bsnm{Srikant},~\bfnm{R.}\binits{R.}}
(\byear{2008}).
\btitle{On the connection-level stability of congestion-controlled
communication networks}.
\bjournal{IEEE Trans. Inform. Theory}
\bvolume{54}
\bpages{2317--2338}.
\bmrnumber{MR2450864}
\end{barticle}
%
\endbibitem

\bibitem{Mas2005}
%
\begin{barticle}[msn]
\bauthor{\bsnm{Massouli{\'e}},~\bfnm{Laurent}\binits{L.}}
(\byear{2007}).
\btitle{Structural properties of proportional fairness: Stability and
insensitivity}.
\bjournal{Ann. Appl. Probab.}
\bvolume{17}
\bpages{809--839}.
\bmrnumber{MR2326233}
\end{barticle}
%
\endbibitem

\bibitem{RM}
%
\begin{barticle}[auto]
\bauthor{\bsnm{Massouli\'{e}},~\bfnm{L.}\binits{L.}} \AND
\bauthor{\bsnm{Roberts},~\bfnm{J.}\binits{J.}}
(\byear{2000}).
\btitle{Bandwidth sharing and admission control for elastic traffic}.
\bjournal{Telecommunication Systems}
\bvolume{15}
\bpages{185--201}.
\end{barticle}
%
\endbibitem

\bibitem{MasRob2002}
%
\begin{barticle}[auto]
\bauthor{\bsnm{Massouli\'{e}},~\bfnm{L.}\binits{L.}} \AND
\bauthor{\bsnm{Roberts},~\bfnm{J.}\binits{J.}}
(\byear{2002}).
\btitle{Bandwidth sharing: Objectives and algorithms}.
\bjournal{IEEE/ACM Transactions on Networking}
\bvolume{10}
\bpages{320--328}.
\end{barticle}
%
\endbibitem

\bibitem{MW}
%
\begin{barticle}[auto]
\bauthor{\bsnm{Mo},~\bfnm{J.}\binits{J.}} \AND
\bauthor{\bsnm{Walrand},~\bfnm{J.}\binits{J.}}
(\byear{2000}).
\btitle{Fair end-to-end window-based congestion control}.
\bjournal{IEEE/ACM Transactions on Networking}
\bvolume{8}
\bpages{556--567}.
\end{barticle}
%
\endbibitem

\bibitem{REV}
%
\begin{bmisc}[unstr]
\bauthor{\bsnm{Roughan},~\bfnm{M.}\binits{M.}},
\bauthor{\bsnm{Erramilli},~\bfnm{A.}\binits{A.}} \AND
\bauthor{\bsnm{Veitch},~\bfnm{D.}\binits{D.}}
(2001). Network performance for TCP networks, Part I: Persistent
sources. In \textit{Proceeding ITC'17}, Brasil, September 2001.
\end{bmisc}
%
\endbibitem

\bibitem{SCH}
%
\begin{barticle}[msn]
\bauthor{\bsnm{Schassberger},~\bfnm{R.}\binits{R.}}
(\byear{1986}).
\btitle{Two remarks on insensitive stochastic models}.
\bjournal{Adv. in Appl. Probab.}
\bvolume{18}
\bpages{791--814}.
\bmrnumber{MR857330}
\end{barticle}
%
\endbibitem

\bibitem{SW}
%
\begin{bmisc}[unstr]
\bauthor{\bsnm{Shah},~\bfnm{D.}\binits{D.}} \AND
\bauthor{\bsnm{Wischik},~\bfnm{D.}\binits{D.}}
\bhowpublished{(2007). Heavy traffic analysis of optimal scheduling algorithms
for switched networks. Preprint}.
\end{bmisc}
%
\endbibitem

\bibitem{Sri2004}
%
\begin{bmisc}[unstr]
\bauthor{\bsnm{Srikant},~\bfnm{R.}\binits{R.}}
\bhowpublished{(2004). On the positive recurrence of a Markov chain describing
file arrivals and departures in a congestion-controlled network. In
\textit{Presented at the IEEE Computer Communications Workshop}}.
\end{bmisc}
%
\endbibitem

\bibitem{TE}
%
\begin{barticle}[msn]
\bauthor{\bsnm{Tassiulas},~\bfnm{Leandros}\binits{L.}} \AND
\bauthor{\bsnm{Ephremides},~\bfnm{Anthony}\binits{A.}}
(\byear{1992}).
\btitle{Stability properties of constrained queueing systems and scheduling
policies for maximum throughput in multihop radio networks}.
\bjournal{IEEE Trans. Automat. Control}
\bvolume{37}
\bpages{1936--1948}.
\bmrnumber{MR1200609}
\end{barticle}
%
\endbibitem

\bibitem{RSS}
%
\begin{barticle}[msn]
\bauthor{\bsnm{Williams},~\bfnm{R.~J.}\binits{R.~J.}}
(\byear{1987}).
\btitle{Reflected {B}rownian motion with skew symmetric data in a~polyhedral
domain}.
\bjournal{Probab. Theory Related Fields}
\bvolume{75}
\bpages{459--485}.
\bmrnumber{MR894900}
\end{barticle}
%
\endbibitem

\bibitem{Wi}
%
\begin{barticle}[msn]
\bauthor{\bsnm{Williams},~\bfnm{R.~J.}\binits{R.~J.}}
(\byear{1998}).
\btitle{Diffusion approximations for open multiclass queueing networks:
Sufficient conditions involving state space collapse}.
\bjournal{Queueing Systems Theory Appl.}
\bvolume{30}
\bpages{27--88}.
\bmrnumber{MR1663759}
\end{barticle}
%
\endbibitem

\bibitem{Wi98}
%
\begin{barticle}[msn]
\bauthor{\bsnm{Williams},~\bfnm{R.~J.}\binits{R.~J.}}
(\byear{1998}).
\btitle{An invariance principle for semimartingale reflecting {B}rownian
motions in an orthant}.
\bjournal{Queueing Systems Theory Appl.}
\bvolume{30}
\bpages{5--25}.
\bmrnumber{MR1663755}
\end{barticle}
%
\endbibitem

\bibitem{YY}
%
\begin{bmisc}[unstr]
\bauthor{\bsnm{Ye},~\bfnm{H.}\binits{H.}} \AND
\bauthor{\bsnm{Yao},~\bfnm{D.~D.}\binits{D.~D.}}
(2008). Heavy-traffic optimality of a stochastic
network under
utility-maximizing resource allocation. \textit{Operations Research}. \textbf{56}
453--470.
\end{bmisc}
%
\endbibitem

\end{thebibliography}
\end{document}